\documentclass[11pt]{article}
\usepackage[letterpaper,margin=1in]{geometry}
\usepackage{amsmath,amssymb,amsthm,bm,mathtools}
\usepackage[T1]{fontenc}
\usepackage{lmodern,microtype}
\usepackage{xcolor}
\usepackage{graphicx}
\graphicspath{{figures/numerical/}{latex/figures/numerical/}}
\usepackage[numbers,sort&compress]{natbib}
\usepackage[colorlinks=true,linkcolor=blue,urlcolor=blue]{hyperref}
\usepackage{algorithm}
\usepackage{algpseudocode}
\numberwithin{equation}{section}
\allowdisplaybreaks[2]
\newtheorem{remark}{Remark}[section]
\newtheorem{corollary}{Corollary}[section]
\newtheorem{theorem}{Theorem}[section]
\title{An Explicit Parallelizable Robin--Robin Splitting Scheme\\
for Navier--Stokes--Naghdi Shell Interaction}
\author{{Yifan Wang}\\
\small Department of Mathematics and Statistics\\
\small Texas Tech University\\
\small Lubbock, TX, USA 79401\\
}
\date{}
\begin{document}
\maketitle
{

\begin{abstract}
{
We develop a novel explicit Robin-Robin splitting scheme for the system of an incompressible
viscous fluid coupled to a time dependent Naghdi shell model. The shell model includes three
displacement components and three independent rotation components, accounting
for membrane, bending, and transverse-shear deformations. Both Robin conditions
depend only on interface data from the previous time step, allowing the fluid
and shell subproblems to be solved independently without coupling subiterations.
For the fixed-domain Stokes-Naghdi problem, we derive a discrete energy identity
that incorporates shell rotational inertia, interface velocities, and interface
loads. We prove that the scheme is unconditionally stable for any fixed positive
common Robin coefficient under homogeneous external data. Under suitable
regularity and initialization assumptions, we establish first-order consistency
and an $O(\sqrt{\Delta t})$ temporal error bound in the mechanical energy norm.
Numerical experiments with a moving-domain pressure pulse compare the proposed
scheme with the kinematically coupled $\beta$-scheme and a monolithic solver.
To demonstrate its applicability to a more realistic vascular configuration,
we simulate a stented artery using a bonded layered shell with stiff strut
bands, capturing spatial variations in the arterial wall response.
}
\end{abstract}



\section{Introduction}\label{sec:introduction}

Fluid-structure interaction between an incompressible fluid and a thin elastic
wall plays an important role in modeling flow through compliant vessels.
In this paper, we consider a viscous incompressible fluid coupled to a
Naghdi-type shell with three Cartesian displacement components and three
independent infinitesimal rotation components. Unlike the conventional
Koiter shell model, the Naghdi shell accounts for membrane, transverse-shear,
and bending deformations using only first spatial derivatives of the unknowns.
The fluid-Naghdi shell model builds on our earlier work
~\cite{TambacaTutek2016, WangTambacaCanic2026}, where we developed a sequential,
non-iterative loosely coupled scheme that incorporates shell inertia through
a Robin boundary condition for the fluid subproblem. In the present study,
we retain the same Naghdi shell formulation but introduce a
different interface partitioning strategy that allows the fluid and shell
subproblems to be solved concurrently at each time step. Our main focus is
the development and analysis of this partitioned scheme on a fixed domain.

Fluid-structure interaction problems are commonly solved using either monolithic or partitioned methods. Monolithic methods solve the fluid and structural unknowns simultaneously, while partitioned methods treat the two subproblems separately, allowing existing fluid and structural solvers to be used. However, partitioned methods can suffer from the added-mass effect, particularly when an incompressible fluid interacts with a thin, light elastic wall. As shown by Causin, Gerbeau, and Nobile~\cite{CausinGerbeauNobile2005}, conventional explicit Dirichlet-Neumann coupling can become unstable in this setting. Strongly coupled partitioned methods can overcome this difficulty, but they generally require multiple subiterations within each time step, increasing the computational cost. This motivates the development of stable, loosely coupled partitioned schemes that advance the fluid and structural subproblems separately without requiring subiterations. Several approaches have been developed to improve the stability of partitioned schemes through the choice of interface conditions. The kinematically coupled $\beta$-scheme provides a non-iterative method for fluid--thin-wall interaction, including models that account for longitudinal wall displacement~\cite{BukacEtAl2013}. For the case $\beta=1$, Buka\v{c} and Muha~\cite{BukacMuha2016} proved unconditional stability and first-order convergence in time for a linear model. Other approaches include stabilized explicit coupling using Nitsche's method~\cite{BurmanFernandez2009} and generalized Robin--Neumann schemes for thick elastic structures~\cite{FernandezMullaertVidrascu2015}. Robin-Robin coupling has also been analyzed, with stability and error estimates established at the time-semidiscrete level~\cite{BurmanDurstGuzman2022} and for fully discrete schemes~\cite{BurmanEtAl2022}. These studies demonstrate how the choice of interface conditions affects the stability and accuracy of partitioned methods.

In this work, we extend the parallelizable Robin-Robin splitting scheme developed by \cite{GuoEtAl2025RobinRobin} for an incompressible fluid interacting with a thick elastic structure to a reduced-dimensional Naghdi shell model. The key feature of the proposed scheme is that both Robin interface conditions depend only on previous-time data, allowing the fluid and shell subproblems to be solved concurrently and in parallel without subiterations. The fluid interacts with the shell through its translational velocity and the corresponding interface force per unit reference area, while the independent shell rotations and angular velocities are governed by the shell equations. At each time step, the Robin conditions use the average of the previous fluid and shell velocities, together with the corresponding average interface loads entering with opposite signs. After the separate subproblems are solved, the updated interface loads are recovered from the Robin conditions and used in the next time step. Although each subproblem is solved implicitly, the interface coupling is explicit. Consequently, the fluid and shell velocities need not match exactly, and the interface forces need not be in exact balance at each time step. The main contribution of this work is the extension of the parallelizable Robin-Robin splitting scheme for fluid thick structure interaction to the fluid-Naghdi shell formulation, together with stability and error analysis that accounts for membrane, transverse-shear, and bending effects, as well as rotational inertia. Our theoretical analysis focuses on the fixed-domain Stokes-Naghdi shell problem. We derive a discrete energy identity and establish unconditional stability for homogeneous external data with any positive time step and fixed positive common Robin coefficient. We also derive stability estimates for nonzero external forces and boundary traction data under appropriate assumptions. Furthermore, we establish first-order temporal consistency and an $O(\sqrt{\Delta t})$ error estimate in the mechanical energy norm under suitable regularity and initialization assumptions. Finally, we assess the proposed scheme for moving-domain fluid-structure interaction through numerical experiments, comparing its results with those obtained from the kinematically coupled $\beta$ scheme and a monolithic solver. 

The paper is organized as follows. Section~\ref{sec:model} introduces the moving-domain fluid-Naghdi shell model. Section~\ref{sec:splitting} presents the explicit Robin-Robin splitting scheme and its numerical implementation. Section~\ref{sec:fixed-domain} describes the fixed-domain problem and the assumptions used in the later analysis. Sections~\ref{sec:energy-stability}, \ref{sec:consistency}, and~\ref{sec:temporal-error} establish the discrete energy identity, consistency estimates, and temporal error bounds, respectively. Finally, Section~\ref{sec:numerical} presents numerical comparisons of the proposed Robin-Robin splitting scheme with the loosely coupled $\beta$ scheme and the monolithic method. It also demonstrates the application of the Fluid-Naghdi shell model and the proposed explicit Robin-Robin splitting scheme to simulate pressure wave propagation through a stented artery.
}

\section{Mathematical model}\label{sec:model}

We consider the interaction between an incompressible Newtonian fluid and a deformable Naghdi shell. The fluid occupies a time-dependent 3D domain, while the shell is described on a fixed 2D parameter domain, with 3D displacement and infinitesimal rotation fields. We assume that the shell midsurface coincides with the fluid-structure interface, neglecting the offset between the physical inner wall and the shell midsurface. The shell is modeled using a linear constitutive law defined on the reference surface, while the fluid domain deforms according to the displacement of the shell midsurface.

\subsection{Reference domains and shell geometry}\label{subsec:geometry}
Let $T>0$ be the final time. Denote the reference fluid domain by
$\widehat{\Omega}_f\subset\mathbb R^3$ and its elastic lateral boundary by
$\widehat{\Gamma}_s$. The remaining boundary is decomposed into a homogeneous
Dirichlet portion and inlet/outlet portions:
\begin{equation}
\partial\widehat{\Omega}_f
=\widehat{\Gamma}_s\cup\widehat{\Sigma}_f^D
\cup\widehat{\Gamma}_f^{\mathrm{in}}
\cup\widehat{\Gamma}_f^{\mathrm{out}}.
\end{equation}
Boundary portions have disjoint relative interiors. 

We consider a cylindrical reference shell of axial length $L>0$, midsurface radius $R>0$, and thickness $\varepsilon>0$. The reference midsurface $\widehat{\Gamma}_s$ is parameterized by the axial coordinate $x$ and the circumferential arc-length coordinate $y_{\mathrm{arc}}$ as
\begin{equation}\label{eq:reference-surface}
\begin{aligned}
\omega=(0,L)\times(0,2\pi R),\qquad
\widehat{\Gamma}_s=\varphi(\omega)=\varphi(x,y_{\mathrm{arc}})
=\left(x,R\cos(y_{\mathrm{arc}}/R),R\sin(y_{\mathrm{arc}}/R)\right).
\end{aligned}
\end{equation}
Here, $y_{\mathrm{arc}}$ denotes the circumferential arc-length coordinate, with periodic identification of the boundaries $y_{\mathrm{arc}}=0$ and $y_{\mathrm{arc}}=2\pi R$. The parameterization is defined on the resulting periodic domain, and its inverse is understood locally or on the cylindrical surface with this identification.

Letting $\theta=y_{\mathrm{arc}}/R$, we define the covariant tangent vectors and the unit normal to the reference midsurface as
\begin{equation}\label{eq:reference-frame}
\begin{aligned}
\boldsymbol a_1=\partial_1\varphi=(1,0,0),\qquad
\boldsymbol a_2=\partial_2\varphi=(0,-\sin\theta,\cos\theta),\qquad
\boldsymbol a_3=\boldsymbol a_1\times\boldsymbol a_2
=(0,-\cos\theta,-\sin\theta),
\end{aligned}
\end{equation}
where $\partial_1=\partial_x$ and
$\partial_2=\partial_{y_{\mathrm{arc}}}$. With this parameterization, $\boldsymbol a_3$ points inward and is opposite to the outward fluid normal on the undeformed lateral boundary. Throughout the paper, fluid tractions are defined using the outward fluid normal $\boldsymbol n_f=-\boldsymbol a_3$. With $a_{\alpha\beta}=\boldsymbol a_\alpha\cdot\boldsymbol a_\beta$,
$\boldsymbol a_c=[a_{\alpha\beta}]$, $\boldsymbol a^c=\boldsymbol a_c^{-1}$,
and $b_{\alpha\beta}=\boldsymbol a_3\cdot\partial_{\alpha\beta}\varphi$,
the reference metric and curvature are
\begin{equation}
\boldsymbol a_c=\boldsymbol a^c=\mathbf I_2,
\qquad
[b_{\alpha\beta}]
=\begin{pmatrix}0&0\\0&1/R\end{pmatrix},
\qquad
dS_0=dx\,dy_{\mathrm{arc}}.
\end{equation}
In particular, no additional factor $R$ belongs in reference shell integrals
when this arc-length parameter is used.

The shell displacement and infinitesimal rotation, both expressed in Cartesian coordinates, are given by
\begin{equation}
\boldsymbol\eta=(\eta_x,\eta_y,\eta_z),
\qquad
\boldsymbol\omega=(\omega_1,\omega_2,\omega_3).
\end{equation}
The deformed midsurface is described by
\begin{equation}\label{eq:current-surface}
\Phi_s(x,y_{\mathrm{arc}},t)
=\varphi(x,y_{\mathrm{arc}})+\boldsymbol\eta(x,y_{\mathrm{arc}},t),
\qquad
\Gamma_s(t)=\Phi_s(\omega,t),
\end{equation}
where $\Gamma_s(t)$ is the moving fluid-shell interface. Assuming that $\Phi_s$ is sufficiently regular and injective, we define the surface Jacobian as
\begin{equation}\label{eq:surface-jacobian}
J_s(x,y_{\mathrm{arc}},t)
=\left|\partial_1\Phi_s(x,y_{\mathrm{arc}},t)
\times\partial_2\Phi_s(x,y_{\mathrm{arc}},t)\right|>0.
\end{equation}
For the moving-domain formulation, we assume that the displacement is sufficiently regular, for example, $\boldsymbol\eta\in W^{1,\infty}$ in space, and that $\Phi_s$ is injective with $J_s$ bounded away from zero. These assumptions ensure a well-defined deformed interface. In contrast, the $H^1$ shell energy space introduced below specifies the weak solution and test spaces but does not, by itself, guarantee an admissible moving fluid domain. Under these assumptions, integration over the moving interface can be expressed on the fixed parameter domain through
\begin{equation}\label{eq:surface-change}
\int_{\Gamma_s(t)} h(\boldsymbol x,t)\,dS
=\int_\omega h(\Phi_s(x,y_{\mathrm{arc}},t),t)
J_s(x,y_{\mathrm{arc}},t)\,dx\,dy_{\mathrm{arc}}.
\end{equation}

\subsection{Arbitrary Lagrangian--Eulerian mapping}\label{subsec:ale}
To describe the moving fluid domain, we introduce the arbitrary Lagrangian-Eulerian (ALE) mapping from the fixed reference domain $\widehat{\Omega}_f$ to the current fluid domain $\Omega_f(t)$:
\begin{equation}\label{eq:ale-map}
\widehat{\varphi}_f(\widehat{\boldsymbol x},t)
=\widehat{\boldsymbol x}+\widehat{\boldsymbol\eta}_f(\widehat{\boldsymbol x},t),
\qquad
\Omega_f(t)=\widehat{\varphi}_f(\widehat{\Omega}_f,t).
\end{equation}
where $\widehat{\boldsymbol\eta}_f$ denotes the fluid mesh displacement. We assume that $\widehat{\varphi}_f$ is an orientation-preserving diffeomorphism, with deformation gradient and Jacobian given by
\begin{equation}
\widehat{\mathbf F}_f=\widehat\nabla\widehat{\varphi}_f,
\qquad
\widehat J_f=\det\widehat{\mathbf F}_f>0.
\end{equation}
Here, $\widehat J_f$ is the volume Jacobian of the ALE mapping, which is distinct from the surface Jacobian $J_s$ defined in~\eqref{eq:surface-jacobian}. Throughout the paper, we use the vector-gradient convention $(\nabla\boldsymbol v)_{ij}=\partial_j v_i$. Under the ALE mapping, the spatial gradient transforms according to
\begin{equation}
\nabla\boldsymbol v
=(\widehat\nabla\widehat{\boldsymbol v})
\widehat{\mathbf F}_f^{-1},
\qquad
\widehat{\boldsymbol v}
=\boldsymbol v\circ\widehat{\varphi}_f,
\end{equation}
where the gradients are evaluated at corresponding points in the current and reference domains.

At the fluid-shell interface, the mesh displacement matches the shell displacement, giving
\begin{equation}\label{eq:ale-interface}
\widehat{\boldsymbol\eta}f(\varphi(x,y{\mathrm{arc}}),t)
=\boldsymbol\eta(x,y_{\mathrm{arc}},t),
\qquad
\widehat{\varphi}f(\varphi(x,y{\mathrm{arc}}),t)
=\Phi_s(x,y_{\mathrm{arc}},t).
\end{equation}
The displacement is then extended into the fluid domain to determine the interior mesh motion. For fixed external boundaries, one possible choice is the harmonic extension
\begin{equation}\label{eq}
\left\{\begin{aligned}
-\widehat\Delta\widehat{\boldsymbol\eta}_f&=\boldsymbol0
&&\text{in }\widehat\Omega_f,\\
\widehat{\boldsymbol\eta}_f\circ\varphi&=\boldsymbol\eta
&&\text{on }\omega, \\ \widehat{\boldsymbol\eta}_f&=\boldsymbol0
&&\text{on }\partial\widehat\Omega_f\setminus\widehat\Gamma_s.
\end{aligned}\right.
\end{equation}
Other mesh extension methods may also be used, provided that the interface displacement is preserved and the ALE mapping remains admissible.

The fluid mesh velocity and the pullback of the fluid velocity to the reference domain are defined by
\begin{equation}\label{eq:mesh-velocity}
\boldsymbol w_f(\boldsymbol x,t)
=\partial_t\widehat{\varphi}_f
(\widehat{\varphi}_f^{-1}(\boldsymbol x,t),t),\qquad
\widehat{\boldsymbol u}(\widehat{\boldsymbol x},t)
=\boldsymbol u(\widehat{\varphi}_f(\widehat{\boldsymbol x},t),t).
\end{equation}
Applying the chain rule gives the ALE time derivative
\begin{equation}\label{eq:ale-derivative}
\left.\partial_t\boldsymbol u\right|_{\widehat\Omega_f}
=\partial_t\boldsymbol u+(\boldsymbol w_f\cdot\nabla)\boldsymbol u.
\end{equation}

\subsection{Fluid equations}
Let $\rho_f>0$ and $\mu_f>0$ denote the fluid density and dynamic viscosity, respectively. The fluid velocity $\boldsymbol u$ and pressure $p$ satisfy the incompressible Navier-Stokes equations in the ALE formulation:
\begin{equation}\label{eq:moving-navier-stokes}
\left\{\begin{aligned}
\rho_f\left(\left.\partial_t\boldsymbol u\right|_{\widehat\Omega_f}
+((\boldsymbol u-\boldsymbol w_f)\cdot\nabla)\boldsymbol u\right)
-\nabla\cdot\boldsymbol\sigma_f(\boldsymbol u,p)
&=\boldsymbol F_f
&&\text{in }\Omega_f(t),\\
\nabla\cdot\boldsymbol u&=0
&&\text{in }\Omega_f(t),
\end{aligned}\right.
\end{equation}
where
\begin{equation}\label{eq:fluid-stress}
\boldsymbol\sigma_f(\boldsymbol u,p)
=-p\mathbf I_3+2\mu_f\mathbb D(\boldsymbol u),
\qquad
\mathbb D(\boldsymbol u)=\frac{\nabla\boldsymbol u+(\nabla\boldsymbol u)^T}{2}.
\end{equation}
Here $\boldsymbol F_f$ is the external force.

Let $\Sigma_f^D(t)$, $\Gamma_f^{\mathrm{in}}(t)$, and $\Gamma_f^{\mathrm{out}}(t)$ denote the current Dirichlet, inlet, and outlet boundaries, respectively. We impose the following boundary conditions:
\begin{equation}\label{eq:fluid-boundaries}
\begin{aligned}
\boldsymbol u&=\boldsymbol0&&\text{on }\Sigma_f^D(t),\\
\boldsymbol\sigma_f(\boldsymbol u,p)\boldsymbol n_f
&=-p_{\mathrm{in}}\boldsymbol n_f&&\text{on }\Gamma_f^{\mathrm{in}}(t),\\
\boldsymbol\sigma_f(\boldsymbol u,p)\boldsymbol n_f
&=-p_{\mathrm{out}}\boldsymbol n_f&&\text{on }\Gamma_f^{\mathrm{out}}(t).
\end{aligned}
\end{equation}
Here, $p_{\mathrm{in}}$ and $p_{\mathrm{out}}$ specify the normal traction at the inlet and outlet, rather than Dirichlet conditions on the fluid pressure. Alternatively, a prescribed inlet velocity may be used in place of the inlet traction condition, with appropriate modifications to the fluid trial and test spaces. In the subsequent analysis, we use the boundary conditions given here.

\subsection{Naghdi shell equation}
In the Naghdi shell model, the infinitesimal rotation $\boldsymbol\omega$ is treated as an independent variable rather than being determined by the displacement gradient. We express the shell strains in the local orthonormal frame $(\boldsymbol a_1,\boldsymbol a_2,\boldsymbol a_3)$, represented by
\begin{equation}
\mathbf Q=\begin{pmatrix}
1&0&0\\0&-\sin\theta&\cos\theta\\0&-\cos\theta&-\sin\theta
\end{pmatrix}.
\end{equation}
The displacement and rotation enter the strain formulation through the vectors
$\partial_i\boldsymbol\eta+\boldsymbol a_i\times\boldsymbol\omega$, $i=1,2$.
Projecting these vectors onto the tangential directions of the reference midsurface gives the membrane-shear strain block
\begin{equation}\label{eq:membrane-block}
\mathbf A(\boldsymbol\eta,\boldsymbol\omega)=\begin{pmatrix}
\boldsymbol a_1\cdot\partial_1\boldsymbol\eta&
\boldsymbol a_1\cdot\partial_2\boldsymbol\eta
{+\boldsymbol\omega\cdot\boldsymbol a_3}\\[1mm]
\boldsymbol a_2\cdot\partial_1\boldsymbol\eta
{-\boldsymbol\omega\cdot\boldsymbol a_3}&
\boldsymbol a_2\cdot\partial_2\boldsymbol\eta
\end{pmatrix},
\end{equation}
The transverse shear components are obtained by projecting the same vectors onto the reference normal $\boldsymbol a_3$:
\begin{equation}\label{eq:shear-row}
\boldsymbol e(\boldsymbol\eta,\boldsymbol\omega)^T
=\begin{pmatrix}
\boldsymbol a_3\cdot\partial_1\boldsymbol\eta
+\boldsymbol\omega\cdot\boldsymbol a_2&
\boldsymbol a_3\cdot\partial_2\boldsymbol\eta
-\boldsymbol\omega\cdot\boldsymbol a_1
\end{pmatrix}.
\end{equation}
Thus, combining the tangential and transverse shear components, we obtain the membrane-shear strain tensor
\begin{equation}\label{eq:membrane-shear-strain}
\widehat{\mathbf E}_{\mathrm{mem}}(\boldsymbol\eta,\boldsymbol\omega)
=\begin{pmatrix}
\mathbf A(\boldsymbol\eta,\boldsymbol\omega)\\
\boldsymbol e(\boldsymbol\eta,\boldsymbol\omega)^T
\end{pmatrix}\in\mathbb R^{3\times2}.
\end{equation}
The full tangential block and transverse shear components are retained in the constitutive law, accounting for the contributions of the independent rotation and its coupling with the displacement gradients.

For test functions $(\boldsymbol v,\boldsymbol w)$, the corresponding membrane-shear strain components are defined in the same way:
\begin{equation}\label{eq:test-membrane-block}
\mathbf D(\boldsymbol v,\boldsymbol w)=\begin{pmatrix}
\boldsymbol a_1\cdot\partial_1\boldsymbol v&
\boldsymbol a_1\cdot\partial_2\boldsymbol v
{+\boldsymbol w\cdot\boldsymbol a_3}\\[1mm]
\boldsymbol a_2\cdot\partial_1\boldsymbol v
{-\boldsymbol w\cdot\boldsymbol a_3}&
\boldsymbol a_2\cdot\partial_2\boldsymbol v
\end{pmatrix},
\end{equation}
\begin{equation}\label{eq:test-shear-row}
\boldsymbol d(\boldsymbol v,\boldsymbol w)^T
=\begin{pmatrix}
\boldsymbol a_3\cdot\partial_1\boldsymbol v+\boldsymbol w\cdot\boldsymbol a_2&
\boldsymbol a_3\cdot\partial_2\boldsymbol v-\boldsymbol w\cdot\boldsymbol a_1
\end{pmatrix},
\qquad
\widehat{\mathbf D}_{\mathrm{mem}}(\boldsymbol v,\boldsymbol w)
=\begin{pmatrix}\mathbf D(\boldsymbol v,\boldsymbol w)\\
\boldsymbol d(\boldsymbol v,\boldsymbol w)^T\end{pmatrix}.
\end{equation}
The bending strain and its corresponding test strain are obtained by projecting the derivatives of the rotation and test functions onto the local reference frame:
\begin{align}
\widehat{\mathbf E}_{\mathrm{bend}}(\boldsymbol\omega)
=
\begin{pmatrix}
\boldsymbol a_1\cdot\partial_1\boldsymbol\omega &
\boldsymbol a_1\cdot\partial_2\boldsymbol\omega\\
\boldsymbol a_2\cdot\partial_1\boldsymbol\omega &
\boldsymbol a_2\cdot\partial_2\boldsymbol\omega\\
\boldsymbol a_3\cdot\partial_1\boldsymbol\omega &
\boldsymbol a_3\cdot\partial_2\boldsymbol\omega
\end{pmatrix},
\label{eq:bending-strain}\quad
\widehat{\mathbf D}_{\mathrm{bend}}(\boldsymbol w)
=
\begin{pmatrix}
\boldsymbol a_1\cdot\partial_1\boldsymbol w &
\boldsymbol a_1\cdot\partial_2\boldsymbol w\\
\boldsymbol a_2\cdot\partial_1\boldsymbol w &
\boldsymbol a_2\cdot\partial_2\boldsymbol w\\
\boldsymbol a_3\cdot\partial_1\boldsymbol w &
\boldsymbol a_3\cdot\partial_2\boldsymbol w
\end{pmatrix}.
\end{align}
For use in the constitutive law, we separate the tangential blocks and normal rows:
\begin{equation}\label{eq:bending-blocks}
\widehat{\mathbf E}_{\mathrm{bend}}(\boldsymbol\omega)
=\begin{pmatrix}\mathbf E_b(\boldsymbol\omega)\\
\boldsymbol e_b(\boldsymbol\omega)^T\end{pmatrix},
\qquad
\widehat{\mathbf D}_{\mathrm{bend}}(\boldsymbol w)
=\begin{pmatrix}\mathbf D_b(\boldsymbol w)\\
\boldsymbol d_b(\boldsymbol w)^T\end{pmatrix},
\end{equation}
where $\mathbf E_b,\mathbf D_b\in\mathbb R^{2\times2}$ and $\boldsymbol e_b,\boldsymbol d_b\in\mathbb R^2$. Since the rotation field is expressed in Cartesian coordinates, its derivatives are taken before projection onto the local reference frame. This distinction is important because differentiating the local rotation components directly would require additional terms arising from variations of the reference frame.

Let $\mathcal C_m$ and $\mathcal C_b$ denote the elasticity tensors
defined on the reference midsurface. Using the matrix contraction
$\mathbf A:\mathbf B=\sum_{i,j}A_{ij}B_{ij}$, we define
\begin{equation}\label{eq:membrane-form}
\begin{aligned}
B_m\big((\boldsymbol\eta,\boldsymbol\omega),(\boldsymbol v,\boldsymbol w)\big)=\varepsilon\int_\omega
\mathcal C_m\widehat{\mathbf E}_{\mathrm{mem}}(\boldsymbol\eta,\boldsymbol\omega)
:\widehat{\mathbf D}_{\mathrm{mem}}(\boldsymbol v,\boldsymbol w)
\,dx\,dy_{\mathrm{arc}},
\end{aligned}
\end{equation}
\begin{equation}\label{eq:bending-form}
\begin{aligned}
B_b\big((\boldsymbol\eta,\boldsymbol\omega),(\boldsymbol v,\boldsymbol w)\big)
=\frac{\varepsilon^3}{12}\int_\omega
\mathcal C_b\widehat{\mathbf E}_{\mathrm{bend}}(\boldsymbol\omega)
:\widehat{\mathbf D}_{\mathrm{bend}}(\boldsymbol w)
\,dx\,dy_{\mathrm{arc}}.
\end{aligned}
\end{equation}
For the cylindrical shell, let $\lambda$ and $\mu$ denote the
Lamé parameters of the shell material. We assume an isotropic material
with Young's modulus $E_s>0$ and Poisson's ratio
$-1<\nu_s\leq 0.5$. The reduced elasticity coefficients are
\begin{equation}\label{eq:material-coefficients}
\frac{2\lambda\mu}{\lambda+2\mu}
=\frac{E_s\nu_s}{1-\nu_s^2},
\qquad
2\mu=\frac{E_s}{1+\nu_s}.
\end{equation}
At $\nu_s=0.5$, these coefficients are defined by continuity
and remain finite, so the singular 3D Lamé parameter $\lambda$
need not be evaluated. The membrane-shear constitutive law is given by
\begin{equation}\label{eq:membrane-material}
\begin{aligned}
&\mathcal C_m
\widehat{\mathbf E}_{\mathrm{mem}}
(\boldsymbol\eta,\boldsymbol\omega)
:\widehat{\mathbf D}_{\mathrm{mem}}
(\boldsymbol v,\boldsymbol w)=
\frac{2\lambda\mu}{\lambda+2\mu}
\operatorname{tr}\mathbf A(\boldsymbol\eta,\boldsymbol\omega)
\operatorname{tr}\mathbf D(\boldsymbol v,\boldsymbol w)\\
&\qquad+
2\mu\,\mathbf A(\boldsymbol\eta,\boldsymbol\omega)
:\mathbf D(\boldsymbol v,\boldsymbol w)+
\mathcal B_m\boldsymbol e(\boldsymbol\eta,\boldsymbol\omega)
\cdot\boldsymbol d(\boldsymbol v,\boldsymbol w),
\end{aligned}
\end{equation}
where $\mathcal B_m=\mu\mathbf I_2$ is the transverse
shear coefficient for the arc-length parameterized cylinder. 

For the bending constitutive law, we introduce
\begin{equation}
\mathbf J=
\begin{pmatrix}
0&-1\\
1&0
\end{pmatrix}.
\end{equation}
The bending contribution is then defined by
\begin{equation}\label{eq:bending-material}
\begin{aligned}
&\mathcal C_b
\widehat{\mathbf E}_{\mathrm{bend}}(\boldsymbol\omega)
:\widehat{\mathbf D}_{\mathrm{bend}}(\boldsymbol w)
=
\frac{2\lambda\mu}{\lambda+2\mu}
\operatorname{tr}\big(
\mathbf J\mathbf E_b(\boldsymbol\omega)\big)
\operatorname{tr}\big(
\mathbf J\mathbf D_b(\boldsymbol w)\big)\\
&\qquad+
2\mu\,
\big(\mathbf J\mathbf E_b(\boldsymbol\omega)\big)
:\big(\mathbf J\mathbf D_b(\boldsymbol w)\big)
+
\mu\,\boldsymbol e_b(\boldsymbol\omega)
\cdot\boldsymbol d_b(\boldsymbol w).
\end{aligned}
\end{equation}
This formulation includes both the tangential and normal
components of the rotation derivatives, with the normal-row
coefficient $\mathcal B_b=\mu\mathbf I_2$. The bending
thickness factor $\varepsilon^3/12$ is included separately
in~\eqref{eq:bending-form}. Both constitutive laws use the
reduced coefficients in~\eqref{eq:material-coefficients},
including their limiting values at $\nu_s=0.5$. The stability and error analysis below relies on the
symmetry and nonnegativity of the resulting bilinear forms.
We also note that replacing $\mathbf J$ by $-\mathbf J$
does not change the bending constitutive law, since
$\mathbf J$ appears in both factors of each tangential term.

We assume that the bilinear forms $B_m$ and $B_b$ are continuous, time-independent, symmetric, and nonnegative:
\begin{equation}\label{eq:form-assumptions}
\begin{aligned}
&B_m\big((\boldsymbol\eta,\boldsymbol\omega),(\boldsymbol v,\boldsymbol w)\big)
=B_m\big((\boldsymbol v,\boldsymbol w),(\boldsymbol\eta,\boldsymbol\omega)\big),\qquad
B_m\big((\boldsymbol\eta,\boldsymbol\omega),(\boldsymbol\eta,\boldsymbol\omega)\big)\geq0,\\
&B_b\big((\boldsymbol\eta,\boldsymbol\omega),(\boldsymbol v,\boldsymbol w)\big)
=B_b\big((\boldsymbol v,\boldsymbol w),(\boldsymbol\eta,\boldsymbol\omega)\big)
,\qquad 
B_b\big((\boldsymbol\eta,\boldsymbol\omega),(\boldsymbol\eta,\boldsymbol\omega)\big)\geq0.
\end{aligned}
\end{equation}
These properties are sufficient for the discrete energy identities derived below. However, controlling the full shell norm requires an additional coercivity estimate on the chosen shell space. Such coercivity does not follow from the positivity of the material coefficients alone.

We impose clamped boundary conditions at both ends of the shell and periodic conditions in the circumferential direction:
\begin{equation}\label{eq:shell-boundaries}
\begin{aligned}
&\boldsymbol\eta(0,y_{\mathrm{arc}},t)
=\boldsymbol\eta(L,y_{\mathrm{arc}},t)=\boldsymbol0,\qquad
\boldsymbol\omega(0,y_{\mathrm{arc}},t)
=\boldsymbol\omega(L,y_{\mathrm{arc}},t)=\boldsymbol0,\\
&\boldsymbol\eta(x,0,t)=\boldsymbol\eta(x,2\pi R,t),\qquad \qquad ~~
\boldsymbol\omega(x,0,t)=\boldsymbol\omega(x,2\pi R,t).
\end{aligned}
\end{equation}
The corresponding shell space is
\begin{equation}\label{eq:shell-space}
\begin{aligned}
V_N(\omega)=\big\{(\boldsymbol v,\boldsymbol w)
\in H^1(\omega;\mathbb R^3)\times H^1(\omega;\mathbb R^3):\;&
\boldsymbol v=\boldsymbol w=\boldsymbol0\text{ at }x=0,L,\\
&\boldsymbol v,\boldsymbol w\text{ have periodic circumferential traces}\big\}.
\end{aligned}
\end{equation}
The boundary conditions in $V_N(\omega)$ are imposed in the trace sense. In particular, periodicity requires matching traces at the circumferential boundaries, without additional conditions on their derivatives.

Let $\rho_s>0$ denote the shell density. We introduce the translational and angular velocities through
$\partial_t\boldsymbol\eta=\dot{\boldsymbol\eta},~
\partial_t\boldsymbol\omega=\dot{\boldsymbol\omega}.$
The weak form of the shell equations is: find $(\boldsymbol\eta,\boldsymbol\omega)$ such that, for every $(\boldsymbol v,\boldsymbol w)\in V_N(\omega)$,
\begin{equation}\label{eq:shell-dynamics}
\begin{aligned}
&\varepsilon\rho_s\int_\omega
\partial_t\dot{\boldsymbol\eta}\cdot\boldsymbol v\,dx\,dy_{\mathrm{arc}}
+\frac{\varepsilon^3\rho_s}{12}\int_\omega
\partial_t\dot{\boldsymbol\omega}\cdot\boldsymbol w\,dx\,dy_{\mathrm{arc}}
+B_m\big((\boldsymbol\eta,\boldsymbol\omega),(\boldsymbol v,\boldsymbol w)\big)
+B_b\big((\boldsymbol\eta,\boldsymbol\omega),(\boldsymbol v,\boldsymbol w)\big)\\
&=\int_\omega\boldsymbol f_s^{\mathrm{ref}}\cdot\boldsymbol v
\,dx\,dy_{\mathrm{arc}}.
\end{aligned}
\end{equation}
Here, $\varepsilon\rho_s$ and $\varepsilon^3\rho_s/12$ represent the translational mass and rotational inertia per unit reference area, respectively. The term $\boldsymbol f_s^{\mathrm{ref}}$ denotes the applied force per unit reference parameter area. No external couple is applied to the rotational equations.

\subsection{Interface coupling conditions}\label{subsec:coupling}

The fluid and shell are coupled through the kinematic condition
\begin{equation}\label{eq:kinematic-coupling}
\boldsymbol u(\Phi_s(x,y_{\mathrm{arc}},t),t)
=\dot{\boldsymbol\eta}(x,y_{\mathrm{arc}},t),
\end{equation}
which enforces no-slip and no-penetration at the moving interface. Differentiating~\eqref{eq:ale-interface} with respect to time gives
\begin{equation}
\boldsymbol w_f\circ\Phi_s=\dot{\boldsymbol\eta},
\qquad
\boldsymbol u=\boldsymbol w_f
\quad\text{on }\Gamma_s(t).
\end{equation}
Thus, the fluid velocity and mesh velocity coincide at the interface. Since $\boldsymbol\omega$ represents an independent infinitesimal shell rotation, no additional condition relating $\dot{\boldsymbol\omega}$ to the fluid velocity is imposed.

The fluid traction on the current interface is defined by
\begin{equation}\label{eq:current-fluid-traction}
\boldsymbol t_f(\boldsymbol x,t)
=\boldsymbol\sigma_f(\boldsymbol u(\boldsymbol x,t),p(\boldsymbol x,t))
\boldsymbol n_f(\boldsymbol x,t),
\qquad \boldsymbol x\in\Gamma_s(t).
\end{equation}
Using the surface change of variables, we express this traction as a force density on the reference parameter domain:
\begin{equation}\label{eq:reference-fluid-traction}
\begin{aligned}
\widehat{\boldsymbol t}_f(x,y_{\mathrm{arc}},t)
={}&J_s(x,y_{\mathrm{arc}},t)
\boldsymbol\sigma_f(\boldsymbol u,p)
(\Phi_s(x,y_{\mathrm{arc}},t),t) \boldsymbol n_f(\Phi_s(x,y_{\mathrm{arc}},t),t).
\end{aligned}
\end{equation}
Here, $\widehat{\boldsymbol t}_f$ includes the surface Jacobian $J_s$ and represents force per unit reference parameter area, rather than simply the pullback of the current traction.

Let $\boldsymbol f_s^\Gamma$ denote the force acting on the shell per unit current surface area. By the balance of interface forces,
\begin{equation}\label{eq:dynamic-coupling}
\boldsymbol f_s^\Gamma=-\boldsymbol t_f,
\qquad
\boldsymbol f_s^{\mathrm{ref}}
=(\boldsymbol f_s^\Gamma\circ\Phi_s)J_s
=-\widehat{\boldsymbol t}_f.
\end{equation}
In the Naghdi shell model, the internal structural forces are represented by the membrane and bending bilinear forms~\eqref{eq:membrane-form} and~\eqref{eq:bending-form}. Therefore, no separate 3D structural traction tensor is introduced.

For compatible fluid and shell test functions satisfying $\boldsymbol{\tilde v}\circ\Phi_s=\boldsymbol v$, the interface work contributions cancel:
\begin{equation}\label{eq:interface-work}
\begin{aligned}
&\int_{\Gamma_s(t)}
\boldsymbol t_f\cdot\boldsymbol{\tilde v}\,dS
+\int_\omega
\boldsymbol f_s^{\mathrm{ref}}\cdot\boldsymbol v
\,dx\,dy_{\mathrm{arc}}
=
\int_\omega
(\widehat{\boldsymbol t}_f+\boldsymbol f_s^{\mathrm{ref}})
\cdot\boldsymbol v\,dx\,dy_{\mathrm{arc}}
=0.
\end{aligned}
\end{equation}
This identity establishes the sign convention and reference surface measure used in the Robin-Robin splitting scheme.

\subsection{Initial conditions and compatibility conditions}

We prescribe the initial fluid velocity, shell displacement, translational velocity, rotation, and angular velocity:
\begin{equation}\label{eq:initial-data}
\begin{aligned}
&\boldsymbol u(\boldsymbol x,0)
=\boldsymbol u^0(\boldsymbol x),\qquad
\boldsymbol\eta(x,y_{\mathrm{arc}},0)
=\boldsymbol\eta^0(x,y_{\mathrm{arc}}),\qquad 
\dot{\boldsymbol\eta}(x,y_{\mathrm{arc}},0)
=\dot{\boldsymbol\eta}^0(x,y_{\mathrm{arc}}),\\
&\boldsymbol\omega(x,y_{\mathrm{arc}},0)
=\boldsymbol\omega^0(x,y_{\mathrm{arc}}),\qquad
\dot{\boldsymbol\omega}(x,y_{\mathrm{arc}},0)
=\dot{\boldsymbol\omega}^0(x,y_{\mathrm{arc}}).
\end{aligned}
\end{equation}
The initial fluid domain is determined by the initial shell displacement through the ALE mapping~\eqref{eq:ale-map} and the interface condition~\eqref{eq:ale-interface}. In particular, when the reference configuration is chosen to coincide with the initial configuration, we take $\boldsymbol\eta^0=\boldsymbol0$ and $\widehat{\boldsymbol\eta}_f(\cdot,0)=\boldsymbol0$.

For sufficiently regular solutions, the initial data must satisfy the incompressibility and kinematic compatibility conditions
\begin{equation}\label{eq:initial-compatibility}
\begin{aligned}
\nabla\cdot\boldsymbol u^0=0
~~\text{in }\Omega_f(0),\qquad
\boldsymbol u^0\circ\Phi_s(\cdot,0)
=\dot{\boldsymbol\eta}^0
~~\text{on }\omega.
\end{aligned}
\end{equation}
The initial data must also satisfy the prescribed external velocity conditions, clamped shell boundary conditions, and circumferential periodicity. The fluid pressure serves as a Lagrange multiplier enforcing incompressibility and is not prescribed independently as an initial evolution variable. The initial interface loads required by the partitioned scheme are introduced together with the time discretization.

\subsection{The weak formulation of the coupled system}\label{subsec:coupled-weak}
For each current domain, we define
\begin{equation}\label{eq:fluid-spaces}
\begin{aligned}
V_f(t)=\{\boldsymbol{\tilde v}\in H^1(\Omega_f(t);\mathbb R^3):
\boldsymbol{\tilde v}=\boldsymbol0\text{ on }\Sigma_f^D(t)\},\qquad
Q_f(t)=L^2(\Omega_f(t)).
\end{aligned}
\end{equation}
For the stated boundary conditions, the inlet and outlet
traction boundaries have positive surface measure and fix the pressure
level. Accordingly, $Q_f(t)=L^2(\Omega_f(t))$ is used without an additional
zero-mean constraint. A different boundary specification that leaves a
pressure nullspace would require a corresponding pressure normalization.

Find $(\boldsymbol\eta,\boldsymbol\omega)\in V_N(\omega)$, their velocities,
and $(\boldsymbol u,p)\in V_f(t)\times Q_f(t)$, satisfying the geometry,
kinematics, and initial data above and
$\boldsymbol u\circ\Phi_s=\dot{\boldsymbol\eta}$, such that for all
$(\boldsymbol{\tilde v},q,\boldsymbol v,\boldsymbol w)$ with
$(\boldsymbol{\tilde v},q)\in V_f(t)\times Q_f(t)$,
$(\boldsymbol v,\boldsymbol w)\in V_N(\omega)$, and
$\boldsymbol{\tilde v}\circ\Phi_s=\boldsymbol v$, one has
\begin{equation}\label{eq:coupled-weak}
\begin{aligned}
&\rho_f\int_{\Omega_f(t)}
\left.\partial_t\boldsymbol u\right|_{\widehat\Omega_f}
\cdot\boldsymbol{\tilde v}\,d\boldsymbol x
+\rho_f\int_{\Omega_f(t)}
[((\boldsymbol u-\boldsymbol w_f)\cdot\nabla)\boldsymbol u]
\cdot\boldsymbol{\tilde v}\,d\boldsymbol x
+2\mu_f\int_{\Omega_f(t)}
\mathbb D(\boldsymbol u):\mathbb D(\boldsymbol{\tilde v})\,d\boldsymbol x \\
&
-\int_{\Omega_f(t)}p\nabla\cdot\boldsymbol{\tilde v}\,d\boldsymbol x
+\int_{\Omega_f(t)}q\nabla\cdot\boldsymbol u\,d\boldsymbol x\\
&+\varepsilon\rho_s\int_\omega
\partial_t\dot{\boldsymbol\eta}\cdot\boldsymbol v\,dx\,dy_{\mathrm{arc}}
+\frac{\varepsilon^3\rho_s}{12}\int_\omega
\partial_t\dot{\boldsymbol\omega}\cdot\boldsymbol w\,dx\,dy_{\mathrm{arc}}
+B_m\big((\boldsymbol\eta,\boldsymbol\omega),(\boldsymbol v,\boldsymbol w)\big)
+B_b\big((\boldsymbol\eta,\boldsymbol\omega),(\boldsymbol v,\boldsymbol w)\big)\\
&=\int_{\Omega_f(t)}\boldsymbol F_f\cdot\boldsymbol{\tilde v}\,d\boldsymbol x
-\int_{\Gamma_f^{\mathrm{in}}(t)}p_{\mathrm{in}}
\boldsymbol n_f\cdot\boldsymbol{\tilde v}\,dS_t
-\int_{\Gamma_f^{\mathrm{out}}(t)}p_{\mathrm{out}}
\boldsymbol n_f\cdot\boldsymbol{\tilde v}\,dS_t.
\end{aligned}
\end{equation}
This formulation follows by integration by parts in the fluid momentum
equation and addition of Naghdi shell equation \eqref{eq:shell-dynamics}, while the interface terms
vanish by \eqref{eq:interface-work}. For lower regularity, boundary traction
work is interpreted in the appropriate trace duality. Its representation
as an $L^2$ surface integral requires additional traction regularity, which
will be stated explicitly for the splitting analysis.

\begin{remark}[Scope of the subsequent analysis]\label{rem:scope}
The model in this section retains the moving fluid domain and nonlinear
fluid convection. Section~3 first states a fully implicit
fluid formulation with the physical traction conditions, then describes
the implemented lagged skew-convection variant separately. The discrete energy and temporal error proofs will then be developed
for the fixed-domain Stokes--Naghdi linearization. Freezing the domain does
not suppress the shell displacement or rotation unknowns; it removes their
feedback on the fluid geometry. The resulting theorems do not, by themselves,
assert stability or convergence for the moving-domain problem.
\end{remark}

\section{Explicit Robin--Robin splitting scheme}\label{sec:splitting}
We now discretize the moving-domain problem \eqref{eq:coupled-weak} in time. During each step, the fluid and shell exchange only data from the previous time level. This
allows the two subproblems to be solved independently, without
subiterations. The word explicit refers to this interface exchange. In this section, we will first present the fully implicit fluid formulation, then
describe the lagged skew-convection form used in the implementation.

% We use the domains, shell space $V_N(\omega)$, and bilinear forms
% $B_m$ and $B_b$ defined in Section~\ref{sec:model}. In particular,
% $B_m$ includes both membrane and transverse-shear effects. All interface
% loads below are force densities per unit reference parameter area.

Let $N\geq1$ be an integer, $\Delta t=T/N$, and $t^n=n\Delta t$.
At time $t^n$, given that
\begin{equation}\label{eq:known-data}
\widehat{\boldsymbol u}^{n},\quad
\widehat p^{n},\quad
\boldsymbol\eta^n,\quad\dot{\boldsymbol\eta}^n,\quad
\boldsymbol\omega^n,\quad\dot{\boldsymbol\omega}^n,\quad
\widehat{\boldsymbol t}_f^n,\quad
\boldsymbol f_s^{\mathrm{ref},n}
\end{equation}
are available, together with the fluid map $\widehat\varphi_f^n$
and mesh displacement $\widehat{\boldsymbol\eta}_f^n$.
Here, the two stored loads have separate roles.
$\widehat{\boldsymbol t}_f^n$ is the reference fluid traction, including
the surface Jacobian, and $\boldsymbol f_s^{\mathrm{ref},n}$ is the load
applied to the shell. For the continuous solution, action and reaction give
$\boldsymbol f_s^{\mathrm{ref}}=-\widehat{\boldsymbol t}_f$.
The discrete update permits a force imbalance as the two subproblems advance independently.

\subsection{Robin conditions at the interface}\label{subsec:robin-construction}
The exact interface conditions on the reference parameter domain are
\begin{equation}\label{eq:exact-interface}
\boldsymbol u\circ\Phi_s=\dot{\boldsymbol\eta},
\qquad
\widehat{\boldsymbol t}_f+\boldsymbol f_s^{\mathrm{ref}}=\boldsymbol0.
\end{equation}
For our splitting scheme, we choose the following two Robin conditions at the interface:
\begin{equation}\label{eq:fluid-robin}
\begin{aligned}
\widehat{\boldsymbol t}_f^{n+1}
+\gamma\widehat{\boldsymbol u}^{n+1}\circ\varphi
={}&\frac{\gamma}{2}
\left(\widehat{\boldsymbol u}^{n}\circ\varphi
+\dot{\boldsymbol\eta}^{n}\right)
+\frac12\left(\widehat{\boldsymbol t}_f^{n}
-\boldsymbol f_s^{\mathrm{ref},n}\right),
\end{aligned}
\end{equation}
\begin{equation}\label{eq:shell-robin}
\begin{aligned}
\boldsymbol f_s^{\mathrm{ref},n+1}
+\gamma\dot{\boldsymbol\eta}^{n+1}
={}&\frac{\gamma}{2}
\left(\widehat{\boldsymbol u}^{n}\circ\varphi
+\dot{\boldsymbol\eta}^{n}\right)
-\frac12\left(\widehat{\boldsymbol t}_f^{n}
-\boldsymbol f_s^{\mathrm{ref},n}\right).
\end{aligned}
\end{equation}
Here, we choose the same constant $\gamma>0$ for both Robin conditions, and $\gamma$ is an independent coupling parameter. Its units in $\mathrm{kg\,m^{-2}\,s^{-1}}$, given that multiplying a velocity by $\gamma$ gives a reference force density. 

Both right-hand sides of \eqref{eq:fluid-robin} and \eqref{eq:shell-robin} use the same average of the previous fluid and
shell velocities, while the previous loads enter with opposite signs, following
the action-reaction convention. Since all these quantities are already known, neither subproblem needs the new solution of the other. To check the relation with the continuous coupling conditions,
set all time levels in \eqref{eq:fluid-robin}--\eqref{eq:shell-robin}
to the same values. The two equations then reduce to
\begin{equation}
\begin{aligned}
\widehat{\boldsymbol t}_f+\boldsymbol f_s^{\mathrm{ref}}
+\gamma(\widehat{\boldsymbol u}\circ\varphi-\dot{\boldsymbol\eta})
&=\boldsymbol0,\\
\widehat{\boldsymbol t}_f+\boldsymbol f_s^{\mathrm{ref}}
-\gamma(\widehat{\boldsymbol u}\circ\varphi-\dot{\boldsymbol\eta})
&=\boldsymbol0.
\end{aligned}
\end{equation}
Adding these equations gives force balance; subtracting them gives
velocity continuity.

\subsection{Fluid domain update}
\label{subsec:geometry-prediction}
In our splitting scheme, the fluid branch needs a domain update before the new shell solution is
available. We predict its boundary from the old shell displacement and velocity:
\begin{equation}\label{eq:predicted-surface}
\Phi_s^{n+1,*}(x,y_{\mathrm{arc}})
=\varphi(x,y_{\mathrm{arc}})
+\boldsymbol\eta^n(x,y_{\mathrm{arc}})
+\Delta t\,\dot{\boldsymbol\eta}^n(x,y_{\mathrm{arc}}).
\end{equation}
This equation will provide the predicted fluid interface. Next, we extend the predicted boundary displacement into the reference fluid domain by solving the harmonic equation
\begin{equation}\label{eq:discrete-extension}
\left\{\begin{aligned}
-\widehat\Delta\widehat{\boldsymbol\eta}_f^{n+1}
&=\boldsymbol0&&\text{in }\widehat\Omega_f,\\
\widehat{\boldsymbol\eta}_f^{n+1}\circ\varphi
&=\boldsymbol\eta^n+\Delta t\,\dot{\boldsymbol\eta}^n
&&\text{on }\omega,\\
\widehat{\boldsymbol\eta}_f^{n+1}
&=\boldsymbol0&&\text{on }\partial\widehat\Omega_f\setminus\widehat\Gamma_s.
\end{aligned}\right.
\end{equation}
This extension defines the new fluid map, its gradient and domain velocity in reference domain:
\begin{equation}\label{eq:numerical-map}
\begin{aligned}
\widehat\varphi_f^{n+1}(\widehat{\boldsymbol x})
&=\widehat{\boldsymbol x}
+\widehat{\boldsymbol\eta}_f^{n+1}(\widehat{\boldsymbol x}),\qquad
\Omega_f^{n+1}=\widehat\varphi_f^{n+1}(\widehat\Omega_f),\qquad
\Gamma_s^{n+1,*}=\Phi_s^{n+1,*}(\omega),\\
&\widehat{\mathbf F}_f^{n+1}
=\widehat\nabla\widehat\varphi_f^{n+1},\qquad
\widehat J_f^{n+1}=\det\widehat{\mathbf F}_f^{n+1},\qquad
\widehat{\boldsymbol w}_f^{n+1}
=\frac{\widehat{\boldsymbol\eta}_f^{n+1}
-\widehat{\boldsymbol\eta}_f^{n}}{\Delta t}.
\end{aligned}
\end{equation}
Assuming that the resulting maps are invertible and preserve orientation. Their displacement increment gives the mesh velocity
\begin{equation}
\boldsymbol w_f^{n+1}
=\widehat{\boldsymbol w}_f^{n+1}
\circ(\widehat\varphi_f^{n+1})^{-1}.
\end{equation}

The fluid interface integrals use the Jacobian of the predicted surface,
\begin{equation}\label{eq:predicted-surface-jacobian}
J_s^{n+1,*}
=\left|\partial_1\Phi_s^{n+1,*}
\times\partial_2\Phi_s^{n+1,*}\right|.
\end{equation}
Noticing that the shell branch will compute its own new surface, which may differ from this prediction.
Including the old velocity gives a first-order predictor and a consistent initial mesh velocity when the wall starts with nonzero velocity. The Jacobian in \eqref{eq:predicted-surface-jacobian} belongs to the
fluid prediction, while the new shell surface is
$\varphi+\boldsymbol\eta^{n+1}$. From the second step onward, the old
fluid surface is also a prediction. Using the predictor and shell
kinematics, we therefore obtain, for $n\geq1$,
\begin{equation}\label{eq:predicted-mesh-trace}
\widehat{\boldsymbol w}_f^{n+1}\circ\varphi
=\frac{\boldsymbol\eta^n+\Delta t\dot{\boldsymbol\eta}^n
-\boldsymbol\eta^{n-1}-\Delta t\dot{\boldsymbol\eta}^{n-1}}{\Delta t}
=2\dot{\boldsymbol\eta}^n-\dot{\boldsymbol\eta}^{n-1}.
\end{equation}
The interface mesh velocity thus involves two shell time levels and
generally differs from either the old or new shell velocity. At the
first step, the fluid map starts from $\boldsymbol\eta^0$ without
prediction, so
$\widehat{\boldsymbol w}_f^1\circ\varphi=\dot{\boldsymbol\eta}^0$.

\begin{remark}[Alternative geometry prediction]
One may instead use the lagged surface
$\Phi_s^{n+1,*}=\varphi+\boldsymbol\eta^n$, which also allows independent
fluid and shell solves. We use \eqref{eq:predicted-surface} to include
the prescribed initial wall velocity in the first mesh update. Both
choices give the same scheme when the fluid domain is fixed, so the
fixed-domain proofs below apply to either choice. Their moving-domain
behavior requires a separate analysis.
\end{remark}

\subsection{Fluid subproblem on the moving domain}\label{subsec:fluid-step}
The fluid fields in the physical domain can be related to their reference representations by
\begin{equation}
\widehat{\boldsymbol u}^{n+1}=\boldsymbol u^{n+1}\circ\widehat\varphi_f^{n+1},
\qquad
\widehat p^{n+1}=p^{n+1}\circ\widehat\varphi_f^{n+1}.
\end{equation}
Both $\widehat{\boldsymbol u}^n$ and
$\widehat{\boldsymbol u}^{n+1}$ are defined on the reference domain.
Their backward difference therefore compares velocities at matching
reference points. In physical coordinates, the old velocity must first
be transported to the new domain. 
% Write this transported field as
% $\boldsymbol u^{n\to n+1}=\widehat{\boldsymbol u}^n
% \circ(\widehat\varphi_f^{n+1})^{-1}$.
The fully implicit fluid subproblem is
\begin{equation}\label{eq:fluid-strong}
\begin{aligned}
\rho_f\left[
\frac{\boldsymbol u^{n+1}-\widehat{\boldsymbol u}^n
\circ(\widehat\varphi_f^{n+1})^{-1}}{\Delta t}
+((\boldsymbol u^{n+1}-\boldsymbol w_f^{n+1})\cdot\nabla)
\boldsymbol u^{n+1}\right]&
-\nabla\cdot\boldsymbol\sigma_f(\boldsymbol u^{n+1},p^{n+1})
=\boldsymbol F_f^{n+1},\\
\nabla\cdot\boldsymbol u^{n+1}&=0
\qquad\text{in }\Omega_f^{n+1}.
\end{aligned}
\end{equation}
The reference form of the divergence constraint is
$\operatorname{tr}(\widehat\nabla\widehat{\boldsymbol u}^{n+1}
(\widehat{\mathbf F}_f^{n+1})^{-1})=0$.
On the predicted interface impose \eqref{eq:fluid-robin}, with the
reference Cauchy traction
\begin{equation}\label{eq:numerical-traction}
\widehat{\boldsymbol t}_f^{n+1}
=J_s^{n+1,*}
\left[\boldsymbol\sigma_f(\boldsymbol u^{n+1},p^{n+1})
\boldsymbol n_f^{n+1}\right]\circ\Phi_s^{n+1,*}.
\end{equation}
The other external conditions are
\begin{equation}\label{eq:external-fluid-step}
\begin{aligned}
\boldsymbol u^{n+1}&=\boldsymbol0
&&\text{on }\Sigma_f^{D,n+1},\\
\boldsymbol\sigma_f(\boldsymbol u^{n+1},p^{n+1})\boldsymbol n_f^{n+1}
&=-p_{\mathrm{in}}^{n+1}\boldsymbol n_f^{n+1}
&&\text{on }\Gamma_f^{\mathrm{in},n+1},\\
\boldsymbol\sigma_f(\boldsymbol u^{n+1},p^{n+1})\boldsymbol n_f^{n+1}
&=-p_{\mathrm{out}}^{n+1}\boldsymbol n_f^{n+1}
&&\text{on }\Gamma_f^{\mathrm{out},n+1}.
\end{aligned}
\end{equation}

To write the weak problem on the reference domain, we use the stress
\begin{equation}\label{eq:reference-stress}
\begin{aligned}
\widehat{\boldsymbol\sigma}_f^{n+1}
={}&-\widehat p^{n+1}\mathbf I
+\mu_f\left[
\widehat\nabla\widehat{\boldsymbol u}^{n+1}(\widehat{\mathbf F}_f^{n+1})^{-1}
+(\widehat{\mathbf F}_f^{n+1})^{-T}
(\widehat\nabla\widehat{\boldsymbol u}^{n+1})^T\right].
\end{aligned}
\end{equation}
Then, the weak fluid problem is: 

For all reference test functions
$\widehat{\boldsymbol{\tilde v}}$ vanishing on the reference Dirichlet boundary
and all admissible $\widehat q$, require
\begin{equation}\label{eq:fluid-weak}
\begin{aligned}
&\rho_f\int_{\widehat\Omega_f}\widehat J_f^{n+1}
\frac{\widehat{\boldsymbol u}^{n+1}-\widehat{\boldsymbol u}^{n}}{\Delta t}
\cdot\widehat{\boldsymbol{\tilde v}}\,d\widehat{\boldsymbol x}
+\rho_f\int_{\widehat\Omega_f}\widehat J_f^{n+1}
\left[\widehat\nabla\widehat{\boldsymbol u}^{n+1}
(\widehat{\mathbf F}_f^{n+1})^{-1}
(\widehat{\boldsymbol u}^{n+1}-\widehat{\boldsymbol w}_f^{n+1})\right]
\cdot\widehat{\boldsymbol{\tilde v}}\,d\widehat{\boldsymbol x}\\
&+\int_{\widehat\Omega_f}\widehat J_f^{n+1}
\widehat{\boldsymbol\sigma}_f^{n+1}(\widehat{\mathbf F}_f^{n+1})^{-T}
:\widehat\nabla\widehat{\boldsymbol{\tilde v}}\,d\widehat{\boldsymbol x}
+\int_{\widehat\Omega_f}\widehat J_f^{n+1}\widehat q\,
\operatorname{tr}\left(\widehat\nabla\widehat{\boldsymbol u}^{n+1}
(\widehat{\mathbf F}_f^{n+1})^{-1}\right)\,d\widehat{\boldsymbol x}\\
&+\gamma\int_\omega
(\widehat{\boldsymbol u}^{n+1}\circ\varphi)
\cdot(\widehat{\boldsymbol{\tilde v}}\circ\varphi)\,dx\,dy_{\mathrm{arc}}\\
&=\int_{\Omega_f^{n+1}}\boldsymbol F_f^{n+1}
\cdot\boldsymbol{\tilde v}\,d\boldsymbol x
-\int_{\Gamma_f^{\mathrm{in},n+1}}p_{\mathrm{in}}^{n+1}
\boldsymbol n_f^{n+1}\cdot\boldsymbol{\tilde v}\,dS
-\int_{\Gamma_f^{\mathrm{out},n+1}}p_{\mathrm{out}}^{n+1}
\boldsymbol n_f^{n+1}\cdot\boldsymbol{\tilde v}\,dS\\
&\quad+\int_\omega\left[
\frac{\gamma}{2}(\widehat{\boldsymbol u}^n\circ\varphi+\dot{\boldsymbol\eta}^n)
+\frac12(\widehat{\boldsymbol t}_f^n-\boldsymbol f_s^{\mathrm{ref},n})
\right]\cdot(\widehat{\boldsymbol{\tilde v}}\circ\varphi)\,dx\,dy_{\mathrm{arc}},
\end{aligned}
\end{equation}
Here
$\boldsymbol{\tilde v}=\widehat{\boldsymbol{\tilde v}}
\circ(\widehat\varphi_f^{n+1})^{-1}$.
The stress already contains pressure, and the last volume term imposes
the transformed divergence constraint. The interface integrals are over
reference parameter area because $\widehat{\boldsymbol t}_f$ already
contains the surface Jacobian.

After the fluid solve, recover $\widehat{\boldsymbol t}_f^{n+1}$ from \eqref{eq:fluid-robin}
\begin{equation}\label{eq:fluid-load-recovery}
\begin{aligned}
\widehat{\boldsymbol t}_f^{n+1}
={}&\frac{\gamma}{2}
(\widehat{\boldsymbol u}^n\circ\varphi+\dot{\boldsymbol\eta}^n)
+\frac12(\widehat{\boldsymbol t}_f^n-\boldsymbol f_s^{\mathrm{ref},n})
-\gamma\widehat{\boldsymbol u}^{n+1}\circ\varphi.
\end{aligned}
\end{equation}
For a sufficiently regular solution, the Robin condition makes
\eqref{eq:fluid-load-recovery} agree with the Cauchy traction in
\eqref{eq:numerical-traction}. In a finite element calculation, we use
the Robin relation to define the exchanged load in the chosen interface
space.

\paragraph{Lagged skew-symmetric convection formulation.}
In our numerical implementation, we linearize the fluid convection term by using a known velocity from the previous time step. Specifically, we define the convection velocity on the updated fluid domain as
\begin{equation}\label{eq:lagged-advector}
\boldsymbol a^{n+1}
=
\widehat{\boldsymbol u}^{n}
\circ(\widehat\varphi_f^{n+1})^{-1}
-\boldsymbol w_f^{n+1},
\end{equation}
where the previous fluid velocity is mapped to the current fluid domain and the mesh velocity is subtracted to account for the moving domain.

Using this known convection velocity, we replace the fully implicit convection term with the following skew-symmetric bilinear form:
\begin{equation}\label{eq:skew-convection}
\begin{aligned}
c^{n+1}(\boldsymbol z,\boldsymbol v)
=\frac{\rho_f}{2}
\int_{\Omega_f^{n+1}}
\Big[
((\boldsymbol a^{n+1}\cdot\nabla)\boldsymbol z)
\cdot\boldsymbol v
-((\boldsymbol a^{n+1}\cdot\nabla)\boldsymbol v)
\cdot\boldsymbol z
\Big]\,d\boldsymbol x.
\end{aligned}
\end{equation}
The main reasoning is that this formulation has the important property of
$c^{n+1}(\boldsymbol z,\boldsymbol z)=0,$
which holds regardless of whether $\boldsymbol a^{n+1}$ is divergence-free or satisfies any particular boundary condition. Thus, the convection term does not contribute to the energy when the test and trial functions are identical.
When the convection term is omitted, this distinction disappears, and both formulations reduce to the Stokes problem considered in the analysis below.

\paragraph{Pressure stabilization.}
In our numerical example, we use continuous piecewise-linear finite elements for both the fluid velocity and pressure ($P_1$-$P_1$) on the tetrahedral mesh. Since this equal-order finite element pair does not satisfy the inf-sup stability condition, we introduce the following pressure stabilization term:
\begin{equation}\label{eq:pressure-stabilization}
s_h^{n+1}(p_h^{n+1},q_h)
=
\delta_p\sum_{K\in\mathcal T_h^{n+1}}h_K^2
\int_K
\nabla p_h^{n+1}\cdot\nabla q_h\,d\boldsymbol x,
\qquad \delta_p>0,
\end{equation}
where $h_K$ denotes the diameter of the current mesh element $K$, and $\delta_p$ is a positive stabilization parameter. With this stabilization, the discrete continuity equation becomes
\[
(q_h,\nabla\cdot\boldsymbol u_h^{n+1})
+s_h^{n+1}(p_h^{n+1},q_h)=0.
\]
When testing the fluid equations with the new velocity and pressure, the mixed pressure terms cancel. In addition, the skew-symmetric convection term vanishes, while the pressure stabilization contributes the nonnegative term
\[
\Delta t\,s_h^{n+1}(p_h^{n+1},p_h^{n+1})\geq 0
\]
to the discrete energy balance. Since the fluid domain geometry and convection velocity are known at each time step, the resulting fluid problem is linear. We emphasize that the analysis presented below focuses on the time discretization of the Stokes problem on a fixed domain. Therefore, the effects of spatial stabilization and errors arising from the transfer of interface quantities are not included in the analysis.

\subsection{Naghdi shell subproblem on the reference parameter domain}
\label{subsec:shell-step}
The Naghdi shell problem is solved on the fixed reference domain $\omega$. We seek the displacement $\boldsymbol\eta^{n+1}$, rotation $\boldsymbol\omega^{n+1}$, and their corresponding velocities $\dot{\boldsymbol\eta}^{n+1}$ and $\dot{\boldsymbol\omega}^{n+1}$ satisfying
\begin{equation}\label{eq:discrete-shell-kinematics}
\frac{\boldsymbol\eta^{n+1}-\boldsymbol\eta^n}{\Delta t}
=\dot{\boldsymbol\eta}^{n+1},
\qquad
\frac{\boldsymbol\omega^{n+1}-\boldsymbol\omega^n}{\Delta t}
=\dot{\boldsymbol\omega}^{n+1},
\end{equation}
and, for every $(\boldsymbol v,\boldsymbol w)\in V_N(\omega)$,
\begin{equation}\label{eq:shell-step}
\begin{aligned}
&\varepsilon\rho_s\int_\omega
\frac{\dot{\boldsymbol\eta}^{n+1}-\dot{\boldsymbol\eta}^n}{\Delta t}
\cdot\boldsymbol v\,dx\,dy_{\mathrm{arc}}
+\frac{\varepsilon^3\rho_s}{12}\int_\omega
\frac{\dot{\boldsymbol\omega}^{n+1}-\dot{\boldsymbol\omega}^n}{\Delta t}
\cdot\boldsymbol w\,dx\,dy_{\mathrm{arc}}\\
&+B_m\big((\boldsymbol\eta^{n+1},\boldsymbol\omega^{n+1}),
(\boldsymbol v,\boldsymbol w)\big)
+B_b\big((\boldsymbol\eta^{n+1},\boldsymbol\omega^{n+1}),
(\boldsymbol v,\boldsymbol w)\big)+\gamma\int_\omega\dot{\boldsymbol\eta}^{n+1}\cdot\boldsymbol v
\,dx\,dy_{\mathrm{arc}}\\
&=\int_\omega\left[
\frac{\gamma}{2}(\widehat{\boldsymbol u}^n\circ\varphi+\dot{\boldsymbol\eta}^n)
-\frac12(\widehat{\boldsymbol t}_f^n-\boldsymbol f_s^{\mathrm{ref},n})
\right]\cdot\boldsymbol v\,dx\,dy_{\mathrm{arc}}.
\end{aligned}
\end{equation}
Once the shell problem is solved, the new shell interface load is recovered from the Robin condition:
\begin{equation}\label{eq:shell-load-recovery}
\begin{aligned}
\boldsymbol f_s^{\mathrm{ref},n+1}
={}&\frac{\gamma}{2}
(\widehat{\boldsymbol u}^n\circ\varphi+\dot{\boldsymbol\eta}^n)
-\frac12(\widehat{\boldsymbol t}_f^n-\boldsymbol f_s^{\mathrm{ref},n})
-\gamma\dot{\boldsymbol\eta}^{n+1}.
\end{aligned}
\end{equation}
Equation~\eqref{eq:shell-load-recovery} recovers the load that
balances the computed shell inertia and elastic forces. As it only uses the
new shell velocity and old interface data, we store this load for the next step. 

Finally, the updated shell midsurface is
\begin{equation}\label{eq:computed-shell-surface}
\Phi_s^{n+1}=\varphi+\boldsymbol\eta^{n+1}.
\end{equation}
Since the fluid domain is constructed using the predicted interface $\Phi_s^{n+1,*}$, the fluid and shell geometries may differ at the end of the time step.
From \eqref{eq:discrete-shell-kinematics}, this difference is
\begin{equation}\label{eq:geometry-gap}
\Phi_s^{n+1}-\Phi_s^{n+1,*}
=\Delta t(\dot{\boldsymbol\eta}^{n+1}-\dot{\boldsymbol\eta}^{n}).
\end{equation}
Equation~\eqref{eq:geometry-gap} gives the geometric difference
between the two at the completed step. For a smooth exact shell motion, the corresponding local geometry prediction error is $O(\Delta t^2)$.

\subsection{Initialization and time-stepping algorithm}
The initial ALE map is constructed by extending the
prescribed shell displacement $\boldsymbol\eta^0$ into
the reference fluid domain, such that
\[
\widehat{\varphi}_f^0\circ\varphi
=\varphi+\boldsymbol\eta^0.
\]
The initial reference fluid velocity is
$\widehat{\boldsymbol u}^0
=\boldsymbol u^0\circ\widehat{\varphi}_f^0$,
with the compatibility conditions given in
Section~\ref{sec:model}. If a compatible initial pressure $p^0$ is available, the initial fluid traction and shell load are
\begin{equation}\label{eq:initial-traction}
\begin{aligned}
\widehat{\boldsymbol t}_f^0
&=J_s^0\left[
\boldsymbol\sigma_f(\boldsymbol u^0,p^0)
\boldsymbol n_f^0
\right]\circ(\varphi+\boldsymbol\eta^0),\qquad
J_s^0
=\left|\partial_1(\varphi+\boldsymbol\eta^0)
\times\partial_2(\varphi+\boldsymbol\eta^0)\right|,\qquad
\boldsymbol f_s^{\mathrm{ref},0}
=-\widehat{\boldsymbol t}_f^0.
\end{aligned}
\end{equation}
If the initial traction is unavailable, it may be
obtained from a compatible initialization solve or
approximated with a controlled error. This error enters
the initial augmented error energy. Both initial loads
may be set to zero when the physical initial interface
traction vanishes.

The complete procedure is given in
Algorithm~\ref{alg:robin-robin}.

\begin{algorithm}[!ht]
\caption{Parallel Robin-Robin fluid-Naghdi shell scheme}
\label{alg:robin-robin}
\begin{algorithmic}[1]

\Require Initial fluid and shell states
$\boldsymbol u^0$, $\boldsymbol\eta^0$,
$\dot{\boldsymbol\eta}^0$, $\boldsymbol\omega^0$,
$\dot{\boldsymbol\omega}^0$;
initial interface loads
$\widehat{\boldsymbol t}_f^0$,
$\boldsymbol f_s^{\mathrm{ref},0}$;
$\Delta t>0$, $\gamma>0$.

\State \textbf{Initialization}
\State Construct the initial ALE map
$\widehat{\varphi}_f^0$.
\State Set
$\widehat{\boldsymbol u}^0
=\boldsymbol u^0\circ\widehat{\varphi}_f^0$.

\For{$n=0,\ldots,N-1$}

    \State \textbf{Step 1: Construct the Robin data.}

    \State ~ Form the fluid Robin right-hand side: $\displaystyle
    \begin{aligned}
    &\frac{\gamma}{2}
    (\widehat{\boldsymbol u}^n\circ\varphi
    +\dot{\boldsymbol\eta}^n)+\frac12
    (\widehat{\boldsymbol t}_f^n
    -\boldsymbol f_s^{\mathrm{ref},n}).
    \end{aligned}$
    % \Statex \hspace{\algorithmicindent}
    % $\displaystyle
    % \begin{aligned}
    % &\frac{\gamma}{2}
    % (\widehat{\boldsymbol u}^n\circ\varphi
    % +\dot{\boldsymbol\eta}^n)+\frac12
    % (\widehat{\boldsymbol t}_f^n
    % -\boldsymbol f_s^{\mathrm{ref},n}).
    % \end{aligned}$

    \State ~ Form the shell Robin right-hand side: $\displaystyle
    \begin{aligned}
    &\frac{\gamma}{2}
    (\widehat{\boldsymbol u}^n\circ\varphi
    +\dot{\boldsymbol\eta}^n)
    -\frac12
    (\widehat{\boldsymbol t}_f^n
    -\boldsymbol f_s^{\mathrm{ref},n}).
    \end{aligned}$
    % \Statex \hspace{\algorithmicindent}
    % $\displaystyle
    % \begin{aligned}
    % &\frac{\gamma}{2}
    % (\widehat{\boldsymbol u}^n\circ\varphi
    % +\dot{\boldsymbol\eta}^n)
    % -\frac12
    % (\widehat{\boldsymbol t}_f^n
    % -\boldsymbol f_s^{\mathrm{ref},n}).
    % \end{aligned}$

    \State \textbf{Step 2: Solve the fluid and shell subproblems concurrently.}

    \State ~ \underline{Fluid subproblem:}

    \State \hspace{\algorithmicindent}
    Predict the interface $\Phi_s^{n+1,*}$
    using \eqref{eq:predicted-surface}.

    \State \hspace{\algorithmicindent}
    Extend the predicted displacement
    using \eqref{eq:discrete-extension}.

    \State \hspace{\algorithmicindent}
    Construct the ALE map
    $\widehat{\varphi}_f^{n+1}$
    and mesh velocity $\boldsymbol w_f^{n+1}$.

    \State \hspace{\algorithmicindent}
    Solve \eqref{eq:fluid-weak} for
    $(\boldsymbol u^{n+1},p^{n+1})$.

    \State \hspace{\algorithmicindent}
    Compute the reference fluid velocity
    $\widehat{\boldsymbol u}^{n+1}
    =\boldsymbol u^{n+1}
    \circ\widehat{\varphi}_f^{n+1}$.

    \State \hspace{\algorithmicindent}
    Recover $\widehat{\boldsymbol t}_f^{n+1}$
    using \eqref{eq:fluid-load-recovery}.

    \State ~ \underline{Shell subproblem:}

    \State \hspace{\algorithmicindent}
    Solve \eqref{eq:discrete-shell-kinematics}
    and \eqref{eq:shell-step} for
    $\boldsymbol\eta^{n+1}$,
    $\boldsymbol\omega^{n+1}$,
    $\dot{\boldsymbol\eta}^{n+1}$,
    and $\dot{\boldsymbol\omega}^{n+1}$.

    \State \hspace{\algorithmicindent}
    Recover $\boldsymbol f_s^{\mathrm{ref},n+1}$
    using \eqref{eq:shell-load-recovery}.

    \State \hspace{\algorithmicindent}
    Update the shell surface
    $\Phi_s^{n+1}
    =\varphi+\boldsymbol\eta^{n+1}$.

    \State \textbf{Step 3:  Exchange the new interface velocities
    and loads after both subproblems finish. Store
    $\widehat{\boldsymbol u}^{n+1}$,
    $\dot{\boldsymbol\eta}^{n+1}$,
    $\widehat{\boldsymbol t}_f^{n+1}$,
    and $\boldsymbol f_s^{\mathrm{ref},n+1}$
    for the next time step.}
    
\EndFor
\end{algorithmic}
\end{algorithm}
% The fluid geometry is updated before the fluid equations
% are solved, but this computation can proceed concurrently
% with the shell solve. At the first time step, the predicted
% geometry uses the initial shell velocity. For a linear
% harmonic extension, the reference mesh velocity
% $\widehat{\boldsymbol w}_f^1$ is the harmonic extension
% of $\dot{\boldsymbol\eta}^0$. Its representation in
% physical coordinates,
% \[
% \boldsymbol w_f^1
% =\widehat{\boldsymbol w}_f^1
% \circ(\widehat{\varphi}_f^1)^{-1},
% \]
% is not necessarily harmonic.

% Both subproblems use only interface data from the previous
% time step and can therefore be solved independently,
% without coupling subiterations. The recovered fluid and
% shell loads are retained separately, with no additional
% interface correction. The fully implicit fluid formulation
% may require nonlinear iterations within the fluid branch,
% but these do not involve the shell problem. For the
% implemented lagged-convection formulation, one linear
% fluid solve is performed per time step.

\begin{remark}[Discrete coupling defects and scope of the analysis]
Adding the two Robin conditions gives the force imbalance
\begin{equation}\label{eq:force-defect}
\begin{aligned}
\widehat{\boldsymbol t}_f^{n+1}+\boldsymbol f_s^{\mathrm{ref},n+1}
=-\gamma\big[&
\widehat{\boldsymbol u}^{n+1}\circ\varphi
-\widehat{\boldsymbol u}^{n}\circ\varphi
+\dot{\boldsymbol\eta}^{n+1}
-\dot{\boldsymbol\eta}^{n}\big].
\end{aligned}
\end{equation}
Subtracting the two conditions gives the velocity mismatch
\begin{equation}\label{eq:velocity-defect}
\begin{aligned}
\gamma(\widehat{\boldsymbol u}^{n+1}\circ\varphi
-\dot{\boldsymbol\eta}^{n+1})
={}&-(\widehat{\boldsymbol t}_f^{n+1}-\widehat{\boldsymbol t}_f^n)
+(\boldsymbol f_s^{\mathrm{ref},n+1}
-\boldsymbol f_s^{\mathrm{ref},n}).
\end{aligned}
\end{equation}
Thus, the interface forces and velocities need not match exactly at each time step. Their differences depend on the corresponding time increments and are controlled through the Robin work identity in the fixed-domain analysis.

For the moving-domain scheme, additional effects arise from the difference between the predicted fluid interface and the computed shell surface, as well as the possible mismatch between fluid and mesh velocities. To isolate the coupling mechanism, the analysis in Section~\ref{sec:fixed-domain} considers a fixed fluid domain without convection, while retaining the full Naghdi shell dynamics. Stability analysis on a moving domain would also require accounting for changes in the fluid volume, discrete geometric conservation, and boundary energy fluxes.
\end{remark}

\section{Fixed-domain Stokes-Naghdi Shell formulation}\label{sec:fixed-domain}

We consider the fixed-domain Stokes-Naghdi shell problem for the stability and error analysis. In this setting, shell deformation affects the interface forces but does not change the fluid geometry. With convection omitted, the fully implicit fluid formulation \eqref{eq:fluid-strong} and the lagged skew-symmetric formulation \eqref{eq:skew-convection} reduce to the same Stokes problem. The analysis considers only time discretization and excludes the pressure stabilization \eqref{eq:pressure-stabilization} used in our numerical implementation.

For every $t\in[0,T]$, take
\begin{equation}\label{eq:fixed-geometry}
\Omega_f(t)=\widehat\Omega_f,\qquad
\Gamma_s(t)=\widehat\Gamma_s,\qquad
\widehat\varphi_f=\operatorname{Id},\qquad
\boldsymbol w_f=\boldsymbol0,\qquad
\widehat{\mathbf F}_f=\mathbf I_3,\qquad
\widehat J_f=1.
\end{equation}
The interface is parameterized by the reference cylindrical map in \eqref{eq:reference-surface}.
The reference surface measure is $dS_0=dx\,dy_{\mathrm{arc}}$, and
fluid traces are evaluated through $\varphi$. The fluid normal and all
integration measures stay fixed as the shell moves. Accordingly, a unit
interface Jacobian describes the frozen coupling surface, rather than
the geometry of the displaced shell.

Since the fluid reference and current coordinates now coincide, we
omit hats on velocity and pressure. We keep the hat on
$\widehat{\boldsymbol t}_f$ to identify a load on the shell parameter
domain:
\begin{equation}\label{eq:fixed-fluid-traction}
\widehat{\boldsymbol t}_f(x,y_{\mathrm{arc}},t)
=\left[\boldsymbol\sigma_f(\boldsymbol u,p)\boldsymbol n_f\right]
(\varphi(x,y_{\mathrm{arc}}),t).
\end{equation}
Here $\boldsymbol n_f$ is the fixed outward fluid normal. The shell
frame normal $\boldsymbol a_3$ keeps its definition from
Section~\ref{sec:model}.

\subsection{Function spaces and analytical assumptions}\label{subsec:fixed-spaces}
Assume that $\widehat\Omega_f$ is a bounded Lipschitz domain with the boundary
decomposition specified in Section~2. Define
\begin{equation}\label{eq:fixed-fluid-space}
V_f=\left\{\boldsymbol{\tilde v}\in H^1(\widehat\Omega_f;\mathbb R^3):
\boldsymbol{\tilde v}=\boldsymbol0\text{ on }\widehat\Sigma_f^D\right\}.
\end{equation}
For the mixed traction formulation with pressure level fixed by the boundary
conditions, take
\begin{equation}\label{eq:fixed-pressure-space}
Q_f=L^2(\widehat\Omega_f).
\end{equation}
Assume the appropriate divergence inf-sup condition: there exists
$\beta_f>0$, such that
\begin{equation}\label{eq:fixed-infsup}
\inf_{\substack{q\in Q_f\\q\ne0}}
\sup_{\substack{\boldsymbol{\tilde v}\in V_f\\\boldsymbol{\tilde v}\ne\boldsymbol0}}
\frac{\displaystyle\int_{\widehat\Omega_f}
q\nabla\cdot\boldsymbol{\tilde v}\,d\boldsymbol x}
{\|q\|_{L^2(\widehat\Omega_f)}
\|\boldsymbol{\tilde v}\|_{H^1(\widehat\Omega_f)}}\geq\beta_f.
\end{equation}
The prescribed traction fixes the pressure level in this setting. If a different boundary condition leaves a pressure nullspace, the pressure space must instead include the corresponding normalization.

We use the shell trial and test space $V_N(\omega)$ defined in \eqref{eq:shell-space}. We retain $B_m$ and $B_b$ from Section~\ref{sec:model} and assume
that they are continuous, symmetric, and nonnegative on $V_N(\omega)$.
These properties are the basis of the energy calculation.

To prove solvability in the stated $H^1$ shell space, we also use
the following shifted coercivity bound, with $c_N>0$ and $C_N\geq0$:
\begin{equation}\label{eq:fixed-shifted-coercivity}
\begin{aligned}
&B_m\big((\boldsymbol v,\boldsymbol w),(\boldsymbol v,\boldsymbol w)\big)
+B_b\big((\boldsymbol v,\boldsymbol w),(\boldsymbol v,\boldsymbol w)\big)
+C_N\left(\|\boldsymbol v\|_{L^2(\omega)}^2
+\|\boldsymbol w\|_{L^2(\omega)}^2\right)\\
&\qquad\geq c_N\left(\|\boldsymbol v\|_{H^1(\omega)}^2
+\|\boldsymbol w\|_{H^1(\omega)}^2\right)
\qquad\text{for every }(\boldsymbol v,\boldsymbol w)\in V_N(\omega).
\end{aligned}
\end{equation}
For the cylindrical material law, this bound follows from the full membrane-shear and bending strains. We give the argument because the
same bounds will later control the shell error. Set
\begin{equation}
c_s=\frac{E_s\nu_s}{1-\nu_s^2},\qquad
k_s=\min\{\mu,2(\mu+c_s)\}>0,
\end{equation}
where $\mu+c_s=E_s/[2(1-\nu_s)]>0$. Decomposition of a
$2\times2$ matrix into its trace and trace-free parts gives
\begin{equation}
c_s(\operatorname{tr}\mathbf S)^2+2\mu|\mathbf S|^2
\geq\min\{2\mu,2(\mu+c_s)\}|\mathbf S|^2.
\end{equation}
Since the reference frame is orthonormal and $\mathbf J$ is orthogonal,
\begin{equation}
\begin{aligned}
B_m\big((\boldsymbol v,\boldsymbol w),(\boldsymbol v,\boldsymbol w)\big)
&\geq\varepsilon k_s\sum_{\alpha=1}^2
\|\partial_\alpha\boldsymbol v+\boldsymbol a_\alpha\times\boldsymbol w\|_{L^2(\omega)}^2,\\
B_b\big((\boldsymbol v,\boldsymbol w),(\boldsymbol v,\boldsymbol w)\big)
&\geq\frac{\varepsilon^3 k_s}{12}
\|\nabla\boldsymbol w\|_{L^2(\omega)}^2.
\end{aligned}
\end{equation}
Moreover,
\begin{equation}
\|\nabla\boldsymbol v\|_{L^2(\omega)}^2
\leq2\sum_{\alpha=1}^2
\|\partial_\alpha\boldsymbol v+\boldsymbol a_\alpha\times\boldsymbol w\|_{L^2(\omega)}^2
+4\|\boldsymbol w\|_{L^2(\omega)}^2.
\end{equation}
The preceding estimates control the shell gradients after adding
sufficiently large positive $L^2$ terms. This proves
\eqref{eq:fixed-shifted-coercivity} for the cylindrical law. For other
constitutive tensors, we retain that inequality as an assumption.

Clamping both ends gives a stronger result for the present cylinder.
Poincar\'e's inequality in the axial direction yields $C_P>0$ such that
\begin{equation}
\|\boldsymbol v\|_{L^2(\omega)}\leq C_P\|\nabla\boldsymbol v\|_{L^2(\omega)},
\qquad
\|\boldsymbol w\|_{L^2(\omega)}\leq C_P\|\nabla\boldsymbol w\|_{L^2(\omega)}
\end{equation}
for every $(\boldsymbol v,\boldsymbol w)\in V_N(\omega)$.
The bending bound first controls $\boldsymbol w$ in $H^1$.
The membrane bound then controls $\nabla\boldsymbol v$, since the
rotation term is already bounded. Applying Poincar\'e's inequality to
$\boldsymbol v$ gives
\begin{equation}\label{eq:fixed-unshifted-coercivity}
\begin{aligned}
\|\boldsymbol v\|_{H^1(\omega)}^2
+\|\boldsymbol w\|_{H^1(\omega)}^2
\leq C_{\mathrm{el}}\Big[&
B_m\big((\boldsymbol v,\boldsymbol w),(\boldsymbol v,\boldsymbol w)\big)
+B_b\big((\boldsymbol v,\boldsymbol w),(\boldsymbol v,\boldsymbol w)\big)\Big].
\end{aligned}
\end{equation}
Here $C_{\mathrm{el}}$ depends on geometry, thickness, and material
parameters, but not on $\Delta t$. Thus the elastic energy is equivalent
to the full shell $H^1$ norm for the clamped cylinder. For a more general
law with only shifted coercivity, the elastic form may define a
seminorm; the arguments below also cover that case.

Throughout the analysis, $\rho_f,\mu_f,\rho_s$, and $\varepsilon$
are positive. To write external work as $L^2$ integrals, we take
$\boldsymbol F_f(t)\in L^2(\widehat\Omega_f;\mathbb R^3)$ and
$p_{\mathrm{in/out}}(t)\in L^2(\widehat\Gamma_f^{\mathrm{in/out}})$.
The same work terms can be understood as continuous functionals for
more general data. We state the time regularity needed for consistency
in Section~\ref{sec:consistency}.

\subsection{Continuous equations in the fixed domain}
The fluid equations are
\begin{equation}\label{eq:fixed-stokes}
\left\{\begin{aligned}
\rho_f\partial_t\boldsymbol u-\nabla\cdot\boldsymbol\sigma_f(\boldsymbol u,p)
&=\boldsymbol F_f&&\text{in }\widehat\Omega_f,\\
\nabla\cdot\boldsymbol u&=0&&\text{in }\widehat\Omega_f,\\
\boldsymbol\sigma_f(\boldsymbol u,p)&=-p\mathbf I_3+2\mu_f\mathbb D(\boldsymbol u),\\
\mathbb D(\boldsymbol u)&=\tfrac12\left(\nabla\boldsymbol u+(\nabla\boldsymbol u)^T\right).
\end{aligned}\right.
\end{equation}
The external boundary conditions are
\begin{equation}\label{eq:fixed-external-bc}
\begin{aligned}
\boldsymbol u&=\boldsymbol0&&\text{on }\widehat\Sigma_f^D,\\
\boldsymbol\sigma_f(\boldsymbol u,p)\boldsymbol n_f
&=-p_{\mathrm{in}}\boldsymbol n_f&&\text{on }\widehat\Gamma_f^{\mathrm{in}},\\
\boldsymbol\sigma_f(\boldsymbol u,p)\boldsymbol n_f
&=-p_{\mathrm{out}}\boldsymbol n_f&&\text{on }\widehat\Gamma_f^{\mathrm{out}}.
\end{aligned}
\end{equation}
The shell equation is:
\begin{equation}\label{eq:fixed-continuous-shell}
\begin{aligned}
&\varepsilon\rho_s\int_\omega\partial_t\dot{\boldsymbol\eta}\cdot\boldsymbol v
\,dx\,dy_{\mathrm{arc}}
+\frac{\varepsilon^3\rho_s}{12}\int_\omega
\partial_t\dot{\boldsymbol\omega}\cdot\boldsymbol w\,dx\,dy_{\mathrm{arc}}
+B_m\big((\boldsymbol\eta,\boldsymbol\omega),(\boldsymbol v,\boldsymbol w)\big)
+B_b\big((\boldsymbol\eta,\boldsymbol\omega),(\boldsymbol v,\boldsymbol w)\big)\\
&\qquad=\int_\omega\boldsymbol f_s^{\mathrm{ref}}\cdot\boldsymbol v
\,dx\,dy_{\mathrm{arc}},
\qquad(\boldsymbol v,\boldsymbol w)\in V_N(\omega).
\end{aligned}
\end{equation}
The coupling conditions are imposed on the fixed reference surface:
\begin{equation}\label{eq:fixed-continuous-coupling}
\boldsymbol u\circ\varphi=\dot{\boldsymbol\eta},\qquad
\widehat{\boldsymbol t}_f+\boldsymbol f_s^{\mathrm{ref}}=\boldsymbol0.
\end{equation}
The shell satisfies the clamped and periodic conditions from Section~\ref{sec:model}, together with compatible initial data.

Using the interface conditions with compatible test functions
$\boldsymbol{\tilde v}\circ\varphi=\boldsymbol v$, we obtain the coupled weak formulation
\begin{equation}\label{eq:fixed-continuous-weak}
\begin{aligned}
&\rho_f\int_{\widehat\Omega_f}\partial_t\boldsymbol u
\cdot\boldsymbol{\tilde v}\,d\boldsymbol x
+2\mu_f\int_{\widehat\Omega_f}
\mathbb D(\boldsymbol u):\mathbb D(\boldsymbol{\tilde v})\,d\boldsymbol x
-\int_{\widehat\Omega_f}p\nabla\cdot\boldsymbol{\tilde v}\,d\boldsymbol x
+\int_{\widehat\Omega_f}q\nabla\cdot\boldsymbol u\,d\boldsymbol x\\
&+\varepsilon\rho_s\int_\omega\partial_t\dot{\boldsymbol\eta}
\cdot\boldsymbol v\,dx\,dy_{\mathrm{arc}}
+\frac{\varepsilon^3\rho_s}{12}\int_\omega\partial_t\dot{\boldsymbol\omega}
\cdot\boldsymbol w\,dx\,dy_{\mathrm{arc}}
+B_m\big((\boldsymbol\eta,\boldsymbol\omega),(\boldsymbol v,\boldsymbol w)\big)
+B_b\big((\boldsymbol\eta,\boldsymbol\omega),(\boldsymbol v,\boldsymbol w)\big)\\
&=\int_{\widehat\Omega_f}\boldsymbol F_f\cdot\boldsymbol{\tilde v}\,d\boldsymbol x
-\int_{\widehat\Gamma_f^{\mathrm{in}}}p_{\mathrm{in}}
\boldsymbol n_f\cdot\boldsymbol{\tilde v}\,dS_0
-\int_{\widehat\Gamma_f^{\mathrm{out}}}p_{\mathrm{out}}
\boldsymbol n_f\cdot\boldsymbol{\tilde v}\,dS_0.
\end{aligned}
\end{equation}
Interface work cancels for these compatible tests. The prescribed
inlet and outlet tractions still contribute external work.

\subsection{Time-semidiscrete Robin-Robin splitting scheme in the fixed domain}\label{subsec:fixed-scheme}
Let $\Delta t=T/N$, $t^n=n\Delta t$, and fix $\gamma>0$.
The fields remain continuous in space. To start the recursion, take
$\boldsymbol u^0\in V_f$,
$(\boldsymbol\eta^0,\boldsymbol\omega^0)\in V_N(\omega)$,
$\dot{\boldsymbol\eta}^0,\dot{\boldsymbol\omega}^0
\in L^2(\omega;\mathbb R^3)$, and
$\widehat{\boldsymbol t}_f^0,\boldsymbol f_s^{\mathrm{ref},0}
\in L^2(\omega;\mathbb R^3)$.
Exact initialization also satisfies the continuous compatibility
conditions. At each step we solve the following two problems with
independent test functions, without imposing equality of their traces.

\paragraph{Fluid subproblem.}
Find $(\boldsymbol u^{n+1},p^{n+1})\in V_f\times Q_f$ such that for every
$(\boldsymbol{\tilde v},q)\in V_f\times Q_f$,
\begin{equation}\label{eq:fixed-discrete-fluid}
\begin{aligned}
&\rho_f\int_{\widehat\Omega_f}
\frac{\boldsymbol u^{n+1}-\boldsymbol u^n}{\Delta t}
\cdot\boldsymbol{\tilde v}\,d\boldsymbol x
+2\mu_f\int_{\widehat\Omega_f}
\mathbb D(\boldsymbol u^{n+1}):\mathbb D(\boldsymbol{\tilde v})\,d\boldsymbol x
-\int_{\widehat\Omega_f}p^{n+1}\nabla\cdot\boldsymbol{\tilde v}\,d\boldsymbol x
+\int_{\widehat\Omega_f}q\nabla\cdot\boldsymbol u^{n+1}\,d\boldsymbol x\\
&+\gamma\int_\omega(\boldsymbol u^{n+1}\circ\varphi)
\cdot(\boldsymbol{\tilde v}\circ\varphi)\,dx\,dy_{\mathrm{arc}}\\
&=\int_{\widehat\Omega_f}\boldsymbol F_f^{n+1}\cdot\boldsymbol{\tilde v}\,d\boldsymbol x
-\int_{\widehat\Gamma_f^{\mathrm{in}}}p_{\mathrm{in}}^{n+1}
\boldsymbol n_f\cdot\boldsymbol{\tilde v}\,dS_0
-\int_{\widehat\Gamma_f^{\mathrm{out}}}p_{\mathrm{out}}^{n+1}
\boldsymbol n_f\cdot\boldsymbol{\tilde v}\,dS_0\\
&\quad+\int_\omega\left[
\frac{\gamma}{2}(\boldsymbol u^n\circ\varphi+\dot{\boldsymbol\eta}^n)
+\frac12(\widehat{\boldsymbol t}_f^n-\boldsymbol f_s^{\mathrm{ref},n})
\right]\cdot(\boldsymbol{\tilde v}\circ\varphi)\,dx\,dy_{\mathrm{arc}}.
\end{aligned}
\end{equation}
The fluid load is updated by
\begin{equation}\label{eq:fixed-discrete-fluid-load}
\widehat{\boldsymbol t}_f^{n+1}
=\frac{\gamma}{2}(\boldsymbol u^n\circ\varphi+\dot{\boldsymbol\eta}^n)
+\frac12(\widehat{\boldsymbol t}_f^n-\boldsymbol f_s^{\mathrm{ref},n})
-\gamma\boldsymbol u^{n+1}\circ\varphi.
\end{equation}

\paragraph{Shell subproblem.}
Find $(\boldsymbol\eta^{n+1},\boldsymbol\omega^{n+1})$ and
$(\dot{\boldsymbol\eta}^{n+1},\dot{\boldsymbol\omega}^{n+1})$ in
$V_N(\omega)$ satisfying
\begin{equation}\label{eq:fixed-discrete-kinematics}
\frac{\boldsymbol\eta^{n+1}-\boldsymbol\eta^n}{\Delta t}
=\dot{\boldsymbol\eta}^{n+1},\qquad
\frac{\boldsymbol\omega^{n+1}-\boldsymbol\omega^n}{\Delta t}
=\dot{\boldsymbol\omega}^{n+1},
\end{equation}
and, for all $(\boldsymbol v,\boldsymbol w)\in V_N(\omega)$,
\begin{equation}\label{eq:fixed-discrete-shell}
\begin{aligned}
&\varepsilon\rho_s\int_\omega
\frac{\dot{\boldsymbol\eta}^{n+1}-\dot{\boldsymbol\eta}^n}{\Delta t}
\cdot\boldsymbol v\,dx\,dy_{\mathrm{arc}}
+\frac{\varepsilon^3\rho_s}{12}\int_\omega
\frac{\dot{\boldsymbol\omega}^{n+1}-\dot{\boldsymbol\omega}^n}{\Delta t}
\cdot\boldsymbol w\,dx\,dy_{\mathrm{arc}}\\
&+B_m\big((\boldsymbol\eta^{n+1},\boldsymbol\omega^{n+1}),
(\boldsymbol v,\boldsymbol w)\big)
+B_b\big((\boldsymbol\eta^{n+1},\boldsymbol\omega^{n+1}),
(\boldsymbol v,\boldsymbol w)\big)+\gamma\int_\omega\dot{\boldsymbol\eta}^{n+1}\cdot\boldsymbol v\,dx\,dy_{\mathrm{arc}}\\
&
=\int_\omega\left[
\frac{\gamma}{2}(\boldsymbol u^n\circ\varphi+\dot{\boldsymbol\eta}^n)
-\frac12(\widehat{\boldsymbol t}_f^n-\boldsymbol f_s^{\mathrm{ref},n})
\right]\cdot\boldsymbol v\,dx\,dy_{\mathrm{arc}}.
\end{aligned}
\end{equation}
The shell load is updated by
\begin{equation}\label{eq:fixed-discrete-shell-load}
\boldsymbol f_s^{\mathrm{ref},n+1}
=\frac{\gamma}{2}(\boldsymbol u^n\circ\varphi+\dot{\boldsymbol\eta}^n)
-\frac12(\widehat{\boldsymbol t}_f^n-\boldsymbol f_s^{\mathrm{ref},n})
-\gamma\dot{\boldsymbol\eta}^{n+1}.
\end{equation}
Both right-hand sides of \eqref{eq:fixed-discrete-fluid} and \eqref{eq:fixed-discrete-shell} use only old interface data, so the solves
can run concurrently. Each branch then recovers its own load for the next step. The trace theorem and the recovery relations ensure that the loads remain in $L^2(\omega;\mathbb R^3)$. For smooth solutions, the
recovered fluid load agrees with \eqref{eq:fixed-fluid-traction}. For weaker solutions, it is defined through the Robin condition.

\subsection{Solvability of the subproblems}
For fixed old data, the fluid subproblem has the symmetric velocity form
\begin{equation}\label{eq:fixed-fluid-coercivity}
\frac{\rho_f}{\Delta t}\|\boldsymbol{\tilde v}\|_{L^2(\widehat\Omega_f)}^2
+2\mu_f\|\mathbb D(\boldsymbol{\tilde v})\|_{L^2(\widehat\Omega_f)}^2
+\gamma\|\boldsymbol{\tilde v}\circ\varphi\|_{L^2(\omega)}^2.
\end{equation}
Korn's inequality with an $L^2$ term gives coercivity in
$H^1(\widehat\Omega_f)$ for every fixed $\Delta t>0$.
Together with \eqref{eq:fixed-infsup} and continuity of the data,
this gives a unique velocity and a pressure with the stated normalization.

For the shell, substitute \eqref{eq:fixed-discrete-kinematics} into
\eqref{eq:fixed-discrete-shell}. The bilinear part for the new velocities is
\begin{equation}\label{eq:fixed-shell-coercivity}
\begin{aligned}
&\frac{\varepsilon\rho_s}{\Delta t}
\|\dot{\boldsymbol\eta}^{n+1}\|_{L^2(\omega)}^2
+\frac{\varepsilon^3\rho_s}{12\Delta t}
\|\dot{\boldsymbol\omega}^{n+1}\|_{L^2(\omega)}^2
+\gamma\|\dot{\boldsymbol\eta}^{n+1}\|_{L^2(\omega)}^2\\
&\quad+\Delta t\,B_m\big(
(\dot{\boldsymbol\eta}^{n+1},\dot{\boldsymbol\omega}^{n+1}),
(\dot{\boldsymbol\eta}^{n+1},\dot{\boldsymbol\omega}^{n+1})\big)
+\Delta t\,B_b\big(
(\dot{\boldsymbol\eta}^{n+1},\dot{\boldsymbol\omega}^{n+1}),
(\dot{\boldsymbol\eta}^{n+1},\dot{\boldsymbol\omega}^{n+1})\big).
\end{aligned}
\end{equation}

The known displacement and rotation contribute
$B_m((\boldsymbol\eta^n,\boldsymbol\omega^n),(\boldsymbol v,\boldsymbol w))$
and $B_b((\boldsymbol\eta^n,\boldsymbol\omega^n),(\boldsymbol v,\boldsymbol w))$
to the right-hand side, with a minus sign after eliminating the kinematics.
The new-velocity form combines positive translational and rotational
mass terms with nonnegative elasticity. To see its coercivity, let
$m_{\Delta t}=\min\{\varepsilon\rho_s/\Delta t,
\varepsilon^3\rho_s/(12\Delta t)\}>0$.
If $C_N>0$, the form in \eqref{eq:fixed-shell-coercivity} controls
the left-hand side of \eqref{eq:fixed-shifted-coercivity} with coefficient
$\min\{\Delta t,m_{\Delta t}/C_N\}$. If $C_N=0$, the elastic terms alone
provide the same bound with coefficient $\Delta t$.
Thus, for every fixed $\Delta t>0$, the shell velocities are uniquely
determined. The displacement and rotation follow from
\eqref{eq:fixed-discrete-kinematics}. The coercivity constant may depend
on $\Delta t$, which is sufficient for solvability at each time step.

%%%%%%%%%%%%

\section{Discrete energy stability estimate}\label{sec:energy-stability}

We now prove stability of the fixed-domain scheme. In the continuous
problem, interface work cancels because the velocities agree and the
loads balance. The Robin-Robin splitting scheme does not enforce these conditions at the new time level. Instead, its Robin relations express the interface work through changes in the numerical interface energy and two nonnegative defect terms.
Combining this relation with the fluid and shell energy balances yields the discrete energy identity.

The mechanical energy at time $t^n$ is
\begin{equation}\label{eq:stability-mechanical-energy}
\begin{aligned}
\mathcal E^n={}&\frac{\rho_f}{2}
\|\boldsymbol u^n\|_{L^2(\widehat\Omega_f)}^2
+\frac{\varepsilon\rho_s}{2}
\|\dot{\boldsymbol\eta}^n\|_{L^2(\omega)}^2
+\frac{\varepsilon^3\rho_s}{24}
\|\dot{\boldsymbol\omega}^n\|_{L^2(\omega)}^2\\
&+\frac12B_m\big((\boldsymbol\eta^n,\boldsymbol\omega^n),
(\boldsymbol\eta^n,\boldsymbol\omega^n)\big)
+\frac12B_b\big((\boldsymbol\eta^n,\boldsymbol\omega^n),
(\boldsymbol\eta^n,\boldsymbol\omega^n)\big).
\end{aligned}
\end{equation}
The shell terms include membrane-shear and bending energy, with
the displacement-rotation coupling retained. The
rotational kinetic coefficient $\varepsilon^3\rho_s/24$ is half the
corresponding mass coefficient in the shell equation.

The Robin conditions introduce additional interface terms in the augmented energy
\begin{equation}\label{eq:stability-augmented-energy}
\begin{aligned}
\mathcal E_{\mathrm{aug}}^n={}&\mathcal E^n
+\frac{\Delta t\,\gamma}{2}
\left(\|\boldsymbol u^n\circ\varphi\|_{L^2(\omega)}^2
+\|\dot{\boldsymbol\eta}^n\|_{L^2(\omega)}^2\right)
+\frac{\Delta t}{2\gamma}
\left(\|\widehat{\boldsymbol t}_f^n\|_{L^2(\omega)}^2
+\|\boldsymbol f_s^{\mathrm{ref},n}\|_{L^2(\omega)}^2\right).
\end{aligned}
\end{equation}
These terms account for the Robin interface exchange and are distinct from the mechanical energy.

\subsection{Discrete energy identity}
The load recovery relations \eqref{eq:fixed-discrete-fluid-load} and
\eqref{eq:fixed-discrete-shell-load} give
\begin{equation}\label{eq:stability-two-robin-relations}
\begin{aligned}
\widehat{\boldsymbol t}_f^{n+1}
+\gamma\boldsymbol u^{n+1}\circ\varphi
&=\frac{\gamma}{2}(\boldsymbol u^n\circ\varphi+\dot{\boldsymbol\eta}^n)
+\frac12(\widehat{\boldsymbol t}_f^n-\boldsymbol f_s^{\mathrm{ref},n}),\\
\boldsymbol f_s^{\mathrm{ref},n+1}
+\gamma\dot{\boldsymbol\eta}^{n+1}
&=\frac{\gamma}{2}(\boldsymbol u^n\circ\varphi+\dot{\boldsymbol\eta}^n)
-\frac12(\widehat{\boldsymbol t}_f^n-\boldsymbol f_s^{\mathrm{ref},n}).
\end{aligned}
\end{equation}
Squaring and adding these relations cancels the cross terms on the right and gives, almost everywhere on $\omega$,
\begin{equation}\label{eq:stability-robin-squares}
\begin{aligned}
&|\widehat{\boldsymbol t}_f^{n+1}+\gamma\boldsymbol u^{n+1}\circ\varphi|^2
+|\boldsymbol f_s^{\mathrm{ref},n+1}+\gamma\dot{\boldsymbol\eta}^{n+1}|^2\\
&=\frac{\gamma^2}{2}
|\boldsymbol u^n\circ\varphi+\dot{\boldsymbol\eta}^n|^2
+\frac12|\widehat{\boldsymbol t}_f^n-\boldsymbol f_s^{\mathrm{ref},n}|^2\\
&=\gamma^2\big(|\boldsymbol u^n\circ\varphi|^2+|\dot{\boldsymbol\eta}^n|^2\big)
+|\widehat{\boldsymbol t}_f^n|^2+|\boldsymbol f_s^{\mathrm{ref},n}|^2
-\frac{\gamma^2}{2}|\boldsymbol u^n\circ\varphi-\dot{\boldsymbol\eta}^n|^2
-\frac12|\widehat{\boldsymbol t}_f^n+\boldsymbol f_s^{\mathrm{ref},n}|^2.
\end{aligned}
\end{equation}
Expanding the left-hand side, dividing by $2\gamma$, and integrating
over the interface yields
\begin{equation}\label{eq:stability-interface-work-identity}
\begin{aligned}
&\int_\omega\left[
\widehat{\boldsymbol t}_f^{n+1}\cdot(\boldsymbol u^{n+1}\circ\varphi)
+\boldsymbol f_s^{\mathrm{ref},n+1}\cdot\dot{\boldsymbol\eta}^{n+1}
\right]dx\,dy_{\mathrm{arc}}\\
&= \frac{\gamma}{2}\left[
\|\boldsymbol u^n\circ\varphi\|_{L^2(\omega)}^2
+\|\dot{\boldsymbol\eta}^n\|_{L^2(\omega)}^2
-\|\boldsymbol u^{n+1}\circ\varphi\|_{L^2(\omega)}^2
-\|\dot{\boldsymbol\eta}^{n+1}\|_{L^2(\omega)}^2\right]\\
&\quad+\frac{1}{2\gamma}\left[
\|\widehat{\boldsymbol t}_f^n\|_{L^2(\omega)}^2
+\|\boldsymbol f_s^{\mathrm{ref},n}\|_{L^2(\omega)}^2
-\|\widehat{\boldsymbol t}_f^{n+1}\|_{L^2(\omega)}^2
-\|\boldsymbol f_s^{\mathrm{ref},n+1}\|_{L^2(\omega)}^2\right]\\
&\quad-\frac{\gamma}{4}
\|\boldsymbol u^n\circ\varphi-\dot{\boldsymbol\eta}^n\|_{L^2(\omega)}^2
-\frac{1}{4\gamma}
\|\widehat{\boldsymbol t}_f^n+\boldsymbol f_s^{\mathrm{ref},n}\|_{L^2(\omega)}^2.
\end{aligned}
\end{equation}
The last two terms are nonpositive and measure the velocity mismatch and force imbalance at $t^n$, the time level used in the Robin data.

\begin{theorem}[Exact discrete energy identity]
Under the assumptions of Section~\ref{sec:fixed-domain},
suppose that the fixed-domain scheme admits a sufficiently
regular solution such that the computed velocities can be
used as test functions in the corresponding variational
formulations. Then, for any $\gamma>0$ and $\Delta t>0$,
the following discrete energy identity holds at each time step:
\begin{equation}\label{eq:stability-exact-energy-identity}
\begin{aligned}
&\mathcal E_{\mathrm{aug}}^{n+1}-\mathcal E_{\mathrm{aug}}^n
+\frac{\rho_f}{2}
\|\boldsymbol u^{n+1}-\boldsymbol u^n\|_{L^2(\widehat\Omega_f)}^2
+2\mu_f\Delta t\,
\|\mathbb D(\boldsymbol u^{n+1})\|_{L^2(\widehat\Omega_f)}^2\\
&+\frac{\varepsilon\rho_s}{2}
\|\dot{\boldsymbol\eta}^{n+1}-\dot{\boldsymbol\eta}^n\|_{L^2(\omega)}^2
+\frac{\varepsilon^3\rho_s}{24}
\|\dot{\boldsymbol\omega}^{n+1}-\dot{\boldsymbol\omega}^n\|_{L^2(\omega)}^2\\
&+\frac12B_m\big(
(\boldsymbol\eta^{n+1}-\boldsymbol\eta^n,\boldsymbol\omega^{n+1}-\boldsymbol\omega^n),
(\boldsymbol\eta^{n+1}-\boldsymbol\eta^n,\boldsymbol\omega^{n+1}-\boldsymbol\omega^n)\big)\\
&+\frac12B_b\big(
(\boldsymbol\eta^{n+1}-\boldsymbol\eta^n,\boldsymbol\omega^{n+1}-\boldsymbol\omega^n),
(\boldsymbol\eta^{n+1}-\boldsymbol\eta^n,\boldsymbol\omega^{n+1}-\boldsymbol\omega^n)\big)\\
&+\frac{\Delta t\,\gamma}{4}
\|\boldsymbol u^n\circ\varphi-\dot{\boldsymbol\eta}^n\|_{L^2(\omega)}^2
+\frac{\Delta t}{4\gamma}
\|\widehat{\boldsymbol t}_f^n+\boldsymbol f_s^{\mathrm{ref},n}\|_{L^2(\omega)}^2\\
&=\Delta t\int_{\widehat\Omega_f}
\boldsymbol F_f^{n+1}\cdot\boldsymbol u^{n+1}\,d\boldsymbol x
-\Delta t\int_{\widehat\Gamma_f^{\mathrm{in}}}
p_{\mathrm{in}}^{n+1}\boldsymbol n_f\cdot\boldsymbol u^{n+1}\,dS_0
-\Delta t\int_{\widehat\Gamma_f^{\mathrm{out}}}
p_{\mathrm{out}}^{n+1}\boldsymbol n_f\cdot\boldsymbol u^{n+1}\,dS_0.
\end{aligned}
\end{equation}
\end{theorem}

\begin{proof}
Testing the fluid equation with
$(\boldsymbol{\tilde v},q)=(\boldsymbol u^{n+1},p^{n+1})$
cancels the pressure terms. Using the recovered fluid load to express
the interface work gives
\begin{equation}\label{eq:stability-fluid-tested}
\begin{aligned}
&\rho_f\int_{\widehat\Omega_f}
\frac{\boldsymbol u^{n+1}-\boldsymbol u^n}{\Delta t}
\cdot\boldsymbol u^{n+1}\,d\boldsymbol x
+2\mu_f\|\mathbb D(\boldsymbol u^{n+1})\|_{L^2(\widehat\Omega_f)}^2\\
&=\int_\omega\widehat{\boldsymbol t}_f^{n+1}
\cdot(\boldsymbol u^{n+1}\circ\varphi)\,dx\,dy_{\mathrm{arc}}
+\int_{\widehat\Omega_f}\boldsymbol F_f^{n+1}
\cdot\boldsymbol u^{n+1}\,d\boldsymbol x\\
&\quad-\int_{\widehat\Gamma_f^{\mathrm{in}}}p_{\mathrm{in}}^{n+1}
\boldsymbol n_f\cdot\boldsymbol u^{n+1}\,dS_0
-\int_{\widehat\Gamma_f^{\mathrm{out}}}p_{\mathrm{out}}^{n+1}
\boldsymbol n_f\cdot\boldsymbol u^{n+1}\,dS_0.
\end{aligned}
\end{equation}
By the polarization identity
\begin{equation}\label{eq:stability-vector-polarization}
2(\boldsymbol a-\boldsymbol b)\cdot\boldsymbol a
=|\boldsymbol a|^2-|\boldsymbol b|^2+|\boldsymbol a-\boldsymbol b|^2,
\end{equation}
the first term becomes the fluid kinetic-energy difference divided by
$\Delta t$, plus the nonnegative increment
$\rho_f\|\boldsymbol u^{n+1}-\boldsymbol u^n\|_{L^2(\widehat\Omega_f)}^2/(2\Delta t)$.

For the shell, choose
$(\boldsymbol v,\boldsymbol w)
=(\dot{\boldsymbol\eta}^{n+1},\dot{\boldsymbol\omega}^{n+1})$
and use the recovered shell load. Then
\begin{equation}\label{eq:stability-shell-tested}
\begin{aligned}
&\varepsilon\rho_s\int_\omega
\frac{\dot{\boldsymbol\eta}^{n+1}-\dot{\boldsymbol\eta}^n}{\Delta t}
\cdot\dot{\boldsymbol\eta}^{n+1}\,dx\,dy_{\mathrm{arc}}
+\frac{\varepsilon^3\rho_s}{12}\int_\omega
\frac{\dot{\boldsymbol\omega}^{n+1}-\dot{\boldsymbol\omega}^n}{\Delta t}
\cdot\dot{\boldsymbol\omega}^{n+1}\,dx\,dy_{\mathrm{arc}}\\
&+B_m\big((\boldsymbol\eta^{n+1},\boldsymbol\omega^{n+1}),
(\dot{\boldsymbol\eta}^{n+1},\dot{\boldsymbol\omega}^{n+1})\big)
+B_b\big((\boldsymbol\eta^{n+1},\boldsymbol\omega^{n+1}),
(\dot{\boldsymbol\eta}^{n+1},\dot{\boldsymbol\omega}^{n+1})\big)\\
&=\int_\omega\boldsymbol f_s^{\mathrm{ref},n+1}
\cdot\dot{\boldsymbol\eta}^{n+1}\,dx\,dy_{\mathrm{arc}}.
\end{aligned}
\end{equation}
The same identity applies to both inertial terms, with kinetic-energy
coefficients $\varepsilon\rho_s/2$ and $\varepsilon^3\rho_s/24$.
For the elastic terms, the discrete kinematic relations express the new
velocities through displacement and rotation differences.

Using symmetry of $B_m$, we obtain
\begin{equation}\label{eq:stability-membrane-polarization}
\begin{aligned}
&2\Delta t\,B_m\big((\boldsymbol\eta^{n+1},\boldsymbol\omega^{n+1}),
(\dot{\boldsymbol\eta}^{n+1},\dot{\boldsymbol\omega}^{n+1})\big)\\
&=B_m\big((\boldsymbol\eta^{n+1},\boldsymbol\omega^{n+1}),
(\boldsymbol\eta^{n+1},\boldsymbol\omega^{n+1})\big)
-B_m\big((\boldsymbol\eta^n,\boldsymbol\omega^n),
(\boldsymbol\eta^n,\boldsymbol\omega^n)\big)\\
&\quad+B_m\big(
(\boldsymbol\eta^{n+1}-\boldsymbol\eta^n,\boldsymbol\omega^{n+1}-\boldsymbol\omega^n),
(\boldsymbol\eta^{n+1}-\boldsymbol\eta^n,\boldsymbol\omega^{n+1}-\boldsymbol\omega^n)\big),
\end{aligned}
\end{equation}
The same calculation for $B_b$ gives
\begin{equation}\label{eq:stability-bending-polarization}
\begin{aligned}
&2\Delta t\,B_b\big((\boldsymbol\eta^{n+1},\boldsymbol\omega^{n+1}),
(\dot{\boldsymbol\eta}^{n+1},\dot{\boldsymbol\omega}^{n+1})\big)\\
&=B_b\big((\boldsymbol\eta^{n+1},\boldsymbol\omega^{n+1}),
(\boldsymbol\eta^{n+1},\boldsymbol\omega^{n+1})\big)
-B_b\big((\boldsymbol\eta^n,\boldsymbol\omega^n),
(\boldsymbol\eta^n,\boldsymbol\omega^n)\big)\\
&\quad+B_b\big(
(\boldsymbol\eta^{n+1}-\boldsymbol\eta^n,\boldsymbol\omega^{n+1}-\boldsymbol\omega^n),
(\boldsymbol\eta^{n+1}-\boldsymbol\eta^n,\boldsymbol\omega^{n+1}-\boldsymbol\omega^n)\big).
\end{aligned}
\end{equation}
Both identities act on the displacement-rotation pair and retain all
elastic cross terms.

Multiply the tested fluid and shell equations by $\Delta t$ and add
them. Their energy differences give
$\mathcal E^{n+1}-\mathcal E^n$, while the remaining interface work
is $\Delta t$ times the left-hand side of
\eqref{eq:stability-interface-work-identity}. Substituting this
identity gives the augmented-energy difference and the two nonnegative
defect terms, proving \eqref{eq:stability-exact-energy-identity}.
\end{proof}

\begin{corollary}[Unconditional stability under homogeneous data]
Suppose that $\boldsymbol F_f^{n+1}=\boldsymbol 0$,
$p_{\mathrm{in}}^{n+1}=0$, and $p_{\mathrm{out}}^{n+1}=0$
for all time steps. Then the augmented energy satisfies
\begin{equation}\label{eq:stability-monotonicity}
\mathcal E_{\mathrm{aug}}^{n+1}
\leq \mathcal E_{\mathrm{aug}}^n
\leq \mathcal E_{\mathrm{aug}}^0.
\end{equation}
Moreover, summing the discrete energy inequality over
$n=0,\ldots,N-1$ and discarding the nonnegative kinetic
and elastic increment terms yields
\begin{equation}\label{eq:stability-cumulative-bound}
\begin{aligned}
&\mathcal E_{\mathrm{aug}}^N
+2\mu_f\Delta t\sum_{n=0}^{N-1}
\|\mathbb D(\boldsymbol u^{n+1})\|_{L^2(\widehat\Omega_f)}^2\\
&+\frac{\Delta t\,\gamma}{4}\sum_{n=0}^{N-1}
\|\boldsymbol u^n\circ\varphi-\dot{\boldsymbol\eta}^n\|_{L^2(\omega)}^2+\frac{\Delta t}{4\gamma}\sum_{n=0}^{N-1}
\|\widehat{\boldsymbol t}_f^n
+\boldsymbol f_s^{\mathrm{ref},n}\|_{L^2(\omega)}^2
\leq \mathcal E_{\mathrm{aug}}^0.
\end{aligned}
\end{equation}
Consequently, for any $\Delta t>0$ and fixed $\gamma>0$, the
scheme is stable without a time-step restriction from the added-mass
effect or the fluid-to-structure density ratio. Since
$\mathcal E^n\leq\mathcal E_{\mathrm{aug}}^n$, the mechanical energy
is bounded but need not decrease at every step. The cumulative estimate
also controls the interface defects. Their behavior under refinement is
addressed in Sections~\ref{sec:consistency} and~\ref{sec:temporal-error}.
\end{corollary}

\subsection{External work and driven stability}\label{subsec:driven-stability}
For nonzero external data, we estimate the work on the right of
\eqref{eq:stability-exact-energy-identity} using
$\boldsymbol F_f^{n+1}=\boldsymbol F_f(t^{n+1})$ and
$p_{\mathrm{in/out}}^{n+1}=p_{\mathrm{in/out}}(t^{n+1})$.
These samples must be well defined, with bounded discrete sums of
squares as $\Delta t$ decreases. Continuity in time with values in the
stated spatial $L^2$ spaces is sufficient. An $L^2$ assumption in time
alone does not justify point evaluation.
To estimate the external work, we use the Korn inequality appropriate to the
fluid test space. We consider separately the cases with and without an $L^2$
velocity term.

\paragraph{Korn inequality without an $L^2$ term.}
Suppose first that the fluid test space satisfies
\begin{equation}\label{eq:stability-first-korn}
\|\boldsymbol{\tilde v}\|_{H^1(\widehat\Omega_f)}
\leq C_K\|\mathbb D(\boldsymbol{\tilde v})\|_{L^2(\widehat\Omega_f)}.
\end{equation}
This holds, for example, when the homogeneous Dirichlet boundary
has positive surface measure. Let $C_{\mathrm{in}}$ and
$C_{\mathrm{out}}$ be trace constants such that
\begin{equation}
\begin{aligned}
\|\boldsymbol{\tilde v}\|_{L^2(\widehat\Gamma_f^{\mathrm{in}})}
&\leq C_{\mathrm{in}}\|\boldsymbol{\tilde v}\|_{H^1(\widehat\Omega_f)},\\
\|\boldsymbol{\tilde v}\|_{L^2(\widehat\Gamma_f^{\mathrm{out}})}
&\leq C_{\mathrm{out}}\|\boldsymbol{\tilde v}\|_{H^1(\widehat\Omega_f)}.
\end{aligned}
\end{equation}

The Cauchy-Schwarz and trace inequalities, followed by Young's
inequality, give
\begin{equation}\label{eq:stability-work-bound-first}
\begin{aligned}
&\left|\int_{\widehat\Omega_f}\boldsymbol F_f^{n+1}\cdot\boldsymbol u^{n+1}
\,d\boldsymbol x
-\int_{\widehat\Gamma_f^{\mathrm{in}}}p_{\mathrm{in}}^{n+1}
\boldsymbol n_f\cdot\boldsymbol u^{n+1}\,dS_0
-\int_{\widehat\Gamma_f^{\mathrm{out}}}p_{\mathrm{out}}^{n+1}
\boldsymbol n_f\cdot\boldsymbol u^{n+1}\,dS_0\right|\\
&\leq\mu_f\|\mathbb D(\boldsymbol u^{n+1})\|_{L^2(\widehat\Omega_f)}^2
+\frac{C_K^2}{4\mu_f}\left(
\|\boldsymbol F_f^{n+1}\|_{L^2(\widehat\Omega_f)}
+C_{\mathrm{in}}\|p_{\mathrm{in}}^{n+1}\|_{L^2(\widehat\Gamma_f^{\mathrm{in}})}
+C_{\mathrm{out}}\|p_{\mathrm{out}}^{n+1}\|_{L^2(\widehat\Gamma_f^{\mathrm{out}})}
\right)^2.
\end{aligned}
\end{equation}
The first term is absorbed into the viscous dissipation. Summing the
resulting inequality and discarding the remaining nonnegative terms
gives

\begin{equation}\label{eq:stability-driven-first}
\begin{aligned}
&\mathcal E_{\mathrm{aug}}^N
+\mu_f\Delta t\sum_{n=0}^{N-1}
\|\mathbb D(\boldsymbol u^{n+1})\|_{L^2(\widehat\Omega_f)}^2\\
&\leq\mathcal E_{\mathrm{aug}}^0
+\frac{C_K^2\Delta t}{4\mu_f}\sum_{n=0}^{N-1}\left(
\|\boldsymbol F_f^{n+1}\|_{L^2(\widehat\Omega_f)}
+C_{\mathrm{in}}\|p_{\mathrm{in}}^{n+1}\|_{L^2(\widehat\Gamma_f^{\mathrm{in}})}
+C_{\mathrm{out}}\|p_{\mathrm{out}}^{n+1}\|_{L^2(\widehat\Gamma_f^{\mathrm{out}})}
\right)^2.
\end{aligned}
\end{equation}
This estimate holds for every positive time step.

\paragraph{Korn inequality with an $L^2$ term.}
The tube benchmark has an elastic lateral wall and traction conditions
at both ends. Its fluid test space has no Dirichlet portion of positive
surface measure, since clamping the shell at the end rings does not
impose a fluid Dirichlet condition. We therefore use
\begin{equation}\label{eq:stability-second-korn}
\|\boldsymbol{\tilde v}\|_{H^1(\widehat\Omega_f)}
\leq C_K\left(
\|\mathbb D(\boldsymbol{\tilde v})\|_{L^2(\widehat\Omega_f)}
+\|\boldsymbol{\tilde v}\|_{L^2(\widehat\Omega_f)}\right).
\end{equation}
Applying Young's inequality to the symmetric-gradient and
velocity terms separately gives the bound
\begin{equation}\label{eq:stability-work-bound-second}
\begin{aligned}
&\mu_f\|\mathbb D(\boldsymbol u^{n+1})\|_{L^2(\widehat\Omega_f)}^2
+\frac{\rho_f}{2T}\|\boldsymbol u^{n+1}\|_{L^2(\widehat\Omega_f)}^2\\
&\quad+C_K^2\left(\frac{1}{4\mu_f}+\frac{T}{2\rho_f}\right)
\left(
\|\boldsymbol F_f^{n+1}\|_{L^2(\widehat\Omega_f)}
+C_{\mathrm{in}}\|p_{\mathrm{in}}^{n+1}\|_{L^2(\widehat\Gamma_f^{\mathrm{in}})}
+C_{\mathrm{out}}\|p_{\mathrm{out}}^{n+1}\|_{L^2(\widehat\Gamma_f^{\mathrm{out}})}
\right)^2.
\end{aligned}
\end{equation}
Since $\rho_f\|\boldsymbol u^{n+1}\|^2/(2T)
\leq\mathcal E_{\mathrm{aug}}^{n+1}/T$, the energy identity implies
\begin{equation}\label{eq:stability-driven-recurrence}
\begin{aligned}
&\left(1-\frac{\Delta t}{T}\right)\mathcal E_{\mathrm{aug}}^{n+1}
+\mu_f\Delta t\|\mathbb D(\boldsymbol u^{n+1})\|_{L^2(\widehat\Omega_f)}^2\\
&\leq\mathcal E_{\mathrm{aug}}^n
+C_K^2\Delta t\left(\frac{1}{4\mu_f}+\frac{T}{2\rho_f}\right)
\left(
\|\boldsymbol F_f^{n+1}\|_{L^2(\widehat\Omega_f)}
+C_{\mathrm{in}}\|p_{\mathrm{in}}^{n+1}\|_{L^2(\widehat\Gamma_f^{\mathrm{in}})}
+C_{\mathrm{out}}\|p_{\mathrm{out}}^{n+1}\|_{L^2(\widehat\Gamma_f^{\mathrm{out}})}
\right)^2.
\end{aligned}
\end{equation}
For $\Delta t\leq T/2$, the inequality
$-\log(1-\Delta t/T)\leq2\Delta t/T$ bounds the product of the
recurrence factors by $\exp(2)$ over $N$ steps. Iteration bounds the
maximum energy, and summation then gives
\begin{equation}\label{eq:stability-driven-second}
\begin{aligned}
&\max_{0\leq n\leq N}\mathcal E_{\mathrm{aug}}^n
+\mu_f\Delta t\sum_{n=0}^{N-1}
\|\mathbb D(\boldsymbol u^{n+1})\|_{L^2(\widehat\Omega_f)}^2\\
&\leq C\Bigg[\mathcal E_{\mathrm{aug}}^0
+C_K^2\Delta t\left(\frac{1}{4\mu_f}+\frac{T}{2\rho_f}\right)
\sum_{n=0}^{N-1}\left(
\|\boldsymbol F_f^{n+1}\|_{L^2(\widehat\Omega_f)}
+C_{\mathrm{in}}\|p_{\mathrm{in}}^{n+1}\|_{L^2(\widehat\Gamma_f^{\mathrm{in}})}
+C_{\mathrm{out}}\|p_{\mathrm{out}}^{n+1}\|_{L^2(\widehat\Gamma_f^{\mathrm{out}})}
\right)^2\Bigg],
\end{aligned}
\end{equation}
Here $C$ is independent of $\Delta t$.
One may take $C=1+2\exp(2)$, since the maximum
energy is bounded by $\exp(2)$ times the bracketed quantity and
the viscous sum by that quantity plus the maximum energy. The condition
$\Delta t\leq T/2$ applies only to this driven-data estimate and is
independent of mesh size, density ratio, and shell stiffness. The
homogeneous energy identity has no such restriction.

%%%%%%%%%%%%%%%%

\section{Consistency analysis}\label{sec:consistency}

To estimate consistency, we insert the exact solution into the discrete equations. Backward differences introduce residuals in the fluid and shell equations and in both shell kinematic relations. The Robin conditions based on previous data also introduce two interface residuals. In the following section, we will estimate each of these terms before deriving the error equations.

We keep $\gamma>0$ fixed, and evaluate
the external data at $t^{n+1}$:
$\boldsymbol F_f^{n+1}=\boldsymbol F_f(t^{n+1})$,
$p_{\mathrm{in}}^{n+1}=p_{\mathrm{in}}(t^{n+1})$, and
$p_{\mathrm{out}}^{n+1}=p_{\mathrm{out}}(t^{n+1})$. Here, the analysis
is time-semidiscrete, so no spatial projection or interpolation is involved. We denote the exact solution at time $n$ as $\boldsymbol u(t^n)$ and numerical values as $\boldsymbol u^n$. The exact solution satisfies the equations in
Section~\ref{sec:fixed-domain}, including
\begin{equation}\label{eq:cons-exact-coupling}
\boldsymbol u(t)\circ\varphi=\dot{\boldsymbol\eta}(t),
\qquad
\boldsymbol f_s^{\mathrm{ref}}(t)=-\widehat{\boldsymbol t}_f(t).
\end{equation}
For the bounds below, assume
\begin{equation}\label{eq:cons-regularity}
\begin{aligned}
&\boldsymbol u\in H^2(0,T;L^2(\widehat\Omega_f;\mathbb R^3))
\cap C([0,T];V_f),\qquad
(\boldsymbol\eta,~\boldsymbol\omega)\in H^2(0,T;V_N(\omega)),\\
&\dot{\boldsymbol\eta}, ~\dot{\boldsymbol\omega}
\in H^2(0,T;L^2(\omega;\mathbb R^3)),\qquad
\widehat{\boldsymbol t}_f\in H^1(0,T;L^2(\omega;\mathbb R^3)).
\end{aligned}
\end{equation}
Here $\dot{\boldsymbol\eta}=\partial_t\boldsymbol\eta$ and
$\dot{\boldsymbol\omega}=\partial_t\boldsymbol\omega$.
Thus the regularity in $\dot{\boldsymbol\eta},~\dot{\boldsymbol\omega}$ include third order time derivatives of displacement and rotation in the $L^2$ spaces. 
We also assume the following continuity in time so that the pressure and external data can be evaluated at the discrete time levels:
\begin{equation}\label{eq:cons-data-regularity}
\begin{aligned}
&p\in C([0,T];Q_f),\qquad
\boldsymbol F_f\in C([0,T];L^2(\widehat\Omega_f;\mathbb R^3)),\\
&p_{\mathrm{in}}\in C([0,T];L^2(\widehat\Gamma_f^{\mathrm{in}})),\qquad
p_{\mathrm{out}}\in C([0,T];L^2(\widehat\Gamma_f^{\mathrm{out}})).
\end{aligned}
\end{equation}
Together with \eqref{eq:cons-regularity}, these assumptions ensure that the weak equations, initially satisfied almost everywhere in time, can be evaluated at each $t^n$. The interface traction requires a separate regularity assumption, since the continuity of pressure in $Q_f$ does not guarantee an $L^2$ boundary traction. We estimate the consistency residuals in discrete $L^2$ time norms, for which the assumed square-integrability of the corresponding time derivatives is sufficient.

\paragraph{Interface residuals}
We define each residual as the discrete left-hand side minus the
right-hand side, evaluated at the exact solution. For the fluid Robin
condition, this gives

\begin{equation}\label{eq:cons-fluid-interface-definition}
\begin{aligned}
\boldsymbol{\mathcal R}_f^{n+1}:={}&
\widehat{\boldsymbol t}_f(t^{n+1})
+\gamma\boldsymbol u(t^{n+1})\circ\varphi
-\frac{\gamma}{2}
\big(\boldsymbol u(t^n)\circ\varphi+\dot{\boldsymbol\eta}(t^n)\big)
-\frac12\big(\widehat{\boldsymbol t}_f(t^n)
-\boldsymbol f_s^{\mathrm{ref}}(t^n)\big).
\end{aligned}
\end{equation}
Using exact velocity continuity and force balance from
\eqref{eq:cons-exact-coupling}, we simplify this expression to
\begin{equation}\label{eq:cons-fluid-interface-simplified}
\begin{aligned}
\boldsymbol{\mathcal R}_f^{n+1}
={}&\gamma\big(\dot{\boldsymbol\eta}(t^{n+1})
-\dot{\boldsymbol\eta}(t^n)\big)
+\widehat{\boldsymbol t}_f(t^{n+1})
-\widehat{\boldsymbol t}_f(t^n)
=\int_{t^n}^{t^{n+1}}
\left(\gamma\partial_t\dot{\boldsymbol\eta}(t)
+\partial_t\widehat{\boldsymbol t}_f(t)\right)dt.
\end{aligned}
\end{equation}
For the shell Robin condition, the same sign convention gives
\begin{equation}\label{eq:cons-shell-interface-definition}
\begin{aligned}
\boldsymbol{\mathcal R}_s^{n+1}:={}&
\boldsymbol f_s^{\mathrm{ref}}(t^{n+1})
+\gamma\dot{\boldsymbol\eta}(t^{n+1})
-\frac{\gamma}{2}
\big(\boldsymbol u(t^n)\circ\varphi+\dot{\boldsymbol\eta}(t^n)\big)
+\frac12\big(\widehat{\boldsymbol t}_f(t^n)
-\boldsymbol f_s^{\mathrm{ref}}(t^n)\big).
\end{aligned}
\end{equation}
It follows that
\begin{equation}\label{eq:cons-shell-interface-simplified}
\begin{aligned}
\boldsymbol{\mathcal R}_s^{n+1}
={}&\gamma\big(\dot{\boldsymbol\eta}(t^{n+1})
-\dot{\boldsymbol\eta}(t^n)\big)
-\widehat{\boldsymbol t}_f(t^{n+1})
+\widehat{\boldsymbol t}_f(t^n)
=\int_{t^n}^{t^{n+1}}
\left(\gamma\partial_t\dot{\boldsymbol\eta}(t)
-\partial_t\widehat{\boldsymbol t}_f(t)\right)dt.
\end{aligned}
\end{equation}
The two residuals contain the same acceleration term and opposite
traction increments.
Squaring and adding cancels the cross terms. Applying the Cauchy-Schwarz in time gives
\begin{equation}\label{eq:cons-interface-squares}
\begin{aligned}
\|\boldsymbol{\mathcal R}_f^{n+1}\|_{L^2(\omega)}^2
+\|\boldsymbol{\mathcal R}_s^{n+1}\|_{L^2(\omega)}^2
=~&2\gamma^2\|\dot{\boldsymbol\eta}(t^{n+1})
-\dot{\boldsymbol\eta}(t^n)\|_{L^2(\omega)}^2
+2\|\widehat{\boldsymbol t}_f(t^{n+1})
-\widehat{\boldsymbol t}_f(t^n)\|_{L^2(\omega)}^2\\
&\leq2\Delta t\int_{t^n}^{t^{n+1}}
\left[\gamma^2\|\partial_t\dot{\boldsymbol\eta}(t)\|_{L^2(\omega)}^2
+\|\partial_t\widehat{\boldsymbol t}_f(t)\|_{L^2(\omega)}^2\right]dt.
\end{aligned}
\end{equation}
Multiplication by $\Delta t$ and summation over all intervals yield
\begin{equation}\label{eq:cons-interface-discrete-l2}
\begin{aligned}
&\Delta t\sum_{n=0}^{N-1}
\left(\|\boldsymbol{\mathcal R}_f^{n+1}\|_{L^2(\omega)}^2
+\|\boldsymbol{\mathcal R}_s^{n+1}\|_{L^2(\omega)}^2\right)
\leq2\Delta t^2\int_0^T
\left[\gamma^2\|\partial_t\dot{\boldsymbol\eta}(t)\|_{L^2(\omega)}^2
+\|\partial_t\widehat{\boldsymbol t}_f(t)\|_{L^2(\omega)}^2\right]dt.
\end{aligned}
\end{equation}
% Divide by $\Delta t$, we have:
% \begin{equation}\label{eq:cons-interface-unweighted}
% \begin{aligned}
% &\sum_{n=0}^{N-1}
% \left(\|\boldsymbol{\mathcal R}_f^{n+1}\|_{L^2(\omega)}^2
% +\|\boldsymbol{\mathcal R}_s^{n+1}\|_{L^2(\omega)}^2\right)
% \leq2\Delta t\int_0^T
% \left[\gamma^2\|\partial_t\dot{\boldsymbol\eta}(t)\|_{L^2(\omega)}^2
% +\|\partial_t\widehat{\boldsymbol t}_f(t)\|_{L^2(\omega)}^2\right]dt.
% \end{aligned}
% \end{equation}
% Therefore, the error estimate of the interface residual has $O (\sqrt{\Delta t})$.

\paragraph{Momentum residuals}
The backward-Euler residuals from the momentum equations satisfying that
\begin{equation}\label{eq:cons-momentum-residual-definitions}
\begin{aligned}
\boldsymbol\tau_u^{n+1}
&=\frac{\boldsymbol u(t^{n+1})-\boldsymbol u(t^n)}{\Delta t}
-\partial_t\boldsymbol u(t^{n+1})
=-\frac1{\Delta t}\int_{t^n}^{t^{n+1}}
(t-t^n)\partial_{tt}\boldsymbol u(t)\,dt,\\
\boldsymbol\tau_{\dot\eta}^{n+1}
&=\frac{\dot{\boldsymbol\eta}(t^{n+1})-\dot{\boldsymbol\eta}(t^n)}{\Delta t}
-\partial_t\dot{\boldsymbol\eta}(t^{n+1})
=-\frac1{\Delta t}\int_{t^n}^{t^{n+1}}
(t-t^n)\partial_{tt}\dot{\boldsymbol\eta}(t)\,dt,\\
\boldsymbol\tau_{\dot\omega}^{n+1}
&=\frac{\dot{\boldsymbol\omega}(t^{n+1})-\dot{\boldsymbol\omega}(t^n)}{\Delta t}
-\partial_t\dot{\boldsymbol\omega}(t^{n+1})
=-\frac1{\Delta t}\int_{t^n}^{t^{n+1}}
(t-t^n)\partial_{tt}\dot{\boldsymbol\omega}(t)\,dt.
\end{aligned}
\end{equation}
% Integrating by parts in time gives the exact representations
% \begin{equation}\label{eq:cons-momentum-integral-representations}
% \begin{aligned}
% \boldsymbol\tau_u^{n+1}
% &=-\frac1{\Delta t}\int_{t^n}^{t^{n+1}}
% (t-t^n)\partial_{tt}\boldsymbol u(t)\,dt,\\
% \boldsymbol\tau_{\dot\eta}^{n+1}
% &=-\frac1{\Delta t}\int_{t^n}^{t^{n+1}}
% (t-t^n)\partial_{tt}\dot{\boldsymbol\eta}(t)\,dt,\\
% \boldsymbol\tau_{\dot\omega}^{n+1}
% &=-\frac1{\Delta t}\int_{t^n}^{t^{n+1}}
% (t-t^n)\partial_{tt}\dot{\boldsymbol\omega}(t)\,dt.
% \end{aligned}
% \end{equation}
Since $\int_{t^n}^{t^{n+1}}(t-t^n)^2dt=\Delta t^3/3$, we have:
\begin{equation}\label{eq:cons-momentum-local-bounds}
\begin{aligned}
\|\boldsymbol\tau_u^{n+1}\|_{L^2(\widehat\Omega_f)}^2
&\leq\frac{\Delta t}{3}\int_{t^n}^{t^{n+1}}
\|\partial_{tt}\boldsymbol u(t)\|_{L^2(\widehat\Omega_f)}^2dt,\\
\|\boldsymbol\tau_{\dot\eta}^{n+1}\|_{L^2(\omega)}^2
&\leq\frac{\Delta t}{3}\int_{t^n}^{t^{n+1}}
\|\partial_{tt}\dot{\boldsymbol\eta}(t)\|_{L^2(\omega)}^2dt,\\
\|\boldsymbol\tau_{\dot\omega}^{n+1}\|_{L^2(\omega)}^2
&\leq\frac{\Delta t}{3}\int_{t^n}^{t^{n+1}}
\|\partial_{tt}\dot{\boldsymbol\omega}(t)\|_{L^2(\omega)}^2dt.
\end{aligned}
\end{equation}
Multiplying by $\Delta t$ and summing gives
\begin{equation}\label{eq:cons-momentum-global-bounds}
\begin{aligned}
\Delta t\sum_{n=0}^{N-1}\|\boldsymbol\tau_u^{n+1}\|_{L^2(\widehat\Omega_f)}^2
&\leq\frac{\Delta t^2}{3}\int_0^T
\|\partial_{tt}\boldsymbol u(t)\|_{L^2(\widehat\Omega_f)}^2dt,\\
\Delta t\sum_{n=0}^{N-1}\|\boldsymbol\tau_{\dot\eta}^{n+1}\|_{L^2(\omega)}^2
&\leq\frac{\Delta t^2}{3}\int_0^T
\|\partial_{tt}\dot{\boldsymbol\eta}(t)\|_{L^2(\omega)}^2dt,\\
\Delta t\sum_{n=0}^{N-1}\|\boldsymbol\tau_{\dot\omega}^{n+1}\|_{L^2(\omega)}^2
&\leq\frac{\Delta t^2}{3}\int_0^T
\|\partial_{tt}\dot{\boldsymbol\omega}(t)\|_{L^2(\omega)}^2dt.
\end{aligned}
\end{equation}

\paragraph{Kinematic residuals in the shell energy}
The displacement and rotation updates have residuals as well. Define
\begin{equation}\label{eq:cons-kinematic-definitions}
\begin{aligned}
\boldsymbol\kappa_\eta^{n+1}
=\frac{\boldsymbol\eta(t^{n+1})-\boldsymbol\eta(t^n)}{\Delta t}
-\dot{\boldsymbol\eta}(t^{n+1}),\qquad
\boldsymbol\kappa_\omega^{n+1}
=\frac{\boldsymbol\omega(t^{n+1})-\boldsymbol\omega(t^n)}{\Delta t}
-\dot{\boldsymbol\omega}(t^{n+1}).
\end{aligned}
\end{equation}
Their exact representations are
\begin{equation}\label{eq:cons-kinematic-integrals}
\begin{aligned}
\boldsymbol\kappa_\eta^{n+1}
=-\frac1{\Delta t}\int_{t^n}^{t^{n+1}}
(t-t^n)\partial_{tt}\boldsymbol\eta(t)\,dt,\qquad
\boldsymbol\kappa_\omega^{n+1}
=-\frac1{\Delta t}\int_{t^n}^{t^{n+1}}
(t-t^n)\partial_{tt}\boldsymbol\omega(t)\,dt.
\end{aligned}
\end{equation}
In particular,
\begin{equation}\label{eq:cons-exact-discrete-kinematics}
\begin{aligned}
\frac{\boldsymbol\eta(t^{n+1})-\boldsymbol\eta(t^n)}{\Delta t}
=\dot{\boldsymbol\eta}(t^{n+1})+\boldsymbol\kappa_\eta^{n+1},\qquad
\frac{\boldsymbol\omega(t^{n+1})-\boldsymbol\omega(t^n)}{\Delta t}
=\dot{\boldsymbol\omega}(t^{n+1})+\boldsymbol\kappa_\omega^{n+1}.
\end{aligned}
\end{equation}
The exact kinematics relate velocities to time derivatives, whereas the scheme uses backward differences. This is why the two residuals in \eqref{eq:cons-exact-discrete-kinematics} remain even for the exact
solution.

Cauchy-Schwarz holds for each symmetric nonnegative elastic form,
including when it defines only a seminorm. Applying it to
\eqref{eq:cons-kinematic-integrals}, followed by Cauchy--Schwarz in
time, gives
\begin{equation}\label{eq:cons-kinematic-membrane-local}
\begin{aligned}
&B_m\big((\boldsymbol\kappa_\eta^{n+1},\boldsymbol\kappa_\omega^{n+1}),
(\boldsymbol\kappa_\eta^{n+1},\boldsymbol\kappa_\omega^{n+1})\big)
\leq\frac{\Delta t}{3}\int_{t^n}^{t^{n+1}}
B_m\big((\partial_{tt}\boldsymbol\eta(t),\partial_{tt}\boldsymbol\omega(t)),
(\partial_{tt}\boldsymbol\eta(t),\partial_{tt}\boldsymbol\omega(t))\big)dt,
\end{aligned}
\end{equation}
\begin{equation}\label{eq:cons-kinematic-bending-local}
\begin{aligned}
&B_b\big((\boldsymbol\kappa_\eta^{n+1},\boldsymbol\kappa_\omega^{n+1}),
(\boldsymbol\kappa_\eta^{n+1},\boldsymbol\kappa_\omega^{n+1})\big)
\leq\frac{\Delta t}{3}\int_{t^n}^{t^{n+1}}
B_b\big((\partial_{tt}\boldsymbol\eta(t),\partial_{tt}\boldsymbol\omega(t)),
(\partial_{tt}\boldsymbol\eta(t),\partial_{tt}\boldsymbol\omega(t))\big)dt.
\end{aligned}
\end{equation}
Consequently,
\begin{equation}\label{eq:cons-kinematic-membrane-global}
\begin{aligned}
&\Delta t\sum_{n=0}^{N-1}
B_m\big((\boldsymbol\kappa_\eta^{n+1},\boldsymbol\kappa_\omega^{n+1}),
(\boldsymbol\kappa_\eta^{n+1},\boldsymbol\kappa_\omega^{n+1})\big)
\leq\frac{\Delta t^2}{3}\int_0^T
B_m\big((\partial_{tt}\boldsymbol\eta(t),\partial_{tt}\boldsymbol\omega(t)),
(\partial_{tt}\boldsymbol\eta(t),\partial_{tt}\boldsymbol\omega(t))\big)dt,
\end{aligned}
\end{equation}
and 
\begin{equation}\label{eq:cons-kinematic-bending-global}
\begin{aligned}
&\Delta t\sum_{n=0}^{N-1}
B_b\big((\boldsymbol\kappa_\eta^{n+1},\boldsymbol\kappa_\omega^{n+1}),
(\boldsymbol\kappa_\eta^{n+1},\boldsymbol\kappa_\omega^{n+1})\big)
\leq\frac{\Delta t^2}{3}\int_0^T
B_b\big((\partial_{tt}\boldsymbol\eta(t),\partial_{tt}\boldsymbol\omega(t)),
(\partial_{tt}\boldsymbol\eta(t),\partial_{tt}\boldsymbol\omega(t))\big)dt.
\end{aligned}
\end{equation}

\paragraph{Residual equations and first-order consistency}
We can now write the exact solution in the same form as the numerical
scheme, with the residuals added. For every
$(\boldsymbol{\tilde v},q)\in V_f\times Q_f$, we obtain
\begin{equation}\label{eq:cons-perturbed-fluid}
\begin{aligned}
&\rho_f\int_{\widehat\Omega_f}
\frac{\boldsymbol u(t^{n+1})-\boldsymbol u(t^n)}{\Delta t}
\cdot\boldsymbol{\tilde v}\,d\boldsymbol x
+2\mu_f\int_{\widehat\Omega_f}
\mathbb D(\boldsymbol u(t^{n+1})):\mathbb D(\boldsymbol{\tilde v})\,d\boldsymbol x\\
&-\int_{\widehat\Omega_f}p(t^{n+1})\nabla\cdot\boldsymbol{\tilde v}\,d\boldsymbol x
+\int_{\widehat\Omega_f}q\nabla\cdot\boldsymbol u(t^{n+1})\,d\boldsymbol x
+\gamma\int_\omega(\boldsymbol u(t^{n+1})\circ\varphi)
\cdot(\boldsymbol{\tilde v}\circ\varphi)\,dx\,dy_{\mathrm{arc}}\\
&=\int_{\widehat\Omega_f}\boldsymbol F_f(t^{n+1})\cdot\boldsymbol{\tilde v}
\,d\boldsymbol x
-\int_{\widehat\Gamma_f^{\mathrm{in}}}p_{\mathrm{in}}(t^{n+1})
\boldsymbol n_f\cdot\boldsymbol{\tilde v}\,dS_0
-\int_{\widehat\Gamma_f^{\mathrm{out}}}p_{\mathrm{out}}(t^{n+1})
\boldsymbol n_f\cdot\boldsymbol{\tilde v}\,dS_0\\
&\quad+\int_\omega\left[
\frac{\gamma}{2}(\boldsymbol u(t^n)\circ\varphi+\dot{\boldsymbol\eta}(t^n))
+\frac12(\widehat{\boldsymbol t}_f(t^n)-\boldsymbol f_s^{\mathrm{ref}}(t^n))
\right]\cdot(\boldsymbol{\tilde v}\circ\varphi)\,dx\,dy_{\mathrm{arc}}\\
&\quad+\rho_f\int_{\widehat\Omega_f}\boldsymbol\tau_u^{n+1}
\cdot\boldsymbol{\tilde v}\,d\boldsymbol x
+\int_\omega\boldsymbol{\mathcal R}_f^{n+1}
\cdot(\boldsymbol{\tilde v}\circ\varphi)\,dx\,dy_{\mathrm{arc}}.
\end{aligned}
\end{equation}
The divergence equation has no residual because the exact velocity
is divergence-free at $t^{n+1}$.

For all $(\boldsymbol v,\boldsymbol w)\in V_N(\omega)$, the exact shell
solution satisfies
\begin{equation}\label{eq:cons-perturbed-shell}
\begin{aligned}
&\varepsilon\rho_s\int_\omega
\frac{\dot{\boldsymbol\eta}(t^{n+1})-\dot{\boldsymbol\eta}(t^n)}{\Delta t}
\cdot\boldsymbol v\,dx\,dy_{\mathrm{arc}}
+\frac{\varepsilon^3\rho_s}{12}\int_\omega
\frac{\dot{\boldsymbol\omega}(t^{n+1})-\dot{\boldsymbol\omega}(t^n)}{\Delta t}
\cdot\boldsymbol w\,dx\,dy_{\mathrm{arc}}\\
&+B_m\big((\boldsymbol\eta(t^{n+1}),\boldsymbol\omega(t^{n+1})),
(\boldsymbol v,\boldsymbol w)\big)
+B_b\big((\boldsymbol\eta(t^{n+1}),\boldsymbol\omega(t^{n+1})),
(\boldsymbol v,\boldsymbol w)\big)
+\gamma\int_\omega\dot{\boldsymbol\eta}(t^{n+1})\cdot\boldsymbol v
\,dx\,dy_{\mathrm{arc}}\\
&=\int_\omega\left[
\frac{\gamma}{2}(\boldsymbol u(t^n)\circ\varphi+\dot{\boldsymbol\eta}(t^n))
-\frac12(\widehat{\boldsymbol t}_f(t^n)-\boldsymbol f_s^{\mathrm{ref}}(t^n))
\right]\cdot\boldsymbol v\,dx\,dy_{\mathrm{arc}}\\
&\quad+\varepsilon\rho_s\int_\omega\boldsymbol\tau_{\dot\eta}^{n+1}
\cdot\boldsymbol v\,dx\,dy_{\mathrm{arc}}
+\frac{\varepsilon^3\rho_s}{12}\int_\omega\boldsymbol\tau_{\dot\omega}^{n+1}
\cdot\boldsymbol w\,dx\,dy_{\mathrm{arc}}
+\int_\omega\boldsymbol{\mathcal R}_s^{n+1}
\cdot\boldsymbol v\,dx\,dy_{\mathrm{arc}}.
\end{aligned}
\end{equation}
We note that the exact shell fields also satisfy \eqref{eq:cons-exact-discrete-kinematics}. The kinematic residuals
will enter the energy estimate when we express velocity errors through
displacement and rotation differences. They are not additional forces
in the momentum equation above.
Moreover, the elastic forms are time-independent and evaluated at the exact
new-time fields, so they introduce no further temporal residual.

\paragraph{Proposition (Temporal consistency on the fixed domain).}
Assume \eqref{eq:cons-regularity} and \eqref{eq:cons-data-regularity},
fix $\gamma>0$, and sample the external data exactly. Then the exact
Stokes--Naghdi solution satisfies the time-semidiscrete equations with
residuals \eqref{eq:cons-fluid-interface-definition},
\eqref{eq:cons-shell-interface-definition},
\eqref{eq:cons-momentum-residual-definitions}, and
\eqref{eq:cons-kinematic-definitions}. The interface and momentum
residuals are $O(\Delta t)$ in the discrete $L^2$ time norm
$\big(\Delta t\sum_{n=0}^{N-1}\|r^{n+1}\|_X^2\big)^{1/2}$.
The paired kinematic residual has the same order in each elastic
seminorm, with the squared norm replaced by the corresponding
bilinear form. These statements concern the weighted time norm;
the stated regularity does not give a uniform pointwise residual bound.

\begin{proof}
Equations~\eqref{eq:cons-perturbed-fluid},
\eqref{eq:cons-perturbed-shell}, and
\eqref{eq:cons-exact-discrete-kinematics} give the required residual
equations. Taking square roots in
\eqref{eq:cons-interface-discrete-l2},
\eqref{eq:cons-momentum-global-bounds},
\eqref{eq:cons-kinematic-membrane-global}, and
\eqref{eq:cons-kinematic-bending-global} proves the stated orders.
Their regularity integrals are finite by assumption, and the constants
are independent of $\Delta t$ for fixed $\gamma$.
\end{proof}

With stronger regularity, Taylor expansion also gives the leading
interface terms

\begin{equation}\label{eq:cons-interface-leading-terms}
\begin{aligned}
\boldsymbol{\mathcal R}_f^{n+1}
&=\Delta t\left[\gamma\partial_t\dot{\boldsymbol\eta}(t^{n+1})
+\partial_t\widehat{\boldsymbol t}_f(t^{n+1})\right]+O(\Delta t^2),\\
\boldsymbol{\mathcal R}_s^{n+1}
&=\Delta t\left[\gamma\partial_t\dot{\boldsymbol\eta}(t^{n+1})
-\partial_t\widehat{\boldsymbol t}_f(t^{n+1})\right]+O(\Delta t^2).
\end{aligned}
\end{equation}

The $O(\Delta t^2)$ remainders here require bounded second time
derivatives of $\dot{\boldsymbol\eta}$ and
$\widehat{\boldsymbol t}_f$ in the interface norm. These expansions
explain the leading terms; the proposition above uses the weaker
integral estimates.

\section{A priori temporal error estimate}\label{sec:temporal-error}
We define the error variables from subtracting numerical solution from the exact solution:
\begin{equation}\label{eq:error-definitions}
\begin{aligned}
\boldsymbol e_u^n&=\boldsymbol u(t^n)-\boldsymbol u^n,
& e_p^n&=p(t^n)-p^n,
&\boldsymbol e_\eta^n&=\boldsymbol\eta(t^n)-\boldsymbol\eta^n,
&\boldsymbol e_{\dot\eta}^n&=\dot{\boldsymbol\eta}(t^n)-\dot{\boldsymbol\eta}^n,\\
\boldsymbol e_\omega^n&=\boldsymbol\omega(t^n)-\boldsymbol\omega^n,
&\boldsymbol e_{\dot\omega}^n&=\dot{\boldsymbol\omega}(t^n)-\dot{\boldsymbol\omega}^n,
&\widehat{\boldsymbol e}_{t_f}^n&=\widehat{\boldsymbol t}_f(t^n)-\widehat{\boldsymbol t}_f^n,
&\boldsymbol e_{f_s}^n&=\boldsymbol f_s^{\mathrm{ref}}(t^n)-\boldsymbol f_s^{\mathrm{ref},n}.
\end{aligned}
\end{equation}
Here, we also include separate errors  $\widehat{\boldsymbol e}_{t_f}^n$ and $\boldsymbol e_{f_s}^n$ for the two recovered interface loads (\eqref{eq:fixed-discrete-fluid-load} and \eqref{eq:fixed-discrete-shell-load}), as our numerical scheme does not impose exact force balance.

The mechanical error energy is
\begin{equation}\label{eq:error-mechanical-energy}
\begin{aligned}
\mathcal E_{\mathrm{err}}^n={}&
\frac{\rho_f}{2}\|\boldsymbol e_u^n\|_{L^2(\widehat\Omega_f)}^2
+\frac{\varepsilon\rho_s}{2}\|\boldsymbol e_{\dot\eta}^n\|_{L^2(\omega)}^2
+\frac{\varepsilon^3\rho_s}{24}\|\boldsymbol e_{\dot\omega}^n\|_{L^2(\omega)}^2\\
&+\frac12B_m\big((\boldsymbol e_\eta^n,\boldsymbol e_\omega^n),
(\boldsymbol e_\eta^n,\boldsymbol e_\omega^n)\big)
+\frac12B_b\big((\boldsymbol e_\eta^n,\boldsymbol e_\omega^n),
(\boldsymbol e_\eta^n,\boldsymbol e_\omega^n)\big).
\end{aligned}
\end{equation}
The augmented error energy is
\begin{equation}\label{eq:error-augmented-energy}
\begin{aligned}
\mathcal E_{\mathrm{err,aug}}^n={}&\mathcal E_{\mathrm{err}}^n
+\frac{\Delta t\,\gamma}{2}\left(
\|\boldsymbol e_u^n\circ\varphi\|_{L^2(\omega)}^2
+\|\boldsymbol e_{\dot\eta}^n\|_{L^2(\omega)}^2\right)
+\frac{\Delta t}{2\gamma}\left(
\|\widehat{\boldsymbol e}_{t_f}^n\|_{L^2(\omega)}^2
+\|\boldsymbol e_{f_s}^n\|_{L^2(\omega)}^2\right).
\end{aligned}
\end{equation}
For exact initialization, we take
\begin{equation}\label{eq:fixed-exact-initialization}
\begin{aligned}
\boldsymbol u^0&=\boldsymbol u(0),\qquad
\boldsymbol\eta^0=\boldsymbol\eta(0),\qquad
\dot{\boldsymbol\eta}^0=\dot{\boldsymbol\eta}(0),\qquad
\boldsymbol\omega^0=\boldsymbol\omega(0),\qquad
\dot{\boldsymbol\omega}^0=\dot{\boldsymbol\omega}(0),\\
\widehat{\boldsymbol t}_f^0&=
\left[\boldsymbol\sigma_f(\boldsymbol u(0),p(0))\boldsymbol n_f\right]\circ\varphi,
\qquad
\boldsymbol f_s^{\mathrm{ref},0}=-\widehat{\boldsymbol t}_f^0.
\end{aligned}
\end{equation}
This requires a compatible initial pressure and sufficiently regular
traction. The initial augmented error then vanishes. If the initial fields or loads are approximated, their errors will enter the initial augmented error energy.

\paragraph{Error equations}\label{subsec:error-equations}
Subtract \eqref{eq:fixed-discrete-fluid} and
\eqref{eq:fixed-discrete-shell} from
\eqref{eq:cons-perturbed-fluid} and
\eqref{eq:cons-perturbed-shell}, respectively.
Using the corresponding Robin load relations (\eqref{eq:fixed-discrete-fluid-load}, \eqref{eq:fixed-discrete-shell-load}, \eqref{eq:cons-fluid-interface-definition} and \eqref{eq:cons-shell-interface-definition}), we express
the interface terms through the new-time load errors. We can derive the error equation for the fluid subproblem:
\begin{equation}\label{eq:error-fluid-equation}
\begin{aligned}
&\rho_f\int_{\widehat\Omega_f}
\frac{\boldsymbol e_u^{n+1}-\boldsymbol e_u^n}{\Delta t}
\cdot\boldsymbol{\tilde v}\,d\boldsymbol x
+2\mu_f\int_{\widehat\Omega_f}
\mathbb D(\boldsymbol e_u^{n+1}):\mathbb D(\boldsymbol{\tilde v})\,d\boldsymbol x
-\int_{\widehat\Omega_f}e_p^{n+1}\nabla\cdot\boldsymbol{\tilde v}\,d\boldsymbol x
+\int_{\widehat\Omega_f}q\nabla\cdot\boldsymbol e_u^{n+1}\,d\boldsymbol x\\
&=\int_\omega\widehat{\boldsymbol e}_{t_f}^{n+1}
\cdot(\boldsymbol{\tilde v}\circ\varphi)\,dx\,dy_{\mathrm{arc}}
+\rho_f\int_{\widehat\Omega_f}\boldsymbol\tau_u^{n+1}
\cdot\boldsymbol{\tilde v}\,d\boldsymbol x,
\qquad(\boldsymbol{\tilde v},q)\in V_f\times Q_f.
\end{aligned}
\end{equation}
The shell error equation is
\begin{equation}\label{eq:error-shell-equation}
\begin{aligned}
&\varepsilon\rho_s\int_\omega
\frac{\boldsymbol e_{\dot\eta}^{n+1}-\boldsymbol e_{\dot\eta}^n}{\Delta t}
\cdot\boldsymbol v\,dx\,dy_{\mathrm{arc}}
+\frac{\varepsilon^3\rho_s}{12}\int_\omega
\frac{\boldsymbol e_{\dot\omega}^{n+1}-\boldsymbol e_{\dot\omega}^n}{\Delta t}
\cdot\boldsymbol w\,dx\,dy_{\mathrm{arc}}\\
&+B_m\big((\boldsymbol e_\eta^{n+1},\boldsymbol e_\omega^{n+1}),
(\boldsymbol v,\boldsymbol w)\big)
+B_b\big((\boldsymbol e_\eta^{n+1},\boldsymbol e_\omega^{n+1}),
(\boldsymbol v,\boldsymbol w)\big)\\
&=\int_\omega\boldsymbol e_{f_s}^{n+1}\cdot\boldsymbol v\,dx\,dy_{\mathrm{arc}}
+\varepsilon\rho_s\int_\omega\boldsymbol\tau_{\dot\eta}^{n+1}
\cdot\boldsymbol v\,dx\,dy_{\mathrm{arc}}
+\frac{\varepsilon^3\rho_s}{12}\int_\omega\boldsymbol\tau_{\dot\omega}^{n+1}
\cdot\boldsymbol w\,dx\,dy_{\mathrm{arc}},
\qquad(\boldsymbol v,\boldsymbol w)\in V_N(\omega).
\end{aligned}
\end{equation}
And the two interface error relations:
\begin{equation}\label{eq:error-fluid-robin}
\begin{aligned}
\widehat{\boldsymbol e}_{t_f}^{n+1}
+\gamma\boldsymbol e_u^{n+1}\circ\varphi
={}&\frac{\gamma}{2}(\boldsymbol e_u^n\circ\varphi+\boldsymbol e_{\dot\eta}^n)
+\frac12(\widehat{\boldsymbol e}_{t_f}^n-\boldsymbol e_{f_s}^n)
+\boldsymbol{\mathcal R}_f^{n+1},
\end{aligned}
\end{equation}
\begin{equation}\label{eq:error-shell-robin}
\begin{aligned}
\boldsymbol e_{f_s}^{n+1}+\gamma\boldsymbol e_{\dot\eta}^{n+1}
={}&\frac{\gamma}{2}(\boldsymbol e_u^n\circ\varphi+\boldsymbol e_{\dot\eta}^n)
-\frac12(\widehat{\boldsymbol e}_{t_f}^n-\boldsymbol e_{f_s}^n)
+\boldsymbol{\mathcal R}_s^{n+1}.
\end{aligned}
\end{equation}
And the kinematic error relations:
\begin{equation}\label{eq:error-kinematic-relations}
\begin{aligned}
\frac{\boldsymbol e_\eta^{n+1}-\boldsymbol e_\eta^n}{\Delta t}
=\boldsymbol e_{\dot\eta}^{n+1}+\boldsymbol\kappa_\eta^{n+1},\qquad
\frac{\boldsymbol e_\omega^{n+1}-\boldsymbol e_\omega^n}{\Delta t}
=\boldsymbol e_{\dot\omega}^{n+1}+\boldsymbol\kappa_\omega^{n+1}.
\end{aligned}
\end{equation}
We note that the interface residuals already enter through \eqref{eq:error-fluid-robin} and \eqref{eq:error-shell-robin}, so they did not appear in the momentum equations.

\paragraph{Mechanical error energy balance}
Test the fluid error equation with
$(\boldsymbol{\tilde v},q)=(\boldsymbol e_u^{n+1},e_p^{n+1})$
and the shell equation with
$(\boldsymbol v,\boldsymbol w)
=(\boldsymbol e_{\dot\eta}^{n+1},\boldsymbol e_{\dot\omega}^{n+1})$.
Pressure cancels, and the three inertial terms give kinetic-energy
differences and nonnegative increments, as in the stability proof.

For the elastic terms $B_m$ and $B_b$, since they are symmetric and bilinear, for example, we have:
\[
B_m(a,a-b)
=\frac12\left[
B_m(a,a)-B_m(b,b)+B_m(a-b,a-b)
\right].
\]
Using \eqref{eq:error-kinematic-relations} gives
\begin{equation}\label{eq:error-membrane-polarization}
\begin{aligned}
&\Delta t\,B_m\big((\boldsymbol e_\eta^{n+1},\boldsymbol e_\omega^{n+1}),
(\boldsymbol e_{\dot\eta}^{n+1},\boldsymbol e_{\dot\omega}^{n+1})\big)
=\frac12 B_m\big((\boldsymbol e_\eta^{n+1},\boldsymbol e_\omega^{n+1}),
(\boldsymbol e_\eta^{n+1},\boldsymbol e_\omega^{n+1})\big)
-\frac12 B_m\big((\boldsymbol e_\eta^n,\boldsymbol e_\omega^n),
(\boldsymbol e_\eta^n,\boldsymbol e_\omega^n)\big)\\
&\quad+\frac12 B_m\big(
(\boldsymbol e_\eta^{n+1}-\boldsymbol e_\eta^n,\boldsymbol e_\omega^{n+1}-\boldsymbol e_\omega^n),
(\boldsymbol e_\eta^{n+1}-\boldsymbol e_\eta^n,\boldsymbol e_\omega^{n+1}-\boldsymbol e_\omega^n)\big)-\Delta t\,B_m\big((\boldsymbol e_\eta^{n+1},\boldsymbol e_\omega^{n+1}),
(\boldsymbol\kappa_\eta^{n+1},\boldsymbol\kappa_\omega^{n+1})\big).
\end{aligned}
\end{equation}
Similarly, the corresponding bending identity is
\begin{equation}\label{eq:error-bending-polarization}
\begin{aligned}
&\Delta t\,B_b\big((\boldsymbol e_\eta^{n+1},\boldsymbol e_\omega^{n+1}),
(\boldsymbol e_{\dot\eta}^{n+1},\boldsymbol e_{\dot\omega}^{n+1})\big)
=\frac12 B_b\big((\boldsymbol e_\eta^{n+1},\boldsymbol e_\omega^{n+1}),
(\boldsymbol e_\eta^{n+1},\boldsymbol e_\omega^{n+1})\big)
-\frac12 B_b\big((\boldsymbol e_\eta^n,\boldsymbol e_\omega^n),
(\boldsymbol e_\eta^n,\boldsymbol e_\omega^n)\big)\\
&\quad+\frac12 B_b\big(
(\boldsymbol e_\eta^{n+1}-\boldsymbol e_\eta^n,\boldsymbol e_\omega^{n+1}-\boldsymbol e_\omega^n),
(\boldsymbol e_\eta^{n+1}-\boldsymbol e_\eta^n,\boldsymbol e_\omega^{n+1}-\boldsymbol e_\omega^n)\big)
-\Delta t\,B_b\big((\boldsymbol e_\eta^{n+1},\boldsymbol e_\omega^{n+1}),
(\boldsymbol\kappa_\eta^{n+1},\boldsymbol\kappa_\omega^{n+1})\big).
\end{aligned}
\end{equation}
Multiply the tested equations
\eqref{eq:error-fluid-equation} and
\eqref{eq:error-shell-equation} by $\Delta t$ and add them.
Apply the kinetic-energy polarization identity and
\eqref{eq:error-membrane-polarization}--%
\eqref{eq:error-bending-polarization}. After dropping
the nonnegative kinetic and elastic increments, we obtain:
\begin{equation}\label{eq:error-mechanical-balance}
\begin{aligned}
&\mathcal E_{\mathrm{err}}^{n+1}-\mathcal E_{\mathrm{err}}^n
+2\mu_f\Delta t\|\mathbb D(\boldsymbol e_u^{n+1})\|_{L^2(\widehat\Omega_f)}^2\\
&\leq\Delta t\int_\omega\left[
\widehat{\boldsymbol e}_{t_f}^{n+1}\cdot(\boldsymbol e_u^{n+1}\circ\varphi)
+\boldsymbol e_{f_s}^{n+1}\cdot\boldsymbol e_{\dot\eta}^{n+1}
\right]dx\,dy_{\mathrm{arc}}
+\rho_f\Delta t\int_{\widehat\Omega_f}
\boldsymbol\tau_u^{n+1}\cdot\boldsymbol e_u^{n+1}\,d\boldsymbol x\\
&\quad+\varepsilon\rho_s\Delta t\int_\omega
\boldsymbol\tau_{\dot\eta}^{n+1}\cdot\boldsymbol e_{\dot\eta}^{n+1}\,dx\,dy_{\mathrm{arc}}
+\frac{\varepsilon^3\rho_s\Delta t}{12}\int_\omega
\boldsymbol\tau_{\dot\omega}^{n+1}\cdot\boldsymbol e_{\dot\omega}^{n+1}\,dx\,dy_{\mathrm{arc}}\\
&\quad+\Delta t\,B_m\big((\boldsymbol e_\eta^{n+1},\boldsymbol e_\omega^{n+1}),
(\boldsymbol\kappa_\eta^{n+1},\boldsymbol\kappa_\omega^{n+1})\big)
+\Delta t\,B_b\big((\boldsymbol e_\eta^{n+1},\boldsymbol e_\omega^{n+1}),
(\boldsymbol\kappa_\eta^{n+1},\boldsymbol\kappa_\omega^{n+1})\big).
\end{aligned}
\end{equation}
The last two terms come from the exact solution's kinematic
residuals. They account for the difference between its time derivatives
and the backward differences used by the scheme.

\paragraph{Momentum and kinematic residual work}
We bound the residual work in \eqref{eq:error-mechanical-balance} using Young's inequality, with the final
time $T>0$ as a fixed scale. For the momentum terms, we have:
\begin{equation}\label{eq:error-momentum-work-bounds}
\begin{aligned}
\rho_f\left|\int_{\widehat\Omega_f}\boldsymbol\tau_u^{n+1}
\cdot\boldsymbol e_u^{n+1}\,d\boldsymbol x\right|
&\leq\frac{\rho_f}{2T}\|\boldsymbol e_u^{n+1}\|_{L^2(\widehat\Omega_f)}^2
+\frac{\rho_fT}{2}\|\boldsymbol\tau_u^{n+1}\|_{L^2(\widehat\Omega_f)}^2,\\
\varepsilon\rho_s\left|\int_\omega\boldsymbol\tau_{\dot\eta}^{n+1}
\cdot\boldsymbol e_{\dot\eta}^{n+1}\,dx\,dy_{\mathrm{arc}}\right|
&\leq\frac{\varepsilon\rho_s}{2T}\|\boldsymbol e_{\dot\eta}^{n+1}\|_{L^2(\omega)}^2
+\frac{\varepsilon\rho_sT}{2}\|\boldsymbol\tau_{\dot\eta}^{n+1}\|_{L^2(\omega)}^2,\\
\frac{\varepsilon^3\rho_s}{12}\left|\int_\omega\boldsymbol\tau_{\dot\omega}^{n+1}
\cdot\boldsymbol e_{\dot\omega}^{n+1}\,dx\,dy_{\mathrm{arc}}\right|
&\leq\frac{\varepsilon^3\rho_s}{24T}\|\boldsymbol e_{\dot\omega}^{n+1}\|_{L^2(\omega)}^2
+\frac{\varepsilon^3\rho_sT}{24}\|\boldsymbol\tau_{\dot\omega}^{n+1}\|_{L^2(\omega)}^2.
\end{aligned}
\end{equation}
Regarding the elastic terms, Cauchy-Schwarz in the shell seminorms gives
\begin{equation}\label{eq:error-kinematic-work-bounds}
\begin{aligned}
&\left|B_m\big((\boldsymbol e_\eta^{n+1},\boldsymbol e_\omega^{n+1}),
(\boldsymbol\kappa_\eta^{n+1},\boldsymbol\kappa_\omega^{n+1})\big)\right|\leq\frac{1}{2T}B_m\big((\boldsymbol e_\eta^{n+1},\boldsymbol e_\omega^{n+1}),
(\boldsymbol e_\eta^{n+1},\boldsymbol e_\omega^{n+1})\big)
+\frac{T}{2}B_m\big((\boldsymbol\kappa_\eta^{n+1},\boldsymbol\kappa_\omega^{n+1}),
(\boldsymbol\kappa_\eta^{n+1},\boldsymbol\kappa_\omega^{n+1})\big),\\
&\left|B_b\big((\boldsymbol e_\eta^{n+1},\boldsymbol e_\omega^{n+1}),
(\boldsymbol\kappa_\eta^{n+1},\boldsymbol\kappa_\omega^{n+1})\big)\right|
\leq\frac{1}{2T}B_b\big((\boldsymbol e_\eta^{n+1},\boldsymbol e_\omega^{n+1}),
(\boldsymbol e_\eta^{n+1},\boldsymbol e_\omega^{n+1})\big)
+\frac{T}{2}B_b\big((\boldsymbol\kappa_\eta^{n+1},\boldsymbol\kappa_\omega^{n+1}),
(\boldsymbol\kappa_\eta^{n+1},\boldsymbol\kappa_\omega^{n+1})\big).
\end{aligned}
\end{equation}
Together, the error-dependent terms on the right are
$\mathcal E_{\mathrm{err}}^{n+1}/T$. We bound this by
$\mathcal E_{\mathrm{err,aug}}^{n+1}/T$ when combining the estimates.

\paragraph{Interface work with residuals}
The Robin square identity \eqref{eq:stability-robin-squares} must also account for the two interface
residuals. We use
\begin{equation}\label{eq:error-young-interface}
|\boldsymbol a+\boldsymbol b|^2
\leq\left(1+\frac{\Delta t}{T}\right)|\boldsymbol a|^2
+\left(1+\frac{T}{\Delta t}\right)|\boldsymbol b|^2
\end{equation}
on the right of each Robin error relation (\eqref{eq:error-fluid-robin} and \eqref{eq:error-shell-robin}), taking $\boldsymbol b$
as its interface residual. Applying the same square-expansion and averaging calculation
as in \eqref{eq:stability-robin-squares} yields
\begin{equation}\label{eq:error-perturbed-squares}
\begin{aligned}
&\|\widehat{\boldsymbol e}_{t_f}^{n+1}
+\gamma\boldsymbol e_u^{n+1}\circ\varphi\|_{L^2(\omega)}^2
+\|\boldsymbol e_{f_s}^{n+1}
+\gamma\boldsymbol e_{\dot\eta}^{n+1}\|_{L^2(\omega)}^2\\
&\leq\left(1+\frac{\Delta t}{T}\right)\Big[
\gamma^2\big(\|\boldsymbol e_u^n\circ\varphi\|_{L^2(\omega)}^2
+\|\boldsymbol e_{\dot\eta}^n\|_{L^2(\omega)}^2\big)
+\|\widehat{\boldsymbol e}_{t_f}^n\|_{L^2(\omega)}^2
+\|\boldsymbol e_{f_s}^n\|_{L^2(\omega)}^2\\
&\hspace{3cm}-\frac{\gamma^2}{2}
\|\boldsymbol e_u^n\circ\varphi-\boldsymbol e_{\dot\eta}^n\|_{L^2(\omega)}^2
-\frac12\|\widehat{\boldsymbol e}_{t_f}^n+\boldsymbol e_{f_s}^n\|_{L^2(\omega)}^2
\Big]\\
&\quad+\left(1+\frac{T}{\Delta t}\right)
\left(\|\boldsymbol{\mathcal R}_f^{n+1}\|_{L^2(\omega)}^2
+\|\boldsymbol{\mathcal R}_s^{n+1}\|_{L^2(\omega)}^2\right).
\end{aligned}
\end{equation}
The square bracket is the sum of squares of the two old-data Robin
right-hand sides before the residuals are added. It is therefore
nonnegative, allowing us to keep the displayed defect terms. Expand the new-time squares and multiply by
$\Delta t/(2\gamma)$. The cross terms give the interface work in
\eqref{eq:error-mechanical-balance}, and the remaining terms give the
new interface error energy. The old interface energy is multiplied by
$1+\Delta t/T$. The coefficient of the residual squares is
\begin{equation}\label{eq:error-interface-source-coefficient}
\frac{T+\Delta t}{2\gamma}
\left(\|\boldsymbol{\mathcal R}_f^{n+1}\|_{L^2(\omega)}^2
+\|\boldsymbol{\mathcal R}_s^{n+1}\|_{L^2(\omega)}^2\right).
\end{equation}
This coefficient has no extra factor of $\Delta t$ multiplying
the residual squares. As a result, the interface residual contributes
one lower power of $\Delta t$ than the momentum and kinematic terms
in the final squared-error bound.

\paragraph{One-step augmented error inequality}
Combining \eqref{eq:error-mechanical-balance},
\eqref{eq:error-momentum-work-bounds},
\eqref{eq:error-kinematic-work-bounds}, and
\eqref{eq:error-perturbed-squares} gives
\begin{equation}\label{eq:error-one-step}
\begin{aligned}
&\left(1-\frac{\Delta t}{T}\right)\mathcal E_{\mathrm{err,aug}}^{n+1}
+2\mu_f\Delta t\|\mathbb D(\boldsymbol e_u^{n+1})\|_{L^2(\widehat\Omega_f)}^2\\
&+\left(1+\frac{\Delta t}{T}\right)\frac{\Delta t\,\gamma}{4}
\|\boldsymbol e_u^n\circ\varphi-\boldsymbol e_{\dot\eta}^n\|_{L^2(\omega)}^2
+\left(1+\frac{\Delta t}{T}\right)\frac{\Delta t}{4\gamma}
\|\widehat{\boldsymbol e}_{t_f}^n+\boldsymbol e_{f_s}^n\|_{L^2(\omega)}^2\\
&\leq\left(1+\frac{\Delta t}{T}\right)\mathcal E_{\mathrm{err,aug}}^n\\
&\quad+\frac{T\Delta t}{2}\Bigg[
\rho_f\|\boldsymbol\tau_u^{n+1}\|_{L^2(\widehat\Omega_f)}^2
+\varepsilon\rho_s\|\boldsymbol\tau_{\dot\eta}^{n+1}\|_{L^2(\omega)}^2
+\frac{\varepsilon^3\rho_s}{12}\|\boldsymbol\tau_{\dot\omega}^{n+1}\|_{L^2(\omega)}^2\\
&\hspace{2cm}+B_m\big((\boldsymbol\kappa_\eta^{n+1},\boldsymbol\kappa_\omega^{n+1}),
(\boldsymbol\kappa_\eta^{n+1},\boldsymbol\kappa_\omega^{n+1})\big)
+B_b\big((\boldsymbol\kappa_\eta^{n+1},\boldsymbol\kappa_\omega^{n+1}),
(\boldsymbol\kappa_\eta^{n+1},\boldsymbol\kappa_\omega^{n+1})\big)\Bigg]\\
&\quad+\frac{T+\Delta t}{2\gamma}
\left(\|\boldsymbol{\mathcal R}_f^{n+1}\|_{L^2(\omega)}^2
+\|\boldsymbol{\mathcal R}_s^{n+1}\|_{L^2(\omega)}^2\right).
\end{aligned}
\end{equation}
To obtain the above inequality in terms of augmented energy, we increased the old
mechanical energy by the factor $1+\Delta t/T$ and replaced
$\mathcal E_{\mathrm{err}}^{n+1}/T$ with
$\mathcal E_{\mathrm{err,aug}}^{n+1}/T$.
Specifically, the momentum and kinematic residual estimates contribute
$(\Delta t/T)\mathcal E_{\mathrm{err}}^{n+1}$ to the
right-hand side. Moving this term to the left and using
\eqref{eq:error-augmented-energy}, we obtain
\[
\left(1-\frac{\Delta t}{T}\right)
\mathcal E_{\mathrm{err,aug}}^{n+1}
\leq\mathcal E_{\mathrm{err,aug}}^{n+1}
-\frac{\Delta t}{T}\mathcal E_{\mathrm{err}}^{n+1}
,
\]
because
$\mathcal E_{\mathrm{err}}^{n+1}
\leq\mathcal E_{\mathrm{err,aug}}^{n+1}$.
At the old time level, \eqref{eq:error-perturbed-squares}
multiplies the interface terms by $1+\Delta t/T$.
Together with the old mechanical energy, these terms equal
\[
\mathcal E_{\mathrm{err}}^n
+\left(1+\frac{\Delta t}{T}\right)
\left(\mathcal E_{\mathrm{err,aug}}^n
-\mathcal E_{\mathrm{err}}^n\right)
=
\left(1+\frac{\Delta t}{T}\right)
\mathcal E_{\mathrm{err,aug}}^n
-\frac{\Delta t}{T}\mathcal E_{\mathrm{err}}^n
\leq
\left(1+\frac{\Delta t}{T}\right)
\mathcal E_{\mathrm{err,aug}}^n.
\]
Using these lower and upper bounds gives
\eqref{eq:error-one-step}.

%%%%%%%%%%%%%%%%%%%%

\begin{theorem}[Temporal error estimate]
\label{thm:temporal-error}
Assume the fixed-domain hypotheses of
Section~\ref{sec:fixed-domain} and the exact-solution regularity
in \eqref{eq:cons-regularity} and
\eqref{eq:cons-data-regularity}.
Let $\Delta t=T/N$ with $N\geq2$, and let $\gamma>0$
be fixed independently of $\Delta t$.
Suppose that the numerical and exact external data agree
at the time levels. Then there exists a constant $C>0$,
independent of $\Delta t$, such that
\begin{equation}\label{eq:error-final-estimate}
\begin{aligned}
&\max_{0\leq n\leq N}\mathcal E_{\mathrm{err,aug}}^n
+2\mu_f\Delta t\sum_{n=0}^{N-1}
\|\mathbb D(\boldsymbol e_u^{n+1})\|_{L^2(\widehat\Omega_f)}^2
+\frac{\Delta t\,\gamma}{4}\sum_{n=0}^{N-1}
\|\boldsymbol e_u^n\circ\varphi
-\boldsymbol e_{\dot\eta}^n\|_{L^2(\omega)}^2
+\frac{\Delta t}{4\gamma}\sum_{n=0}^{N-1}
\|\widehat{\boldsymbol e}_{t_f}^n
+\boldsymbol e_{f_s}^n\|_{L^2(\omega)}^2\\
&\leq C\Bigg\{
\mathcal E_{\mathrm{err,aug}}^0
+\frac{T\Delta t^2}{6}\int_0^T\Bigg[
\rho_f\|\partial_{tt}\boldsymbol u(t)\|_{L^2(\widehat\Omega_f)}^2
+\varepsilon\rho_s
\|\partial_{tt}\dot{\boldsymbol\eta}(t)\|_{L^2(\omega)}^2
+\frac{\varepsilon^3\rho_s}{12}
\|\partial_{tt}\dot{\boldsymbol\omega}(t)\|_{L^2(\omega)}^2\\
&\qquad+B_m\big(
(\partial_{tt}\boldsymbol\eta(t),
 \partial_{tt}\boldsymbol\omega(t)),
(\partial_{tt}\boldsymbol\eta(t),
 \partial_{tt}\boldsymbol\omega(t))\big)
 +B_b\big(
(\partial_{tt}\boldsymbol\eta(t),
 \partial_{tt}\boldsymbol\omega(t)),
(\partial_{tt}\boldsymbol\eta(t),
 \partial_{tt}\boldsymbol\omega(t))\big)
\Bigg]\,dt\\
&\qquad+\frac{(T+\Delta t)\Delta t}{\gamma}
\int_0^T\left[
\gamma^2
\|\partial_t\dot{\boldsymbol\eta}(t)\|_{L^2(\omega)}^2
+\|\partial_t\widehat{\boldsymbol t}_f(t)\|_{L^2(\omega)}^2
\right]\,dt
\Bigg\}.
\end{aligned}
\end{equation}
The physical parameters and $\gamma$ are fixed.
This estimate does not imply a uniform unweighted bound
as the shell thickness or density tends to zero.
\end{theorem}

\begin{proof}
The tested error equations and residual-work estimates give
\eqref{eq:error-one-step}.
For brevity, let $\mathcal R^{n+1}$ denote the sum of the
momentum, kinematic, and interface residual terms on the
right-hand side of that inequality, excluding
$(1+\Delta t/T)\mathcal E_{\mathrm{err,aug}}^n$.
All terms in $\mathcal R^{n+1}$ are nonnegative.

We first bound the maximum augmented error energy.
Since $N\geq2$, we have $\Delta t/T\leq1/2$.
Dropping the nonnegative viscous and interface terms in
\eqref{eq:error-one-step} and dividing by $1-\Delta t/T$
gives
\[
\mathcal E_{\mathrm{err,aug}}^{n+1}
\leq
\frac{1+\Delta t/T}{1-\Delta t/T}
\mathcal E_{\mathrm{err,aug}}^n
+\frac{\mathcal R^{n+1}}{1-\Delta t/T}.
\]
For $0\leq\Delta t/T\leq1/2$,
\[
\log\left(
\frac{1+\Delta t/T}{1-\Delta t/T}
\right)
\leq\frac{4\Delta t}{T},
\qquad
\frac{1}{1-\Delta t/T}\leq2.
\]
Since $N\Delta t=T$, iteration over at most $N$ steps yields
\begin{equation}\label{eq:error-gronwall-residual-bound}
\max_{0\leq n\leq N}\mathcal E_{\mathrm{err,aug}}^n
\leq
2\exp(4)\left(
\mathcal E_{\mathrm{err,aug}}^0
+\sum_{n=0}^{N-1}\mathcal R^{n+1}
\right).
\end{equation}

We next bound the accumulated viscous dissipation and
interface defects. Summing \eqref{eq:error-one-step}
over $n=0,\ldots,N-1$ and replacing $1+\Delta t/T$
by $1$ in front of the nonnegative interface defects gives
\begin{align}
&2\mu_f\Delta t\sum_{n=0}^{N-1}
\|\mathbb D(\boldsymbol e_u^{n+1})\|_{L^2(\widehat\Omega_f)}^2
+\frac{\Delta t\,\gamma}{4}
\sum_{n=0}^{N-1}
\|\boldsymbol e_u^n\circ\varphi
-\boldsymbol e_{\dot\eta}^n\|_{L^2(\omega)}^2
+\frac{\Delta t}{4\gamma}
\sum_{n=0}^{N-1}
\|\widehat{\boldsymbol e}_{t_f}^n
+\boldsymbol e_{f_s}^n\|_{L^2(\omega)}^2\\
&\leq
\mathcal E_{\mathrm{err,aug}}^0
-\mathcal E_{\mathrm{err,aug}}^N
+\frac{\Delta t}{T}\sum_{n=0}^{N-1}
\left(
\mathcal E_{\mathrm{err,aug}}^{n+1}
+\mathcal E_{\mathrm{err,aug}}^n
\right)
+\sum_{n=0}^{N-1}\mathcal R^{n+1}\\
&\leq
\mathcal E_{\mathrm{err,aug}}^0
+2\max_{0\leq n\leq N}
\mathcal E_{\mathrm{err,aug}}^n
+\sum_{n=0}^{N-1}\mathcal R^{n+1}.
\end{align}
The viscous and interface sums are then bounded by the same residual
quantity, with a constant independent of $\Delta t$. Here we use nonnegativity of the final energy and
$N\Delta t=T$. Combining this inequality with
\eqref{eq:error-gronwall-residual-bound} bounds the maximum
energy and all accumulated terms by
\[
\bigl(1+6\exp(4)\bigr)
\left(
\mathcal E_{\mathrm{err,aug}}^0
+\sum_{n=0}^{N-1}\mathcal R^{n+1}
\right).
\]
It remains to estimate the summed residual terms.
For the momentum and kinematic residuals, apply
\eqref{eq:cons-momentum-global-bounds},
\eqref{eq:cons-kinematic-membrane-global}, and
\eqref{eq:cons-kinematic-bending-global}.
Their coefficient $T\Delta t/2$ in
\eqref{eq:error-one-step} gives the factor
$T\Delta t^2/6$ in \eqref{eq:error-final-estimate}.

For the interface residuals, divide
\eqref{eq:cons-interface-discrete-l2} by $\Delta t$
to obtain
\[
\sum_{n=0}^{N-1}
\left(
\|\boldsymbol{\mathcal R}_f^{n+1}\|_{L^2(\omega)}^2
+\|\boldsymbol{\mathcal R}_s^{n+1}\|_{L^2(\omega)}^2
\right)
\leq
2\Delta t\int_0^T
\left[
\gamma^2
\|\partial_t\dot{\boldsymbol\eta}(t)\|_{L^2(\omega)}^2
+\|\partial_t\widehat{\boldsymbol t}_f(t)\|_{L^2(\omega)}^2
\right]\,dt.
\]
Multiplying by $(T+\Delta t)/(2\gamma)$ gives the final
residual term in \eqref{eq:error-final-estimate}.
Absorbing the Gronwall factors into C completes the proof.
\end{proof}

%%%%%%%%%%%%%%%%%%%%

\subsection{Convergence rate and initialization}
With exact initialization of the fluid and shell variables and both interface loads, $\mathcal E_{\mathrm{err,aug}}^0=0$.
For fixed $T$, $\gamma$, and an exact solution with the stated
regularity, the momentum and kinematic terms in
\eqref{eq:error-final-estimate} are $O(\Delta t^2)$, while the
interface term is $O(\Delta t)$. Hence
\begin{equation}\label{eq:error-squared-rate}
\max_{0\leq n\leq N}\mathcal E_{\mathrm{err,aug}}^n=O(\Delta t),
\qquad
\max_{0\leq n\leq N}\mathcal E_{\mathrm{err}}^n=O(\Delta t).
\end{equation}
Taking square roots gives the mechanical energy-norm bound
$O(\sqrt{\Delta t})$. Its component estimates are
\begin{equation}\label{eq:error-component-rates}
\begin{aligned}
\max_{0\leq n\leq N}\sqrt{\rho_f}
\|\boldsymbol e_u^n\|_{L^2(\widehat\Omega_f)}&=O(\sqrt{\Delta t}),\\
\max_{0\leq n\leq N}\sqrt{\varepsilon\rho_s}
\|\boldsymbol e_{\dot\eta}^n\|_{L^2(\omega)}&=O(\sqrt{\Delta t}),\\
\max_{0\leq n\leq N}\sqrt{\varepsilon^3\rho_s/12}
\|\boldsymbol e_{\dot\omega}^n\|_{L^2(\omega)}&=O(\sqrt{\Delta t}),\\
\max_{0\leq n\leq N}\Big[
B_m\big((\boldsymbol e_\eta^n,\boldsymbol e_\omega^n),
(\boldsymbol e_\eta^n,\boldsymbol e_\omega^n)\big)
+B_b\big((\boldsymbol e_\eta^n,\boldsymbol e_\omega^n),
(\boldsymbol e_\eta^n,\boldsymbol e_\omega^n)\big)
\Big]^{1/2}&=O(\sqrt{\Delta t}).
\end{aligned}
\end{equation}
The square root of the accumulated symmetric-gradient error has
the same order in its discrete time norm.

The last line controls the shell elastic seminorm. For the clamped
cylinder, \eqref{eq:fixed-unshifted-coercivity} directly gives the
full $H^1$ error bound. Under only shifted coercivity, we also need
to control the shell $L^2$ error. Summing the kinematic error equations
and applying Cauchy-Schwarz in time gives
\begin{equation}\label{eq:error-shell-l2-from-kinematics}
\begin{aligned}
\max_{0\leq n\leq N}\left(
\|\boldsymbol e_\eta^n\|_{L^2(\omega)}^2
+\|\boldsymbol e_\omega^n\|_{L^2(\omega)}^2\right)
&\leq3\left(\|\boldsymbol e_\eta^0\|_{L^2(\omega)}^2
+\|\boldsymbol e_\omega^0\|_{L^2(\omega)}^2\right)
+3T^2\max_{0\leq n\leq N}\left(
\|\boldsymbol e_{\dot\eta}^n\|_{L^2(\omega)}^2
+\|\boldsymbol e_{\dot\omega}^n\|_{L^2(\omega)}^2\right)\\
&\quad+3T\Delta t\sum_{n=0}^{N-1}\left(
\|\boldsymbol\kappa_\eta^{n+1}\|_{L^2(\omega)}^2
+\|\boldsymbol\kappa_\omega^{n+1}\|_{L^2(\omega)}^2\right).
\end{aligned}
\end{equation}
The final sum, with its $\Delta t$ weight, is $O(\Delta t^2)$
by \eqref{eq:cons-kinematic-integrals} and the regularity assumptions.
With exact initialization and fixed physical parameters, the squared
shell $L^2$ error is therefore $O(\Delta t)$. Combining this with
shifted coercivity and the elastic bound gives
\begin{equation}\label{eq:error-shell-h1-rate}
\max_{0\leq n\leq N}\left(
\|\boldsymbol e_\eta^n\|_{H^1(\omega)}^2
+\|\boldsymbol e_\omega^n\|_{H^1(\omega)}^2\right)^{1/2}
=O(\sqrt{\Delta t}).
\end{equation}

For approximate initialization under only shifted coercivity, this
argument additionally requires
$\|\boldsymbol e_\eta^0\|_{L^2(\omega)}^2+
\|\boldsymbol e_\omega^0\|_{L^2(\omega)}^2=O(\Delta t)$.
For the clamped cylinder, unshifted coercivity provides this control
through the initial augmented error bound, so no separate
initial $L^2$ condition is needed.

The mechanical energy-norm rate also holds for approximate initialization
provided $\mathcal E_{\mathrm{err,aug}}^0=O(\Delta t)$.
This condition includes the initial mechanical error, the fluid and
shell velocity trace errors, and both interface load errors.
For fixed $\gamma$, uniformly bounded initial errors in the interface
traces and loads contribute $O(\Delta t)$ to the augmented initial
error energy. To obtain the stated convergence rate, the initial
mechanical error energy must also be $O(\Delta t)$.
Larger initial errors remain in the final estimate.
The full shell $H^1$ rate then follows directly for the clamped
cylinder, or with the additional initial $L^2$ condition under
shifted coercivity.

Because the exact solution satisfies velocity continuity and force balance,
\begin{equation}\label{eq:error-defects-versus-numerical}
\begin{aligned}
\boldsymbol e_u^n\circ\varphi-\boldsymbol e_{\dot\eta}^n
=-(\boldsymbol u^n\circ\varphi-\dot{\boldsymbol\eta}^n),\qquad
\widehat{\boldsymbol e}_{t_f}^n+\boldsymbol e_{f_s}^n
=-(\widehat{\boldsymbol t}_f^n+\boldsymbol f_s^{\mathrm{ref},n}).
\end{aligned}
\end{equation}
Equation~\eqref{eq:error-final-estimate} therefore also controls the
numerical velocity mismatch and force imbalance. For fixed $\gamma$
and either exact initialization or
$\mathcal E_{\mathrm{err,aug}}^0=O(\Delta t)$, their weighted discrete
$L^2$ time norms are $O(\sqrt{\Delta t})$. These estimates use the old
time levels $n=0,\ldots,N-1$ at which the interface data are exchanged.
They follow from the error equations and consistency bounds in
addition to stability.

\begin{remark}[Scope of the error estimate]
The individual traction errors in the augmented energy carry
a factor of $\Delta t$. Consequently, the estimate does not
establish unweighted maximum-in-time convergence of the
individual tractions. The augmented energy alone also does
not establish unweighted maximum-in-time convergence of
the fluid velocity trace errors. The force-imbalance estimate
controls the sum of the two interface load errors.

The $O(\sqrt{\Delta t})$ convergence rate follows from the
present energy argument and does not rule out faster convergence
for particular solutions. Establishing a sharper rate would
require improved estimates of the interface residual work
or additional properties of the coupled problem.

The pressure terms cancel when testing the error equations,
so the present analysis does not provide an estimate for
$e_p^n$. Deriving a pressure error bound would require
applying the assumed inf-sup condition together with
additional estimates for the momentum and interface errors.
Extending the analysis to moving domains would require
accounting for geometry prediction, domain transport,
the associated Jacobians, mesh velocity, and convective
effects. A fully discrete error analysis would also require
incorporating spatial projection errors and compatible
interface-transfer operators.
\end{remark}

\section{Numerical examples}\label{sec:numerical}

\subsection{Pressure-wave propagation in a straight tube}\label{straight_tube}
We evaluate the proposed scheme using a moving-domain benchmark problem involving a pressure pulse in a straight compliant tube. We compare its results with those of the $\beta$-scheme \cite{WangTambacaCanic2026} and a monolithic solver. We first examine the effect of varying the Robin parameter $\gamma$ at a fixed time step, then investigate temporal refinement with $\gamma$ held fixed. The comparisons focus on fluid pressure and radial wall displacement. 

The reference tube has length $L=0.06\,\mathrm{m}$ and radius
$R_0=0.005\,\mathrm{m}$. Its wall motion is governed by the Naghdi shell model described in
Section~\ref{sec:model}. All three solvers use the same membrane,
shear, bending, and inertia operators and the parameters as listed in
Table~\ref{tab:numerical-parameters}.  The pressure-pulse setting is based
on straight-tube $\beta$-scheme benchmarks~\cite{BukacEtAl2013}.
\begin{table}[htbp]
\centering
\caption{Physical parameters of the pressure-pulse in a straight tube study.}
\label{tab:numerical-parameters}
\begin{tabular}{llllll}
\hline
Quantity & Symbol & Value & Quantity & Symbol & Value\\
\hline
Fluid density & $\rho_f$ & $1000\,\mathrm{kg\,m^{-3}}$ &
Dynamic viscosity & $\mu_f$ & $3.5\times10^{-3}\,\mathrm{Pa\,s}$\\
Shell density & $\rho_s$ & $1100\,\mathrm{kg\,m^{-3}}$&
Young's modulus & $E$ & $7.5\times10^4\,\mathrm{Pa}$\\
Poisson ratio & $\nu$ & $0.5$&
Shell thickness & $h_s$ & $10^{-3}\,\mathrm{m}$\\
Reference radius & $R_0$ & $5\times10^{-3}\,\mathrm{m}$&
Tube length & $L$ & $6\times10^{-2}\,\mathrm{m}$\\
\hline
\end{tabular}
\end{table}
Initially, the fluid and shell are at rest, with zero displacement
and interface loads. The shell is clamped at both ends and periodic
in its circumferential parameter. The inlet traction is
$-p_{\mathrm{in}}(t)\boldsymbol{n}_x$, with
\begin{equation}\label{eq:numerical-pulse}
p_{\mathrm{in}}(t)=
\begin{cases}
\displaystyle\frac{p_{\max}}{2}
 \left(1-\cos\frac{2\pi t}{t_p}\right),&0\leq t\leq t_p,\\
0,&t>t_p,
\end{cases}
\qquad
p_{\max}=2000\,\mathrm{Pa},\quad t_p=5\,\mathrm{ms}.
\end{equation}
The outlet traction is zero. Each simulation runs from
$0$ to $12\,\mathrm{ms}$. We exam the pressure and displacement profiles at
$t_k\in\{2,4,6,8,10,12\}\,\mathrm{ms}$ 
for our comparisons.

\begin{figure}[htbp]
\centering
\includegraphics[width=3in]{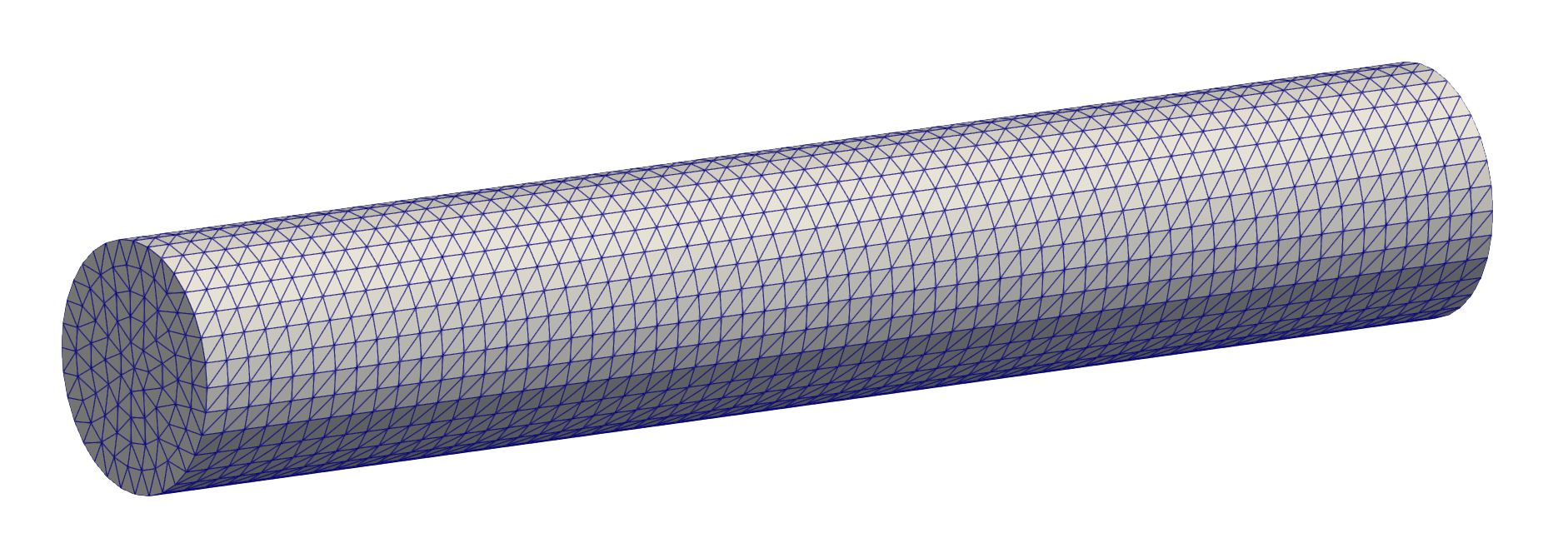}
\includegraphics[width=3in]{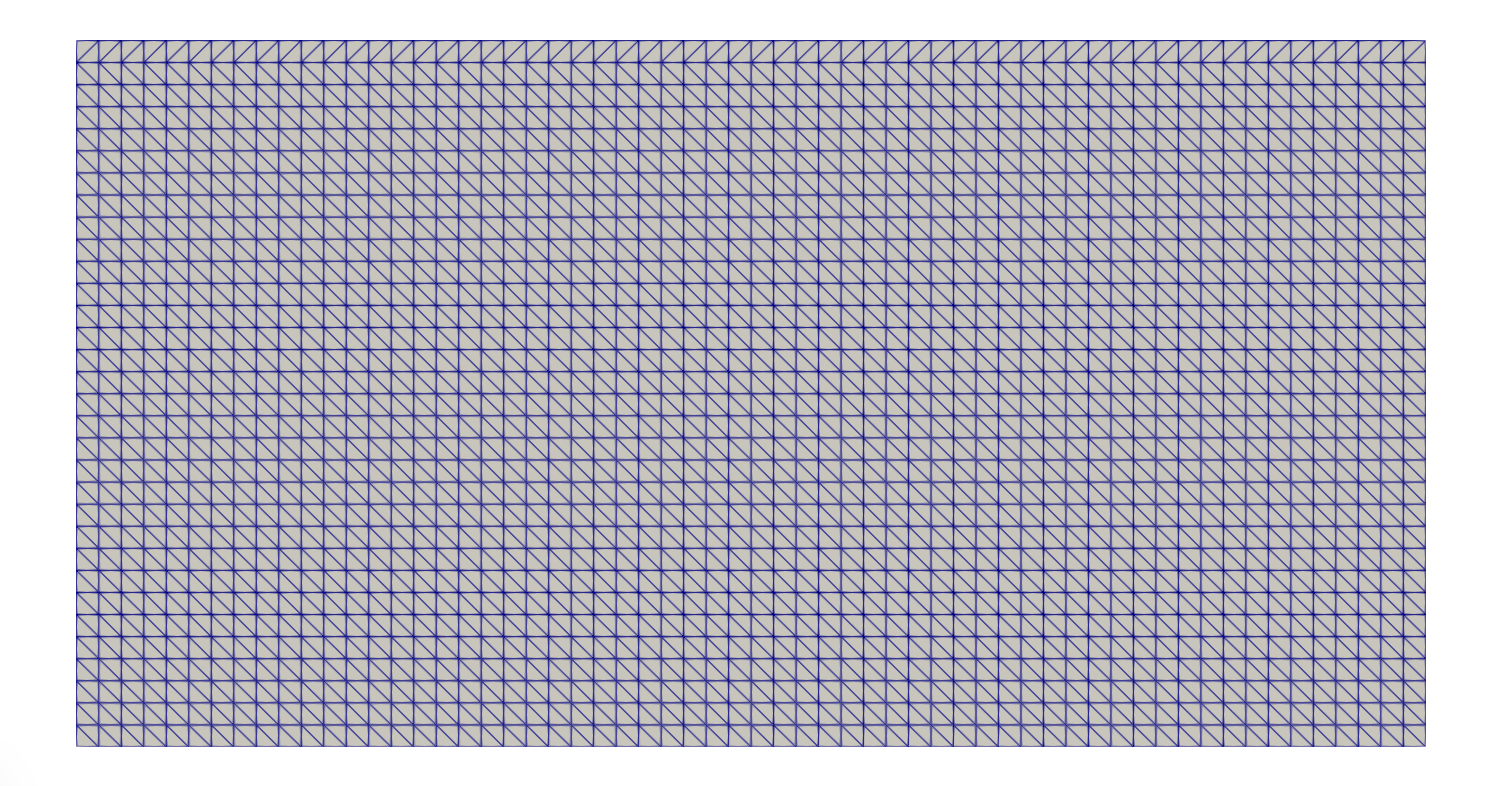}
\caption{Reference meshes: (a) fluid exterior, (b) shell mesh in the parametric domain.}
\label{fig:numerical-mesh}
\end{figure}

Figure~\ref{fig:numerical-mesh} shows the meshes used in our test. The fluid mesh
contains $5978$ vertices and $29160$ tetrahedra. The shell parameter
mesh has $3840$ triangles and $2013$ vertices. Fluid velocity and
pressure use continuous $P_1$ elements, with Brezzi-Pitk\"aranta
stabilization term $0.01\sum_K h_K^2(\nabla p_h,\nabla q_h)_K$.
All six shell components use continuous $P_1$ elements.
The two partitioned solvers use the same lagged skew-convection
form and pressure stabilization, but advance geometry and loads
differently. The $\beta$ scheme uses the new shell-stage displacement
for mesh motion and projects Cauchy traction. Our explicit scheme
predicts geometry from the old state and recovers the load from its
Robin condition. The monolithic solver instead iterates the current
geometry and convecting velocity. 
For our comparison, $\beta$ scheme uses $\beta=1$ and $\Delta t=10^{-4}\,\mathrm{s}$. The monolithic solver uses the same
time step, and within each step,
Picard iterations update the ALE geometry and skew convection until
both velocity increments meet a relative tolerance of $10^{-7}$. A step is rejected if it has not
converged within $30$ iterations.

For result comparison, we divide the reference axial interval into $61$ equal bins. Fluid cells are assigned by their reference centroids, and we use $B_i(t)$
to denote the current cells in bin $i$, and $S_i$ to denote the corresponding shell parameter cells. We analyze the following data profiles:
\begin{equation}\label{eq:numerical-profiles}
\overline p_i(t)=\frac{1}{|B_i(t)|}\int_{B_i(t)}p_h(t)\,dx,
\qquad
\overline\eta_{r,i}(t)=\frac{1}{|S_i|}\int_{S_i}
 \boldsymbol\eta_h(t)\cdot
(0,\cos(y/R_0),\sin(y/R_0))\,dy\,dx.
\end{equation}
The shell integral uses axial and arc-length coordinates $x$ and
$y$. Both quantities average over a slab of cells, rather than sample
an exact cross-section. The radial displacement is also obtained
from the stored mean diameter through
$\overline\eta_{r,i}=(\overline D_i-2R_0)/2$.

We take the $\beta$ solution as the numerical comparison standard.
For either $q=\overline p$ or $q=\overline\eta_r$, define
\begin{equation}\label{eq:numerical-discrepancy}
\mathcal D_q
=100\left(
\frac{\displaystyle\sum_{k=1}^{6}\sum_{i=1}^{61}
 |q_i(t_k)-q_i^\beta(t_k)|^2}
{\displaystyle\sum_{k=1}^{6}\sum_{i=1}^{61}|q_i^\beta(t_k)|^2}
\right)^{1/2}\!\%.
\end{equation}
Here, $k$ indexes the time snapshots and $i$ indexes the equally
spaced bins. This measure compares six profiles at matching physical
times, without phase shifts or amplitude rescaling, and is distinct
from a full-field error or the analytical energy norm. We also compute
discrepancies relative to the monolithic solution by replacing
$q^\beta$ with the corresponding monolithic values. The $\beta$
and monolithic solutions differ by $4.09\%$ in pressure and $4.97\%$
in radial displacement, providing a numerical cross-check, although
neither is established as an exact or mesh-converged reference.

\begin{figure}[!htbp]
\centering
\includegraphics[width=\linewidth]{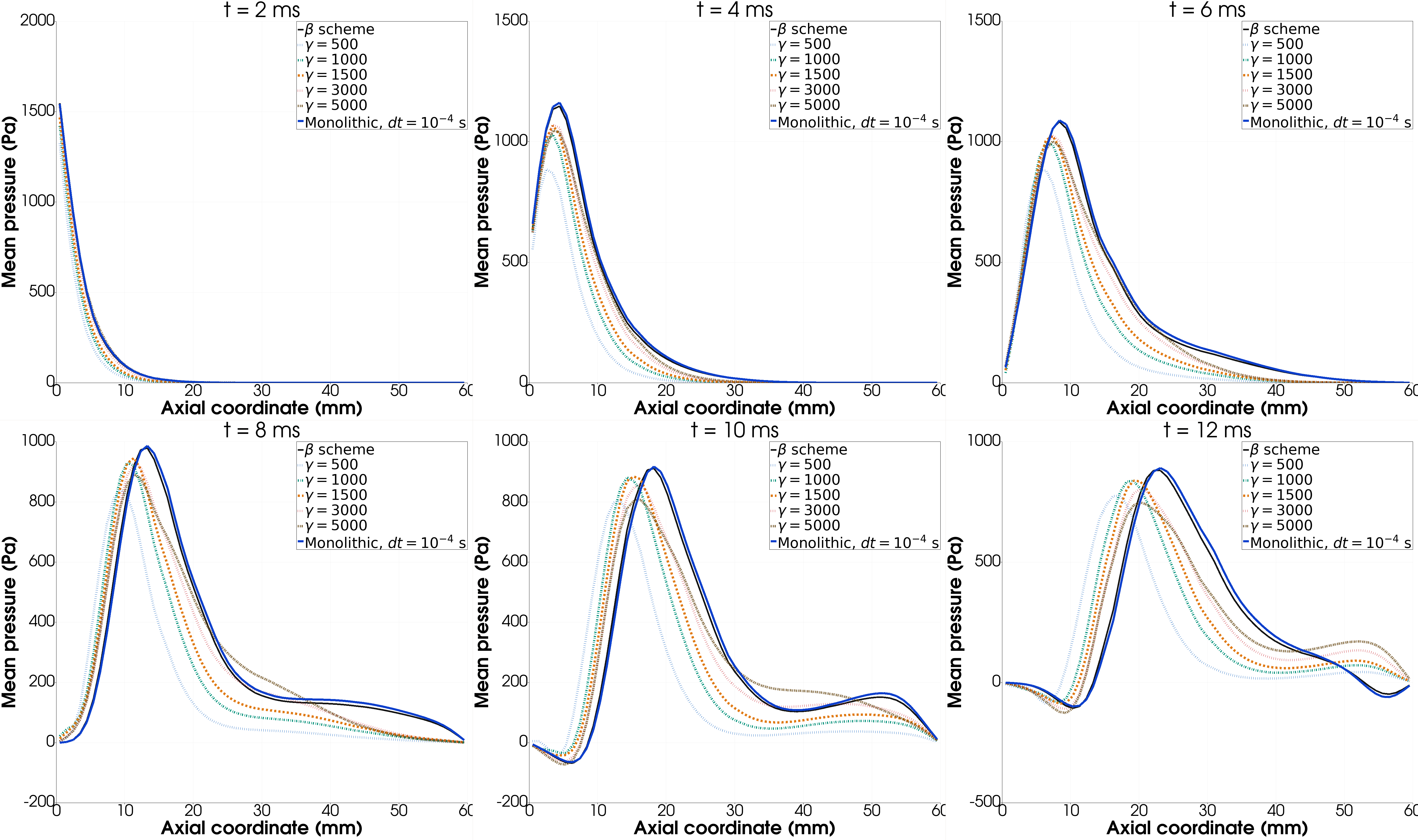}
\caption{Robin parameter sensitivity study: slab-averaged pressure profiles at
$t=2,4,6,8,10,12$ ms. All curves use $\Delta t=10^{-4}\,\mathrm{s}$, and
the $\beta$ scheme and monolithic solution are included for comparison.}
\label{fig:numerical-robin-pressure}
\end{figure}

\begin{figure}[!htbp]
\centering
\includegraphics[width=\linewidth]{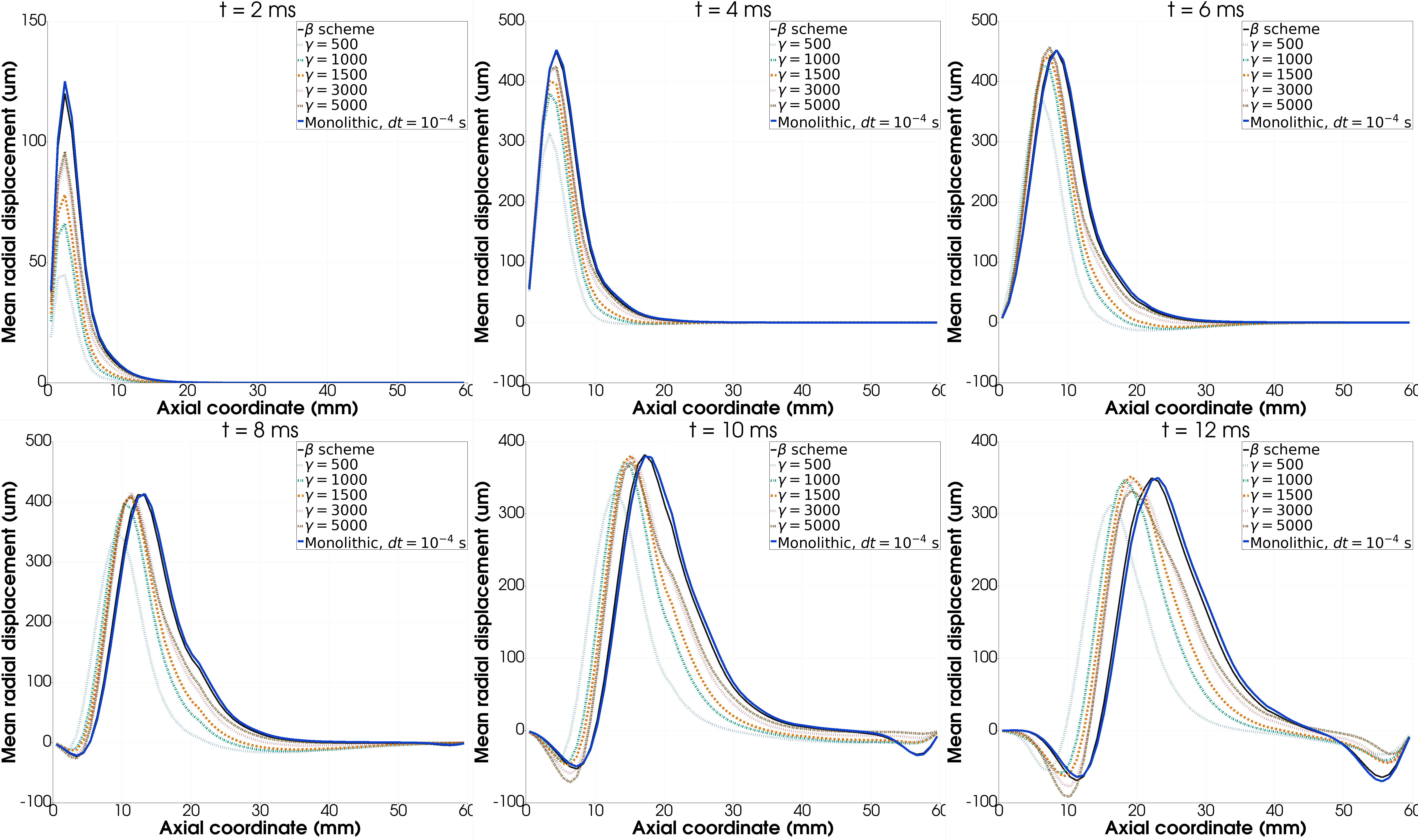}
\caption{Robin parameter sensitivity study: mean radial wall displacement at six different time moments. All curves use $\Delta t=10^{-4}\,\mathrm{s}$, and
the $\beta$ scheme and monolithic solution are included for comparison. Displacement is reported in micrometres.}
\label{fig:numerical-robin-radial}
\end{figure}

\paragraph{Sensitivity study of the Robin parameter}
\label{subsec:numerical-robin}
We hold $\Delta t=10^{-4}\,\mathrm{s}$ and take
$\gamma\in\{500,1000,1500,3000,5000\}$ $\mathrm{kg\,m^{-2}\,s^{-1}}$.
Figures~\ref{fig:numerical-robin-pressure} and
\ref{fig:numerical-robin-radial} show the pressure and wall motion, and
Table~\ref{tab:numerical-robin} summarizes their discrepancies.
For smaller $\gamma$, the results obtained from the explicit scheme differ substantially
from both solutions obtained from the $\beta$ scheme and monolithic solver. However, increasing $\gamma$ reduces the discrepancies in both quantities, and the improvement becomes marginal as $\gamma$ increases from $3000$ to $5000$. Thus, $\gamma = 3000$ is a reasonable choice for the current mesh and time step considered.
\begin{table}[htbp]
\centering
\caption{Robin parameter sensitivity study with $\Delta t=10^{-4},\mathrm{s}$. Discrepancies relative to the $\beta$ scheme solutions are defined by \eqref{eq:numerical-discrepancy} and evaluated using six snapshots at $2$, $4$, $6$, $8$, $10$, and $12,\mathrm{ms}$.}
\vspace{0.1in}
\label{tab:numerical-robin}
\begin{tabular}{rrr}
\hline
$\gamma\ (\mathrm{kg\,m^{-2}\,s^{-1}})$ & $\mathcal D_{\overline p}$ &
$\mathcal D_{\overline\eta_r}$\\
\hline
500  & $56.35\%$ & $67.19\%$\\
1000 & $37.05\%$ & $45.13\%$\\
1500 & $27.65\%$ & $34.37\%$\\
3000 & $18.28\%$ & $23.98\%$\\
5000 & $18.06\%$ & $23.89\%$\\
\hline
\end{tabular}
\end{table}

\paragraph{Time refinement at fixed Robin parameter}
\label{subsec:numerical-time}

\begin{figure}[!htbp]
\centering
\includegraphics[width=\linewidth]{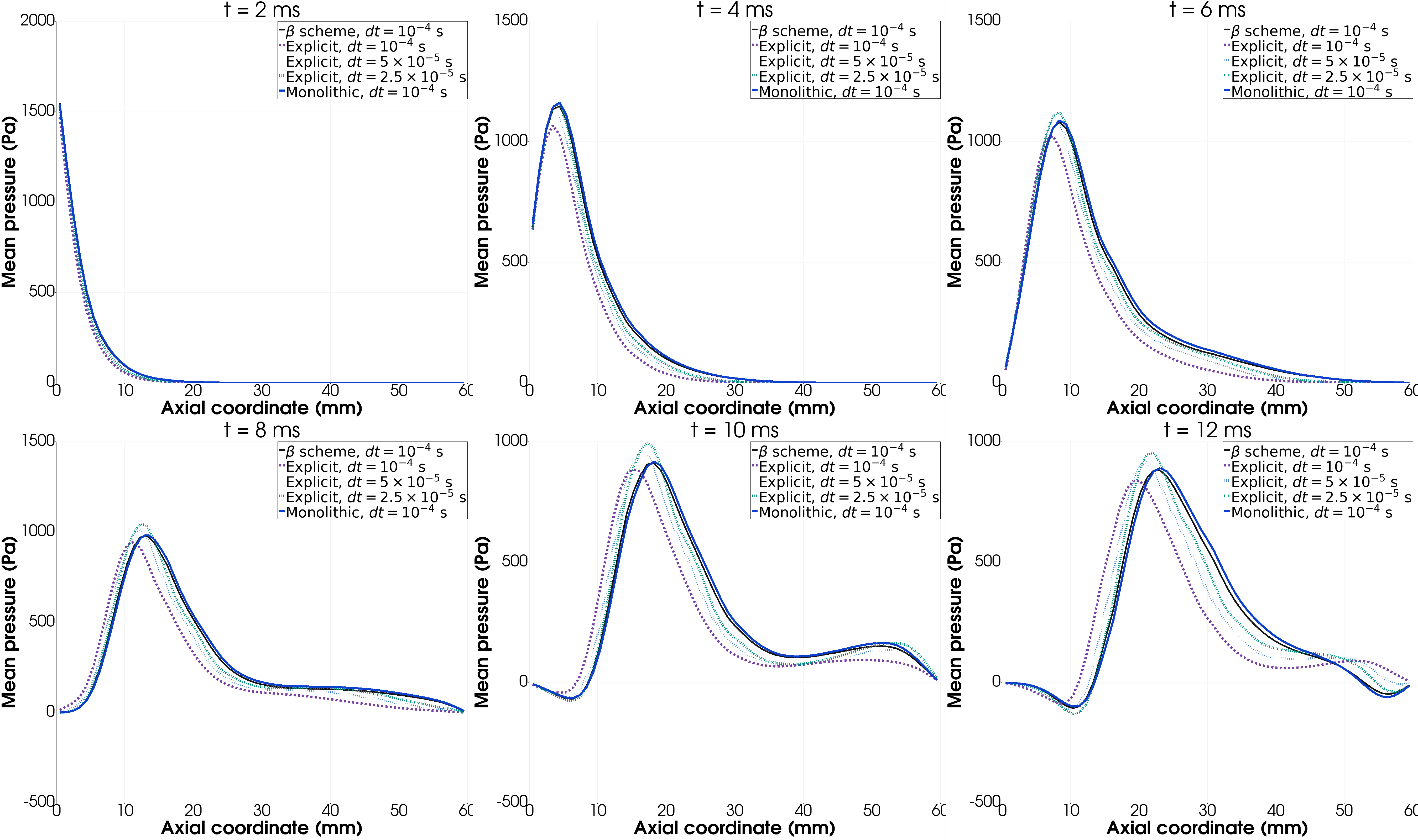}
\caption{Time refinement at $\gamma=1500$: pressure profiles at the six
snapshot times. Only the explicit splitting scheme refine time step from $t=$ ($10^{-4}$, $5\times 10^{-5}$, $2.5\times 10^{-5}$ ). Both $\beta$ scheme and monolithic comparison
solvers use $\Delta t=10^{-4}\,\mathrm{s}$.}
\label{fig:numerical-time1500-pressure}
\end{figure}

\begin{figure}[!htbp]
\centering
\includegraphics[width=\linewidth]{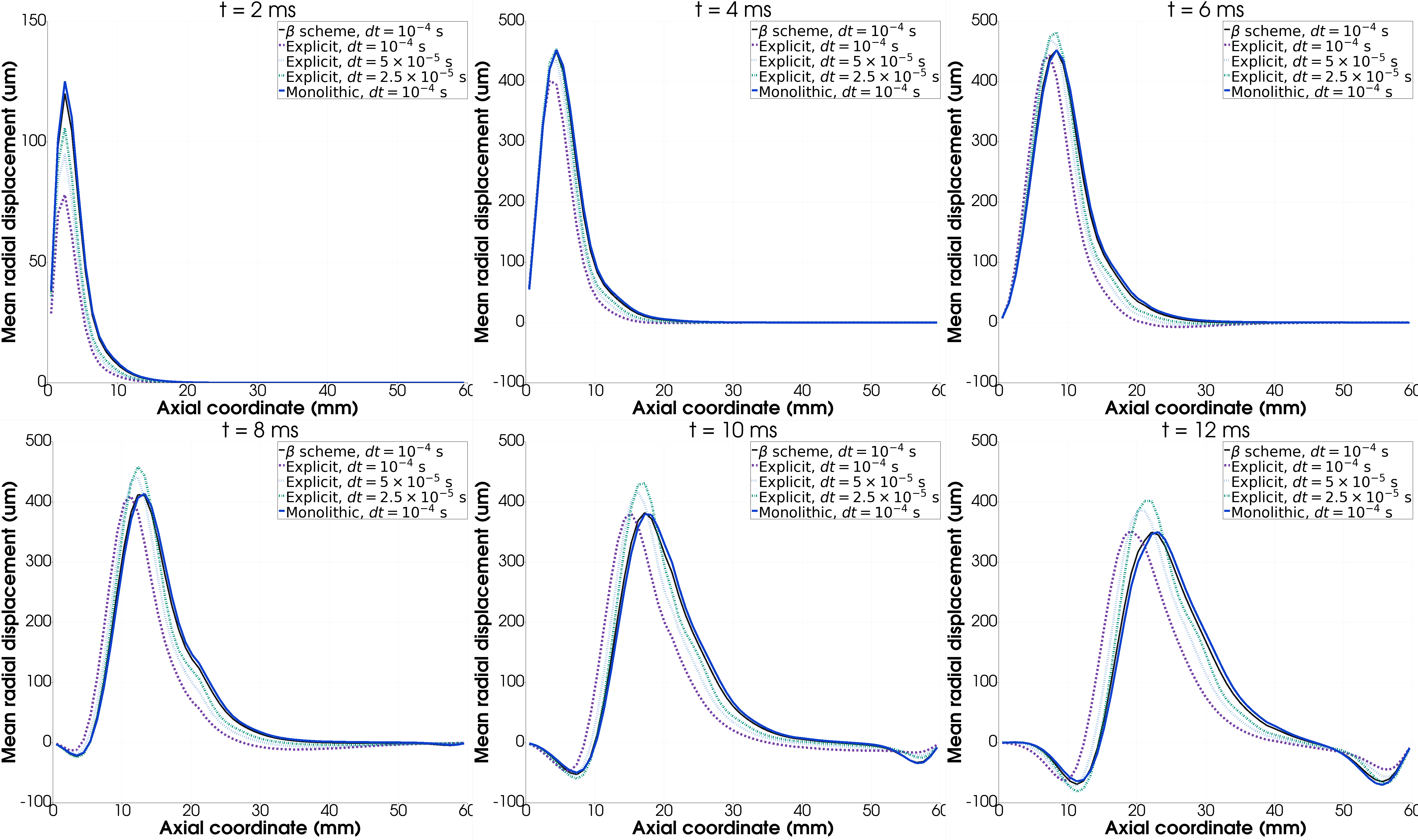}
\caption{Time refinement study with $\gamma=1500$: mean radial displacement at six selected times. Only the explicit splitting scheme is refined in time, using $\Delta t=10^{-4}$, $5\times10^{-5}$, and $2.5\times10^{-5},\mathrm{s}$. Both the $\beta$ scheme and the monolithic solver use $\Delta t=10^{-4},\mathrm{s}$.}
\label{fig:numerical-time1500-radial}
\end{figure}

\begin{figure}[!htbp]
\centering
\includegraphics[width=\linewidth]{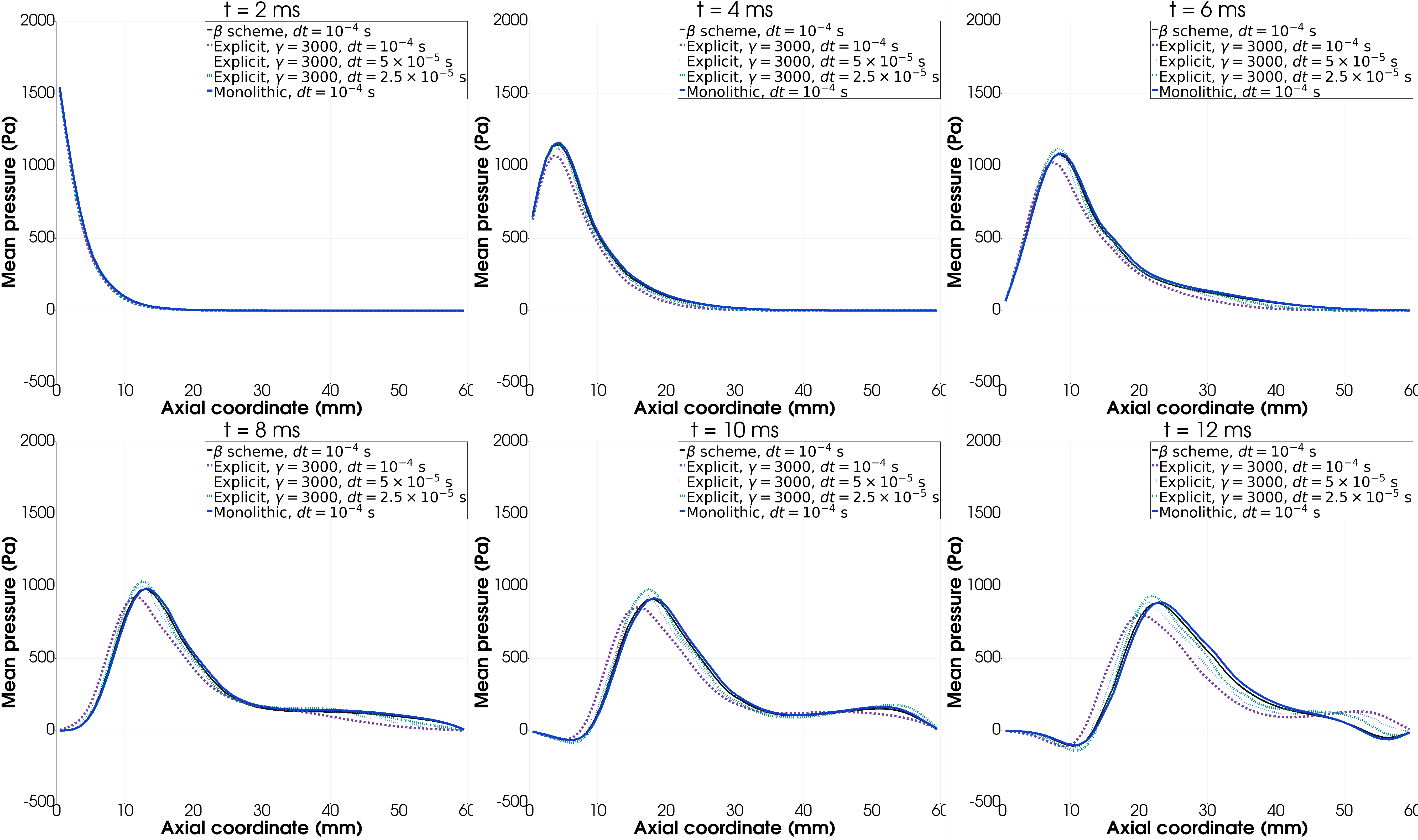}
\caption{Time refinement study with $\gamma=3000$: pressure profiles at six selected times. The time steps used for the explicit scheme and comparison solvers are the same as those in Figure~\ref{fig:numerical-time1500-pressure}.}
\label{fig:numerical-time3000-pressure}
\end{figure}

\begin{figure}[!htbp]
\centering
\includegraphics[width=\linewidth]{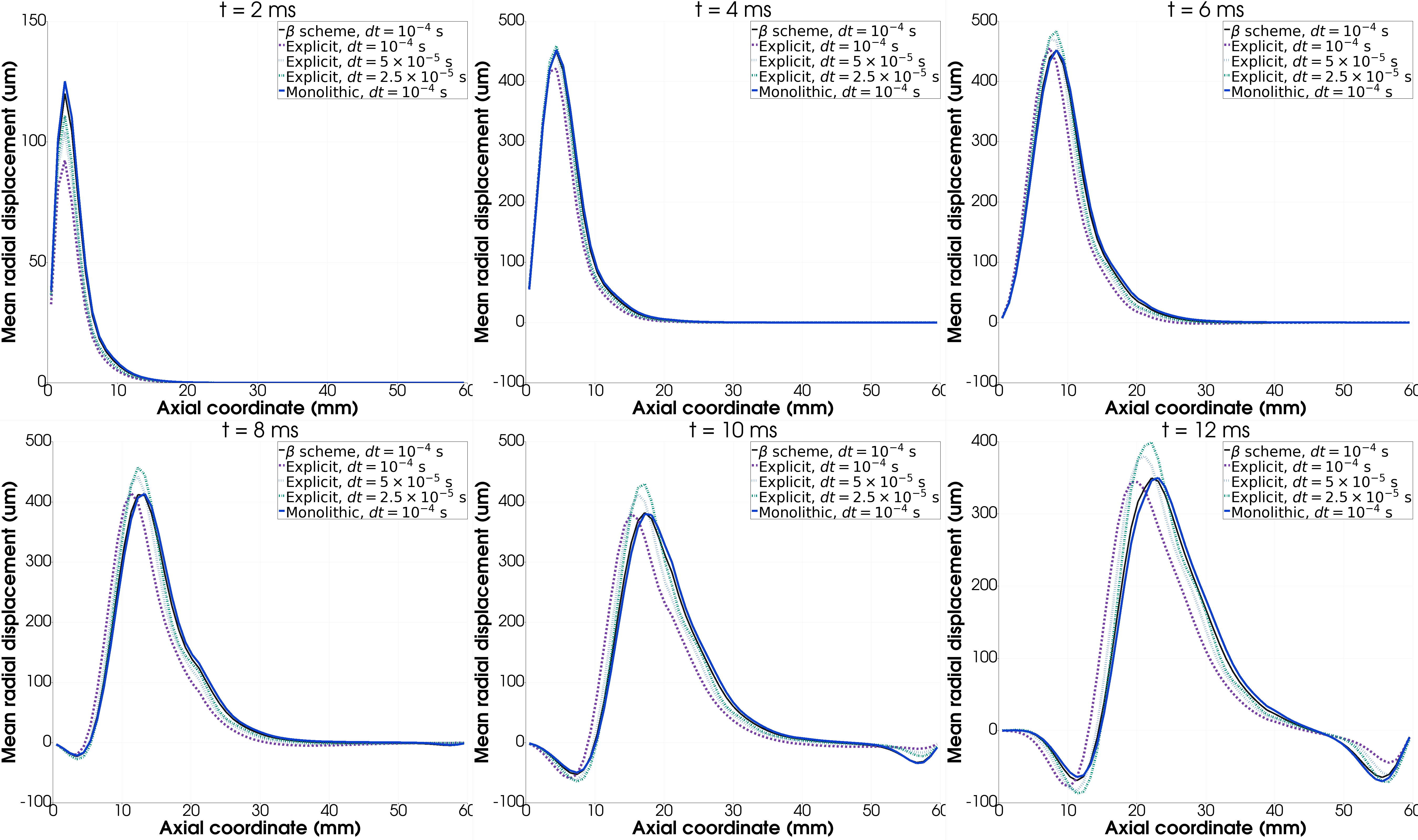}
\caption{Time refinement study with $\gamma=3000$: mean radial displacement,
comparing with Figure~\ref{fig:numerical-time1500-radial}.}
\label{fig:numerical-time3000-radial}
\end{figure}

We study $\gamma=1500$ and $3000$ separately, using explicit
time steps $10^{-4}$, $5\times10^{-5}$, and
$2.5\times10^{-5}\,\mathrm{s}$. All other settings are unchanged.
The $\beta$ and monolithic solutions stay at
$\Delta t=10^{-4}\,\mathrm{s}$.
Figures~\ref{fig:numerical-time1500-pressure}--
\ref{fig:numerical-time3000-radial} show the six profile times,
and Table~\ref{tab:numerical-time} reports the corresponding
discrepancies against the solution from $\beta$-scheme. The table uses only the displayed
$2$ to $12$ ms snapshots.

\begin{table}[htbp]
\centering
\caption{Explicit time refinement at $\gamma=1500$ and $3000$. The $\beta$-scheme
comparison uses $\Delta t=10^{-4}\,\mathrm{s}$. Pressure and radial-displacement discrepancies use the six
snapshots at $2$--$12$ ms and $61$ axial bins.}
\label{tab:numerical-time}
\begin{tabular}{rrrr}
\hline
$\gamma$ & $\Delta t\ (\mathrm{s})$ & $\mathcal D_{\overline p}$ & $\mathcal D_{\overline\eta_r}$\\
\hline
1500 & $10^{-4}$ & $27.65\%$ & $34.37\%$\\
1500 & $5\times10^{-5}$ & $15.44\%$ & $19.42\%$\\
1500 & $2.5\times10^{-5}$ & $8.69\%$ & $12.08\%$\\
\hline
3000 & $10^{-4}$ & $18.28\%$ & $23.98\%$\\
3000 & $5\times10^{-5}$ & $9.97\%$ & $13.28\%$\\
3000 & $2.5\times10^{-5}$ & $6.28\%$ & $9.64\%$\\
\hline
\end{tabular}
\end{table}

Both pressure and displacement discrepancies decrease at each
refinement in Table~\ref{tab:numerical-time}. The same trend holds
against the monolithic solution. For $\gamma=1500$, the pressure discrepancy changes from $27.65\%$ to $8.69\%$, and the displacement discrepancy from $34.37\%$ to $12.08\%$. For $\gamma=3000$, the pressure discrepancy changes from $18.28\%$ to $6.28\%$, and the displacement discrepancy from $23.98\%$ to $9.64\%$.
Both measures decrease at each refinement for both Robin parameters.
At each tested time step, $\gamma=3000$ gives smaller discrepancies than $\gamma=1500$. Together with the fixed-step study, this shows
that the useful choice of $\gamma$ depends on the time step and the quantity being compared. Noticing that the computations include interface mismatch, load imbalance, and geometry staging, so they do not separate a single
mechanism responsible for the remaining differences.

\subsection{Pressure-wave propagation in a stented artery}
\label{sec:numerical-stented}

We next apply the explicit splitting scheme to an artery supported with a stent. 
The reference artery has a length of $33\,\mathrm{mm}$ and a radius of
$1.5\,\mathrm{mm}$. Let $s\in[0,33]\,\mathrm{mm}$ denote distance
from the inlet and $\theta$ the circumferential angle. The connected
strut network lies in the central region, approximately
$6\leq s\leq27\,\mathrm{mm}$. It has eight axial cells and eight
circumferential repeats, with strut width $0.13\,\mathrm{mm}$.

We reported all the parameters used in this example in Table~\ref{tab:stented-materials}. 
\begin{table}[htbp]
\centering
\caption{Layer parameters used in the stented-artery study.}
\label{tab:stented-materials}
\begin{tabular}{lcc}
\hline
Parameter & Artery wall & Outer steel layer\\
\hline
Young's modulus ($\mathrm{Pa}$) & $4\times10^5$ & $2.43\times10^{10}$\\
Thickness ($\mathrm{mm}$) & $0.4$ & $0.15$\\
Density ($\mathrm{kg\,m^{-3}}$) & $1100$ & $1100$\\
Poisson ratio & $0.5$ & $0.5$\\
\hline
\end{tabular}
\end{table}
The shell stiffness is assembled by integrating through both layers,
rather than replacing the wall by steel. The total shell thickness is
$0.4\,\mathrm{mm}$ at the two ends and $0.55\,\mathrm{mm}$ at the stented region. 
To specify this extension of
the homogeneous model, let $z$ measure distance outward from the artery
midsurface, let $h_a$ and $h_s$ be the artery and stent thicknesses,
and let $\chi_s$ be the indicator of the strut bands. The elastic
moments used in the implementation are
\begin{equation}
\begin{split}
 A_j(s,\theta)&=E_a\int_{-h_a/2}^{h_a/2}z^j\,dz
+\chi_s(s,\theta)E_s\int_{h_a/2}^{h_a/2+h_s}z^j\,dz,\\
 m_j(s,\theta)&=\rho_a\int_{-h_a/2}^{h_a/2}z^j\,dz
+\chi_s(s,\theta)\rho_s\int_{h_a/2}^{h_a/2+h_s}z^j\,dz,
 \qquad j=0,1,2.
\end{split}
\label{eq:stented-layer-moments}
\end{equation}
With the common Poisson ratio, $A_0$, $A_1$, and $A_2$ supply the
membrane, membrane-bending coupling, and bending coefficients in the Naghdi equation. The corresponding inertial terms use $m_0$, $m_1$,
and $m_2$. The first moments are retained because the added layer is outside the artery midsurface. This is a layered shell approximation.

\begin{figure}[!htbp]
\centering
\includegraphics[width=0.9\linewidth]{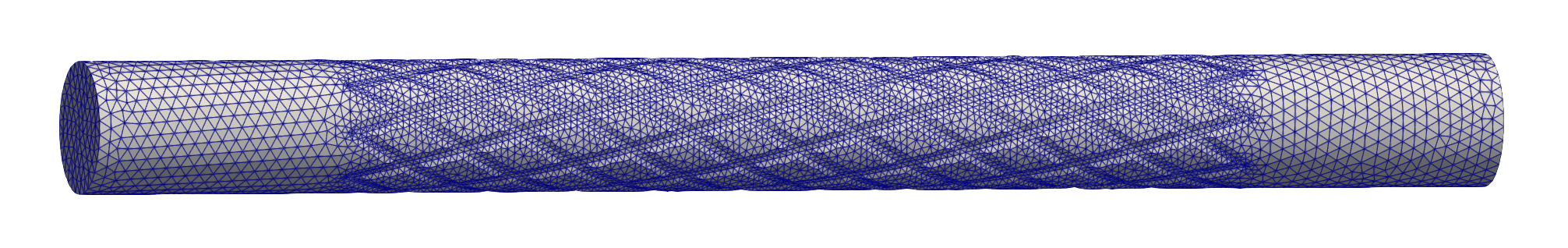}
\includegraphics[width=0.9\linewidth]{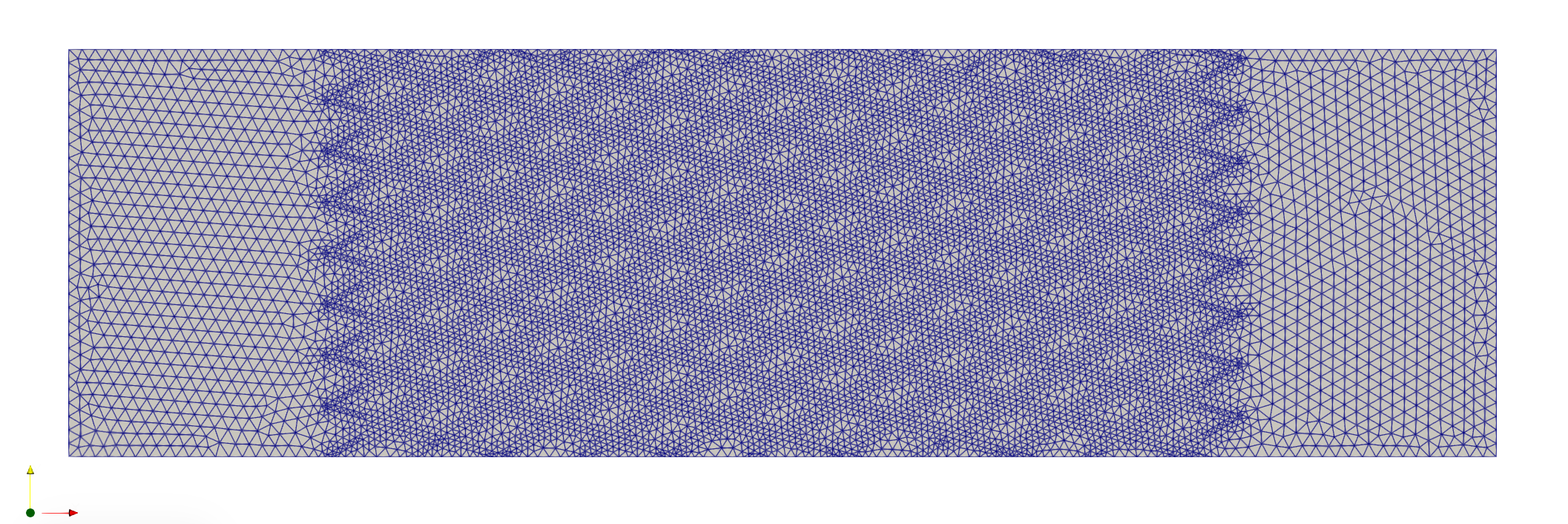}
\caption{Stented-artery mesh: (Top) fluid mesh; (Bottom) shell mesh in the parametric domain.}
\label{fig:stented-mesh}
\end{figure}

Figure~\ref{fig:stented-mesh} shows the fitted material mesh. The shell
has $21,032$ triangles and $10,705$ parameter-mesh vertices. The fluid mesh has $156,320$ tetrahedra and $30,620$ vertices. 
Fluid velocity and pressure use continuous $P_1$ elements with pressure
stabilization coefficient $0.01$. The six shell components also use
continuous $P_1$ elements, with periodicity around the circumference and clamped ends. Harmonic extension supplies the ALE displacement.
The fluid properties, inlet pressure pulse, and outlet boundary conditions are the same as those in our first numerical example in Section~\ref{straight_tube}. The initial state is at rest. The explicit calculation uses
$\gamma=2000$, $\Delta t=10^{-4}\,\mathrm{s}$, and the implemented
lagged skew-convection form. 

\begin{figure}[!htbp]
\centering
\includegraphics[width=\linewidth]{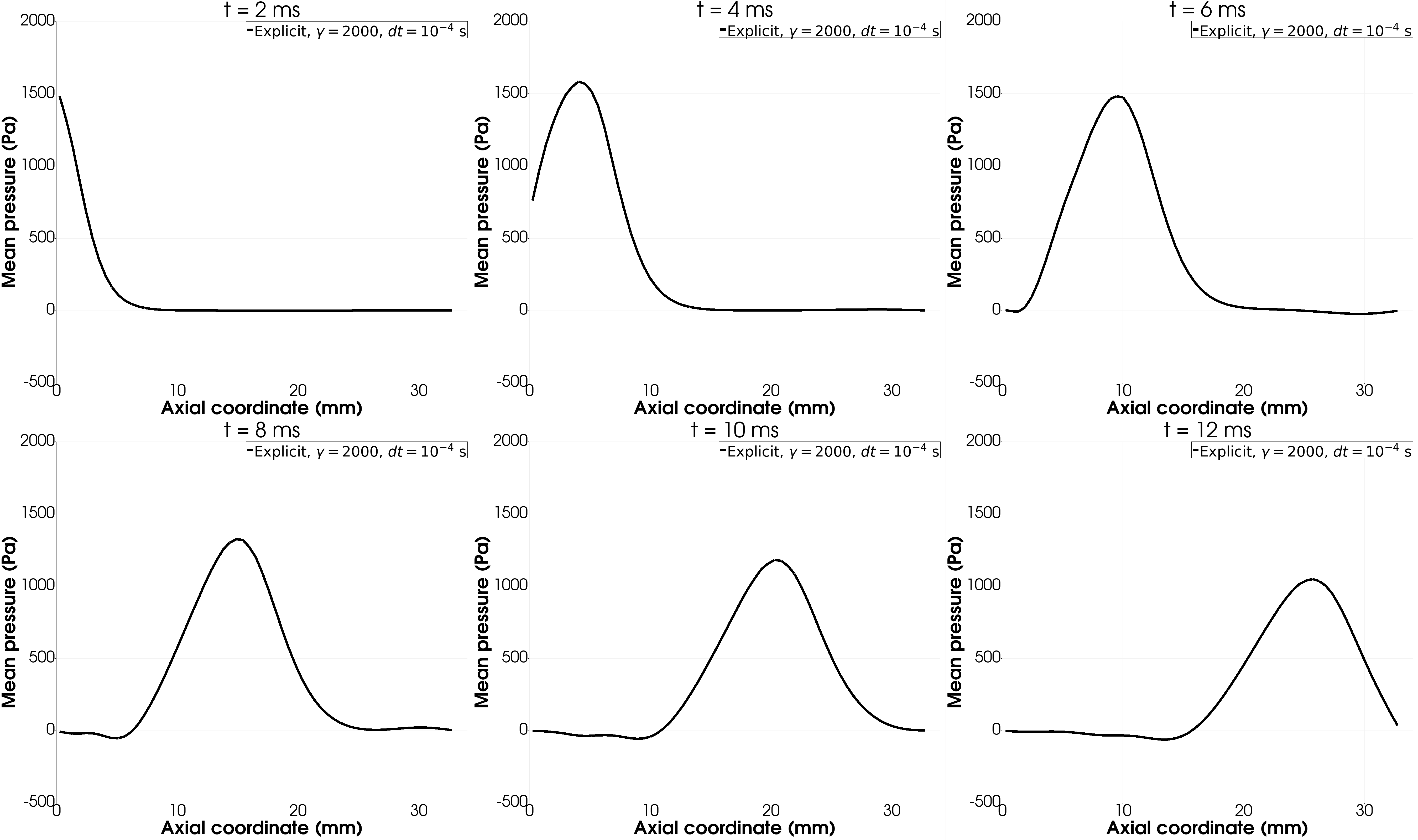}
\caption{Stented artery pressure profiles along the center axis at six different time moments, obtained from using the explicit splitting scheme with $\gamma=2000$ and $\Delta t=10^{-4}\,\mathrm{s}$.}
\label{fig:stented-pressure}
\end{figure}

\begin{figure}[htbp]
    \centering
        \includegraphics[width=0.9\textwidth]{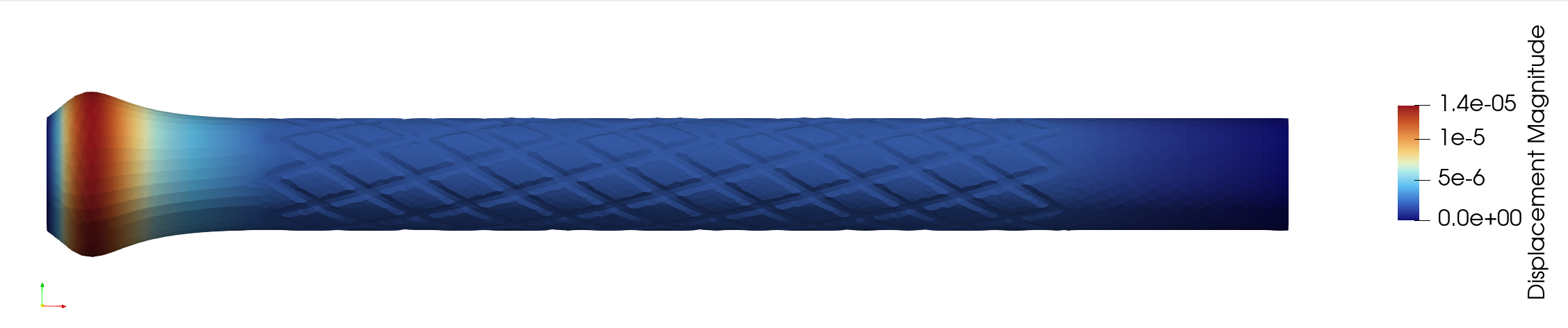}
        \includegraphics[width=0.9\textwidth]{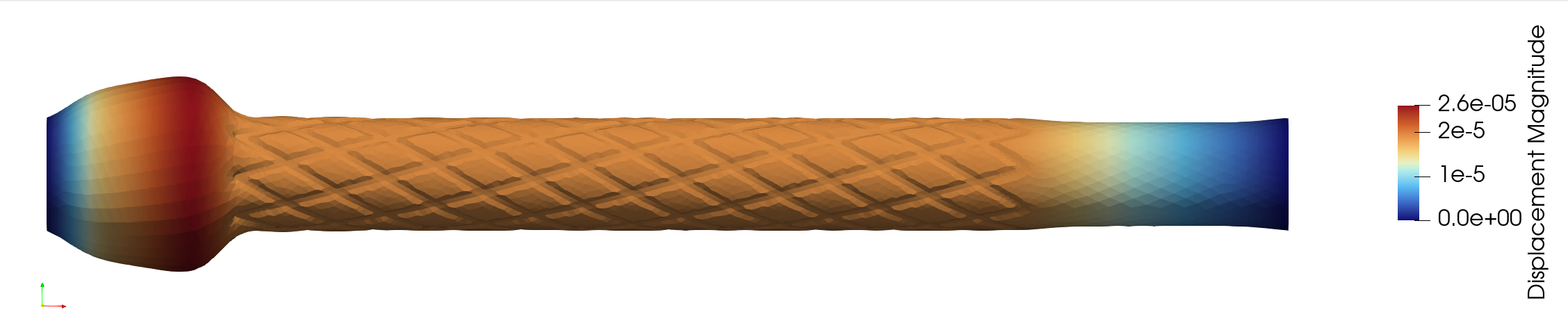}
        \includegraphics[width=0.9\textwidth]{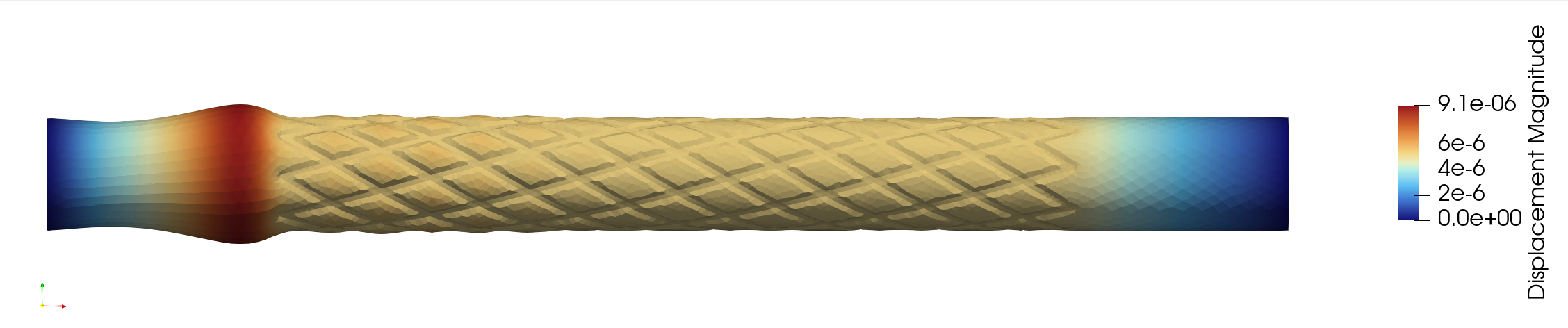}
 \includegraphics[width=0.9\textwidth]{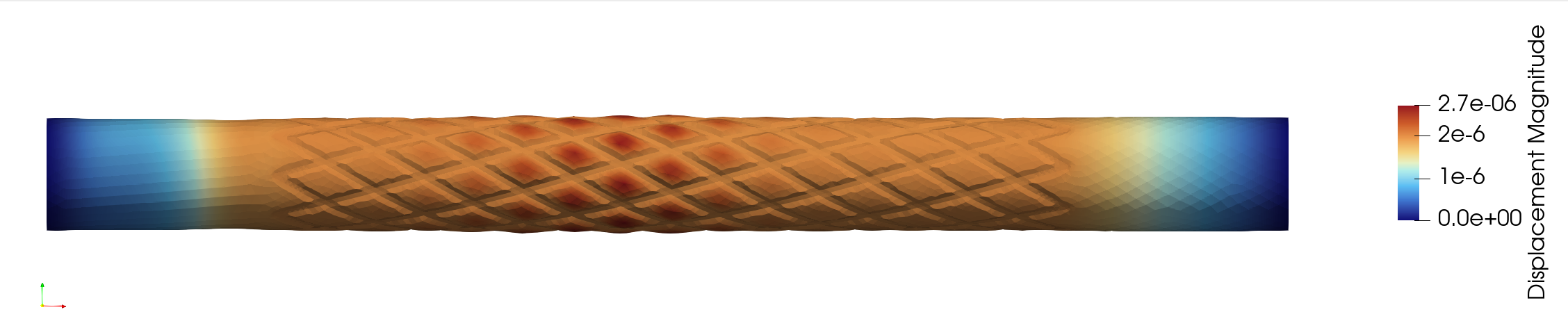}
        \includegraphics[width=0.9\textwidth]{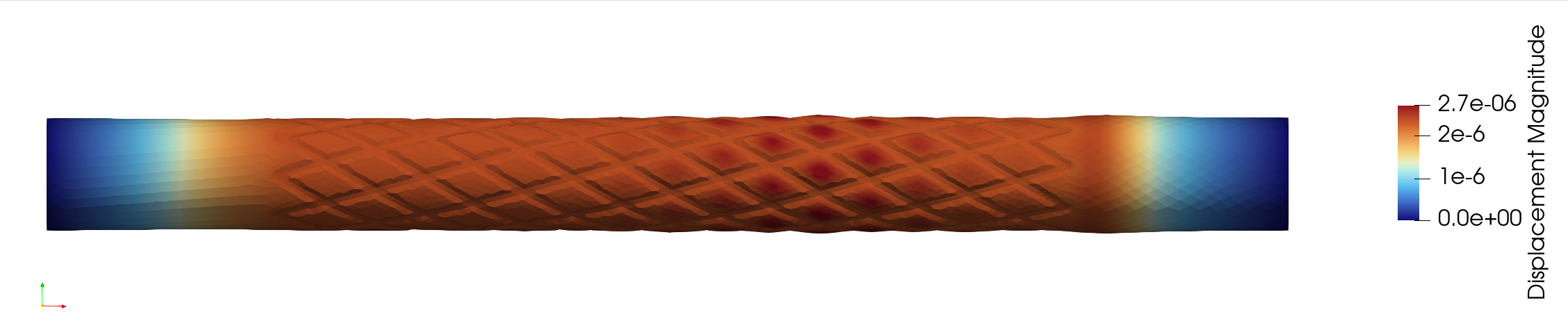}
        \includegraphics[width=0.9\textwidth]{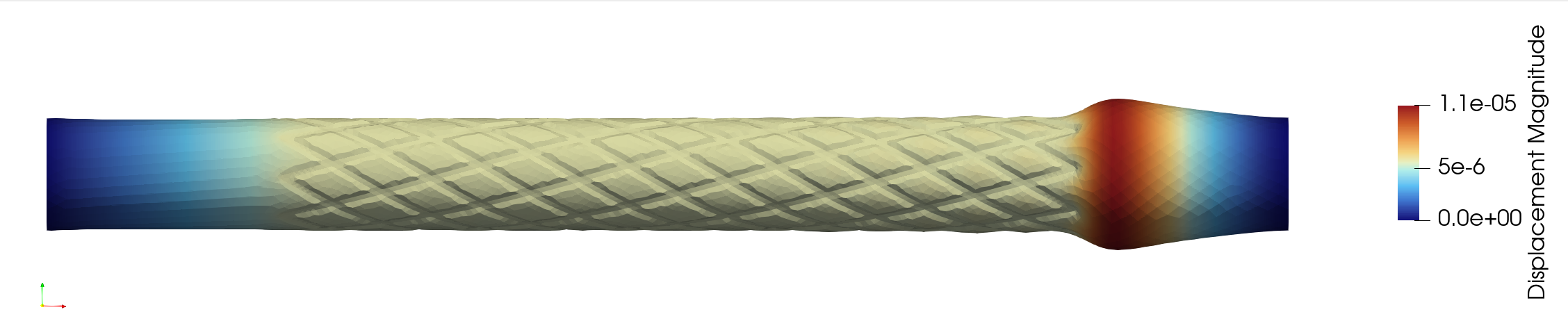}
\caption{Displacement of the stented artery at six time snapshots,
$t=2,4,6,8,10,$ and $12\,\mathrm{ms}$, obtained from using the explicit
Robin-Robin splitting scheme. The displacement is magnified by 10 times for better visualization of the arterial wall deformation.}
    \label{fig:stented_displacement}
\end{figure}

\begin{figure}[!htbp]
\centering
\includegraphics[width=\linewidth]{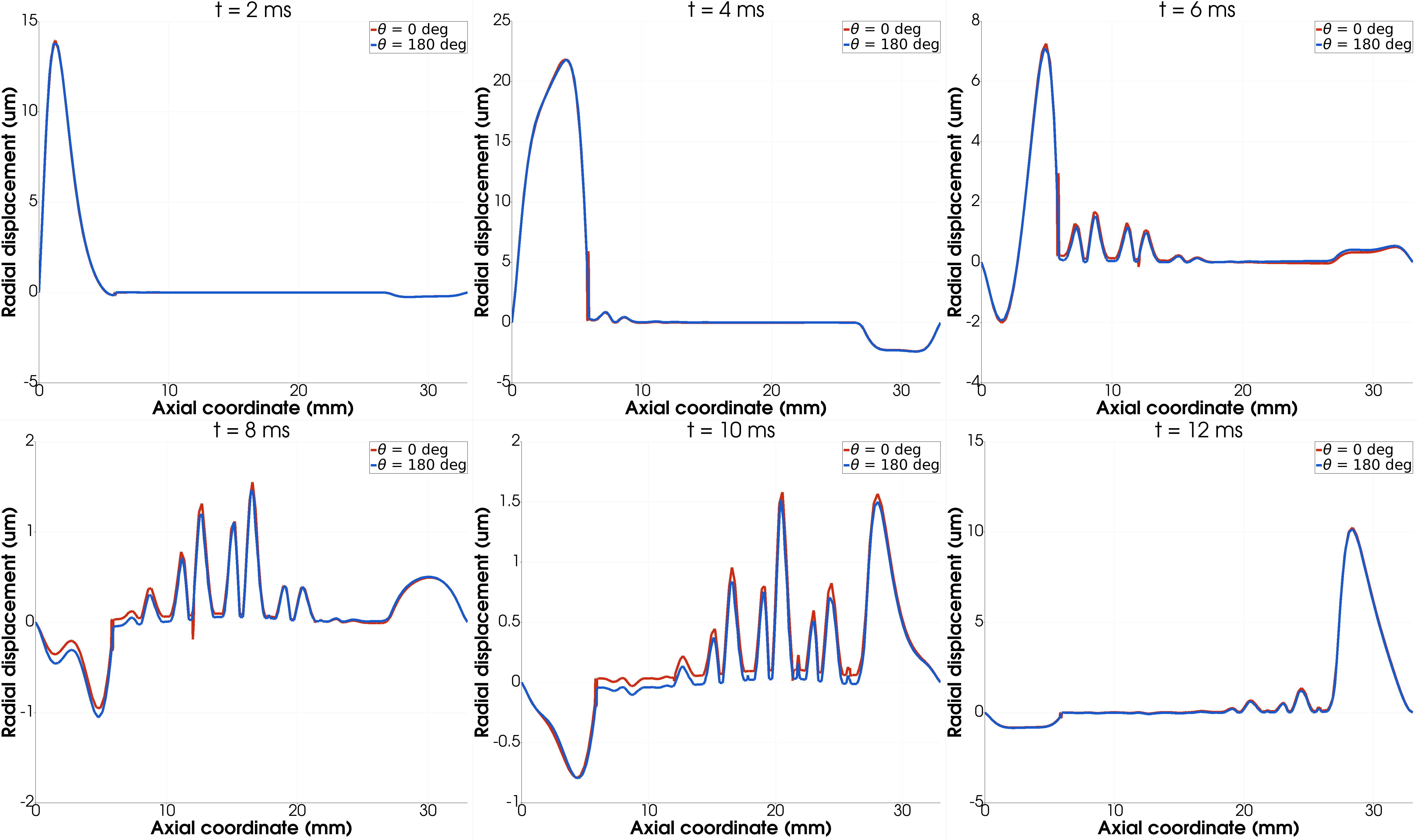}
\caption{Local radial displacement at the interface along the longitudinal
slicing $\theta=0,\pi$. The horizontal coordinate is distance from the
inlet. The two curves are evaluated on opposite sides of the artery.}
\label{fig:stented-longitudinal}
\end{figure}

The figures \ref{fig:stented-pressure}  shows pressure traveling along the center axis at six different times $2$, $4$, $6$, $8$, $10$, and
$12\,\mathrm{ms}$. Figure \ref{fig:stented_displacement} shows the corresponding stented artery
displacements. 
Figure~\ref{fig:stented-longitudinal} represents both sides of the
longitudinal plane $\theta=0,\pi$. In the central region, local
peaks and troughs occur along the interface, and thus the response is not spatially uniform. This is also why both axial and circumferential resolution are needed to describe this heterogeneous shell response.

\section{Conclusions}\label{sec:conclusions}

We have developed an explicit Robin-Robin splitting scheme for fluid interaction with a Naghdi shell. The coupling exchanges translational velocity and interface forces, while the shell formulation retains independent rotations and rotational inertia. Both subproblems use previous-time interface data, allowing them to be solved independently and in parallel without coupling subiterations.

For the fixed-domain Stokes-Naghdi problem, we derive a discrete energy identity that accounts for the mechanical energy, additional interface terms, and nonnegative contributions from velocity mismatch and force imbalance. We establish unconditional stability for homogeneous external data with any fixed positive common Robin coefficient, along with stability estimates for nonzero external data. We further prove first-order temporal consistency and an $O(\sqrt{\Delta t})$ error estimate in the mechanical energy norm under suitable regularity and initialization assumptions.

The straight-tube numerical experiments investigate the effects of the Robin parameter and time step on the computed pressure and displacement. At a fixed time step, increasing the Robin parameter reduces the discrepancies relative to the reference solutions, although the improvement becomes marginal for larger values of $\gamma$. Time refinement with $\gamma=1500$ and $3000$ further reduces the pressure and displacement discrepancies relative to the saved $\beta$ solution. These results demonstrate that the accuracy of the proposed scheme depends on both the Robin parameter and the time step.

The stented-artery example extends the numerical application to a bonded layered shell with a large stiffness contrast between the steel struts and the surrounding arterial wall. The pressure profiles show the pulse propagation, while the shell displacement exhibits spatial variations in both the axial and circumferential directions. These variations motivate the use of a Naghdi shell formulation that accounts for membrane, bending, and transverse-shear deformation. The independent shell rotations provide additional flexibility in representing localized bending near stiffness transitions. More importantly, the low-regularity Naghdi formulation avoids requiring second derivatives of the reference midsurface in its geometric coefficients, as discussed in~\cite{WangTambacaCanic2026}. 

More generally, the proposed Robin-Robin splitting scheme provides a framework for fluid interaction with thin structures represented by 2D surfaces embedded in a 3D fluid domain, rather than 3D elastic bodies. Since the coupling involves only the structural translational velocity and the corresponding interface force, the same splitting strategy can be extended to other reduced-dimensional structural models, including membrane, plate, and shell models. This allows the fluid and structural subproblems to be advanced independently without introducing a 3D structural domain. The stability and convergence analysis for such extensions would depend on the properties of the corresponding structural model and remains a subject for future investigation.

\section*{Data Availability}
The data that support the findings of this study are
available from the corresponding author upon reasonable request.

\newpage

\bibliographystyle{unsrtnat}
\bibliography{references}

@article{TambacaTutek2016,
  author = {Tamba{\v{c}}a, Josip and Tutek, Zvonimir},
  title = {A new linear {Naghdi} type shell model for shells with little regularity},
  journal = {Applied Mathematical Modelling},
  volume = {40},
  number = {23--24},
  pages = {10549--10562},
  year = {2016},
  doi = {10.1016/j.apm.2016.07.007}
}

@misc{WangTambacaCanic2026,
  author = {Wang, Yifan and Tamba{\v{c}}a, Josip and {\v{C}}ani{\'{c}}, Sun{\v{c}}ica},
  title = {A Coupled Fluid {Naghdi} Shell {Windkessel} Model for Pulmonary Artery Banding: Partitioned Numerical Scheme and Parameter Estimation},
  year = {2026},
  howpublished = {\href{https://arxiv.org/abs/2610.20966}{arXiv:2610.20966}},
  note = {Preprint},
  doi = {10.48550/arXiv.2610.20966}
}

@article{CausinGerbeauNobile2005,
  author = {Causin, Paola and Gerbeau, Jean-Fr{\'{e}}d{\'{e}}ric and Nobile, Fabio},
  title = {Added-mass effect in the design of partitioned algorithms for fluid--structure problems},
  journal = {Computer Methods in Applied Mechanics and Engineering},
  volume = {194},
  number = {42--44},
  pages = {4506--4527},
  year = {2005},
  doi = {10.1016/j.cma.2004.12.005}
}

@article{BukacEtAl2013,
  author = {Buka{\v{c}}, Martina and {\v{C}}ani{\'{c}}, Sun{\v{c}}ica and Glowinski, Roland and Tamba{\v{c}}a, Josip and Quaini, Annalisa},
  title = {Fluid--structure interaction in blood flow capturing non-zero longitudinal structure displacement},
  journal = {Journal of Computational Physics},
  volume = {235},
  pages = {515--541},
  year = {2013},
  doi = {10.1016/j.jcp.2012.08.033}
}

@article{BukacMuha2016,
  author = {Buka{\v{c}}, Martina and Muha, Boris},
  title = {Stability and Convergence Analysis of the Extensions of the Kinematically Coupled Scheme for the Fluid--Structure Interaction},
  journal = {SIAM Journal on Numerical Analysis},
  volume = {54},
  number = {5},
  pages = {3032--3061},
  year = {2016},
  doi = {10.1137/16M1055396}
}

@article{BurmanFernandez2009,
  author = {Burman, Erik and Fern{\'{a}}ndez, Miguel A.},
  title = {Stabilization of explicit coupling in fluid--structure interaction involving fluid incompressibility},
  journal = {Computer Methods in Applied Mechanics and Engineering},
  volume = {198},
  number = {5--8},
  pages = {766--784},
  year = {2009},
  doi = {10.1016/j.cma.2008.10.012}
}

@article{FernandezMullaertVidrascu2015,
  author = {Fern{\'{a}}ndez, Miguel A. and Mullaert, Jimmy and Vidrascu, Marina},
  title = {Generalized {Robin--Neumann} explicit coupling schemes for incompressible fluid-structure interaction: Stability analysis and numerics},
  journal = {International Journal for Numerical Methods in Engineering},
  volume = {101},
  number = {3},
  pages = {199--229},
  year = {2015},
  doi = {10.1002/nme.4785}
}

@article{BurmanDurstGuzman2022,
  author = {Burman, Erik and Durst, Rebecca and Guzm{\'{a}}n, Johnny},
  title = {Stability and error analysis of a splitting method using {Robin--Robin} coupling applied to a fluid--structure interaction problem},
  journal = {Numerical Methods for Partial Differential Equations},
  volume = {38},
  number = {5},
  pages = {1396--1406},
  year = {2022},
  doi = {10.1002/num.22840}
}

@article{BurmanEtAl2022,
  author = {Burman, Erik and Durst, Rebecca and Fern{\'{a}}ndez, Miguel A. and Guzm{\'{a}}n, Johnny},
  title = {Fully discrete loosely coupled {Robin--Robin} scheme for incompressible fluid--structure interaction: stability and error analysis},
  journal = {Numerische Mathematik},
  volume = {151},
  number = {4},
  pages = {807--840},
  year = {2022},
  doi = {10.1007/s00211-022-01295-y}
}

@misc{GuoEtAl2025RobinRobin,
  author = {Guo, Shihan and Lin, Ping and Wang, Yifan and Yue, Xiaohe and Zheng, Haibiao},
  title = {An Unconditionally Stable Explicit {Robin--Robin} Partitioned Scheme for Fluid--Structure Interaction},
  year = {2025},
  howpublished = {\href{https://arxiv.org/abs/2510.13096}{arXiv:2510.13096}},
  note = {Preprint, version 2 revised September 2026},
  doi = {10.48550/arXiv.2510.13096}
}

\end{document}